\documentclass[a4paper]{scrartcl}

\usepackage{amsmath,amsthm,amssymb,dsfont,extarrows,amscd}

\usepackage{tikz}
	\usetikzlibrary{shapes}
	\usetikzlibrary{matrix}
\usepackage{tikz-cd}
\usetikzlibrary{babel}
\usepackage{enumitem}
\usepackage[margin=1.2in]{geometry}
\usepackage{scalerel}
\usepackage{microtype} 

\usepackage[backend=biber, 
style=alphabetic,
giveninits=true,
isbn=false, 
url=false,
minalphanames=1,
maxalphanames=5,
maxbibnames=99]{biblatex}
\AtEveryBibitem{\clearfield{month}}
\DeclareNameAlias{author}{given-family}
\renewbibmacro{in:}{}
\DeclareFieldFormat{extraalpha}{#1}

\DeclareLabelalphaTemplate{
  \labelelement{
    \field[final]{shorthand}
    \field{label}
    \field[strwidth=3,strside=left,ifnames=1]{labelname}
    \field[strwidth=1,strside=left]{labelname}
  }
}

\usepackage{tabstackengine}
\TABstackMath
\usepackage{multirow}
\usepackage{graphicx}
\definecolor{Myblue}{rgb}{0,0,0.6}  
\usepackage[colorlinks,citecolor=Myblue,linkcolor=Myblue,urlcolor=Myblue,pdfpagemode=None, linktocpage]{hyperref}
\usepackage{todonotes}
	\setuptodonotes{size=\footnotesize, color=red!40}
	\newcommand{\todoNC}[1]{\todo[inline, color=orange!40]{Nils: #1}}
	\newcommand{\NC}[1]{\todo[color=orange!40]{N: #1}}
	\newcommand{\todoBH}[1]{\todo[inline, color=blue!40]{Benni: #1}}
	\newcommand{\BH}[1]{\todo[color=blue!40]{B: #1}}

\newcommand{\pic}[2][0.75]{
	\begin{tikzpicture}[scale=0.5,baseline={([yshift=-.5ex]current bounding box.center)}]
	\node at (0,0) {\includegraphics[scale=#1]{figures/#2}};
	\end{tikzpicture}
}

\newcommand{\B}{\mathcal{C}_1}
\usepackage{tikz-3dplot}
\usepackage{pgfplots}
\pgfplotsset{width=7cm,compat=1.8}
\usetikzlibrary{decorations.pathreplacing}
\usetikzlibrary{decorations.markings}
\usetikzlibrary{patterns}
\usepgflibrary{shapes.geometric}
\usepackage{subcaption}

\newcommand{\be}{\begin{equation}}
\newcommand{\ee}{\end{equation}}
\newcommand\tikzzbox[1]
{#1}

\tikzset{
	string/.style={draw=#1, postaction={decorate}, decoration={markings,mark=at position .51 with {\arrow[draw=#1]{>}}}},
	costring/.style={draw=#1, postaction={decorate}, decoration={markings,mark=at position .51 with {\arrow[draw=#1]{<}}}},
	ostring/.style={draw=#1, postaction={decorate}, decoration={markings,mark=at position .47 with {\arrow[draw=#1]{>}}}},
	ustring/.style={draw=#1, postaction={decorate}, decoration={markings,mark=at position .56 with {\arrow[draw=#1]{>}}}},
	oostring/.style={draw=#1, postaction={decorate}, decoration={markings,mark=at position .43 with {\arrow[draw=#1]{>}}}},
	uustring/.style={draw=#1, postaction={decorate}, decoration={markings,mark=at position .59 with {\arrow[draw=#1]{>}}}},
	directed/.style={string=blue!50!black}, 
	odirected/.style={ostring=blue!50!black}, 
	udirected/.style={ustring=blue!50!black}, 
	oodirected/.style={oostring=blue!50!black}, 
	uudirected/.style={uustring=blue!50!black},     
	redirected/.style={costring= blue!50!black},
	redirectedgreen/.style={costring= green!50!black},
	directedgreen/.style={string= green!50!black},
	redirectedred/.style={costring= red!50!black},
	directedred/.style={string= red!50!black}%
}

\tikzset{-dot-/.style={decoration={
			markings,
			mark=at position 0.5 with {\fill circle (2pt);}},postaction={decorate}}}

\tikzset{
	Fdot/.style={circle, draw, fill, inner sep=0pt}, 
	Odot/.style={circle, draw, inner sep=0.1pt, minimum size=0.1cm}
}

\def\nicecolourscheme{\shadedraw[top color=orange!29, bottom color=orange!29, draw=orange!22, draw=white, draw opacity=0]}
\def\nicepalecolourscheme{\shadedraw[top color=orange!12, bottom color=orange!12, draw=white, draw opacity=0]}

\newcommand{\chiXYinv}[5]{
	\draw[redirectedred] (#2) .. controls +(0,1) and +(0,1) .. ($(#1)+(1.5,0)$);
	\draw[redirectedred] ($(#3)+(.5,0)$) .. controls +(0,-1) and +(0,-1) .. (#4)  node [pos=0.35](psi){} node [pos=0.3, below] {\tiny$\psi^2$};
	\draw[redirectedred] (#1) .. controls +(0,1.5) and +(0,-1.5) .. (#3);
	\draw[redirectedred] ($(#5)+(1,0)$) .. controls +(0,-1.5) and +(0,1.5) .. ($(#2)+(1,0)$); 
	\draw[color=blue!50!black] ($(#4)+(.5,0)$) .. controls +(0,-1.5) and +(0,1.5) .. ($(#1)+(.5,0)$);
	\draw[redirectedred] ($(#4)+(1,0)$) .. controls +(0,-1.5) and +(0,1.5) .. ($(#2)+(.5,0)$);
	\draw[color=white, line width=4pt] ($(#5)+(.5,0)$) .. controls +(0,-1.5) and +(0,1.5) .. ($(#1)+(1,0)$);
	\draw[color=blue!50!black] ($(#5)+(.5,0)$) .. controls +(0,-1.5) and +(0,1.5) .. ($(#1)+(1,0)$);
	\draw[color=white, line width=4pt] ($(#3)+(1,0)$) .. controls +(0,-1) and +(0,-1) .. (#5);
	\draw[postaction={decorate}, decoration={markings,mark=at position .61 with {\arrow[color=red!50!black]{<}}}] ($(#3)+(1,0)$) .. controls +(0,-1) and +(0,-1) .. (#5)  node [pos=0.8](psi1){} node [pos=0.75, above] {\tiny$\psi^2$};
	\fill (psi) circle (2.9pt);
	\fill (psi1) circle (2.9pt);}
	
\newcommand{\omegaghk}[6]{
		\draw[directedred] (#4) .. controls +(0,-1.5) and +(0,1.5) .. (#1);
		\draw[directedred] ($(#3)+(.5,0)$) .. controls +(0,1.5) and +(0,-1.5) .. ($(#6)+(1,0)$); 
		\draw[directedred] ($(#1)+(.5,0)$) .. controls +(0,1) and +(0,1) .. (#2);
		\draw[postaction={decorate}, decoration={markings,mark=at position .60 with {\arrow[color=red!50!black]{>}}}] (#5) .. controls +(0,-1) and +(0,-1) .. ($(#4)+(.5,0)$)  node [pos=0.75](psi2){} node [pos=0.7, below] {\tiny$\psi^2$};
		\fill (psi2) circle (2.9pt);
		\draw[directedred] ($(#2)+(.5,0)$) .. controls +(0,1.5) and +(0,-1.5) .. ($(#5)+(.5,0)$);
		\draw[directedred] ($(#2)+(1,0)$) .. controls +(0,1.5) and +(0,-1.5) ..  ($(#6)+(.5,0)$);
		\draw[color=white, line width=4pt] ($(#1)+(1,0)$) .. controls +(0,1) and +(0,1) .. (#3);
		\draw[postaction={decorate}, decoration={markings,mark=at position .41 with {\arrow[color=red!50!black]{>}}}] ($(#1)+(1,0)$) .. controls +(0,1) and +(0,1) .. (#3);
		\draw[color=white, line width=4pt] (#6) .. controls +(0,-1) and +(0,-1) .. ($(#4)+(1,0)$);
		\draw[postaction={decorate}, decoration={markings,mark=at position .6 with {\arrow[color=red!50!black]{>}}}] (#6) .. controls +(0,-1) and +(0,-1) .. ($(#4)+(1,0)$)  node [pos=0.25](psi1){} node [pos=0.2, below] {\tiny$\psi^2$};
		\fill (psi1) circle (2.9pt);}
		
\newcommand{\boxtimestwoa}[8]{
		\draw[directedred] ($(#3)+(.5,0)$) .. controls +(0,1.5) and +(0,-1.5) .. ($(#6)+(.5,0)$); 
		\draw[color=blue!50!black] ($(#1)+(.5,0)$) .. controls +(0,1.5) and +(0,-1.5) .. ($(#5)+(.5,0)$); 
		\draw[redirectedred] (#1) .. controls +(0,1.5) and +(0,-1.5) .. (#5); 
		\draw[directedred] ($(#4)+(.5,0)$) .. controls +(0,1.5) and +(0,-1.5) .. ($(#8)+(.5,0)$); 
		\draw[directedred] ($(#4)+(1,0)$) .. controls +(0,1.5) and +(0,-1.5) .. ($(#8)+(1,0)$);
		\draw[directedred] ($(#1)+(1,0)$) .. controls +(0,1) and +(0,1) .. (#3);
		\draw[color=white, line width=4pt] ($(#2)+(1,0)$) .. controls +(0,1) and +(0,1) .. (#4);
		\draw[directedred] ($(#2)+(1,0)$) .. controls +(0,1) and +(0,1) .. (#4);
	\draw[directedred] (#8) .. controls +(0,-1) and +(0,-1) .. ($(#7)+(1,0)$)  node [pos=0.65](psi1){} node [pos=0.6, below] {\tiny$\psi^2$};
	\draw[directedred] (#6) .. controls +(0,-1) and +(0,-1) .. ($(#5)+(1,0)$)  node [pos=0.65](psi2){} node [pos=0.6, below] {\tiny$\psi^2$}; 
		\draw[color=white, line width=4pt] (#2) .. controls +(0,1.5) and +(0,-1.5) .. (#7);
		\draw[redirectedred] (#2) .. controls +(0,1.5) and +(0,-1.5) .. (#7); 
		\draw[color=white, line width=4pt] ($(#2)+(.5,0)$) .. controls +(0,1.5) and +(0,-1.5) .. ($(#7)+(.5,0)$); 
		\draw[color=blue!50!black] ($(#2)+(.5,0)$) .. controls +(0,1.5) and +(0,-1.5) .. ($(#7)+(.5,0)$); 
		\fill (psi1) circle (2.9pt);
		\fill (psi2) circle (2.9pt);
}

\newcommand\boxtimestwob[1]{%
    \def\tempa{#1}%
    \boxtimestwobcontinued
}
\newcommand{\boxtimestwobcontinued}[9]{
		\draw[redirectedred] (#1) .. controls +(0,1.5) and +(0,-1.5) .. (#6); 
		\draw[color=blue!50!black] ($(#2)+(.5,0)$) .. controls +(0,1.5) and +(0,-1.5) .. ($(#8)+(.5,0)$);
		\draw[color=blue!50!black] ($(#4)+(.5,0)$) .. controls +(0,1.5) and +(0,-1.5) .. ($(#9)+(.5,0)$);
		\draw[color=blue!50!black] ($(#5)+(.5,0)$) .. controls +(0,1.5) and +(0,-1.5) .. ($(\tempa)+(.5,0)$);
		\draw[directedred] ($(#2)+(1,0)$) .. controls +(0,1.5) and +(0,-1.5) .. ($(#8)+(1,0)$); 
		\draw[directedred] ($(#4)+(1,0)$) .. controls +(0,1.5) and +(0,-1.5) .. ($(#9)+(1,0)$); 
		\draw[directedred] ($(#5)+(1,0)$) .. controls +(0,1.5) and +(0,-1.5) .. ($(\tempa)+(1,0)$); 
		\draw[directedred] ($(#1)+(.5,0)$) .. controls +(0,1) and +(0,1) .. (#2);
		\draw[directedred] ($(#3)+(.5,0)$) .. controls +(0,1) and +(0,1) .. (#4);
	\draw[directedred] (#8) .. controls +(0,-1) and +(0,-1) .. ($(#6)+(.5,0)$)  node [pos=0.75](psi1){} node [pos=0.75, below] {\tiny$\psi^2$}; 
		\draw[color=white, line width=4pt] (#3) .. controls +(0,1.5) and +(0,-1.5) .. (#7);
		\draw[redirectedred] (#3) .. controls +(0,1.5) and +(0,-1.5) .. (#7); 
		\draw[color=white, line width=4pt] ($(#3)+(1,0)$) .. controls +(0,1) and +(0,1) .. (#5);
		\draw[directedred] ($(#3)+(1,0)$) .. controls +(0,1) and +(0,1) .. (#5);
	\draw[color=white, line width=4pt] (#9) .. controls +(0,-1) and +(0,-1) .. ($(#7)+(.5,0)$); 
	\draw[directedred] (#9) .. controls +(0,-1) and +(0,-1) .. ($(#7)+(.5,0)$)  node [pos=0.75](psi2){} node [pos=0.7, below] {\tiny$\psi^2$}; 
	\draw[color=white, line width=4pt] (\tempa) .. controls +(0,-1) and +(0,-1) .. ($(#7)+(1,0)$); 
	\draw[directedred] (\tempa) .. controls +(0,-1) and +(0,-1) .. ($(#7)+(1,0)$)  node [pos=0.2](psi3){} node [pos=0.2, below] {\tiny$\psi^2$}; 
		\fill (psi1) circle (2.9pt);
		\fill (psi2) circle (2.9pt);
		\fill (psi3) circle (2.9pt);
}
\newcommand{\boxtimestwoc}[8]{
		\draw[directedred] ($(#3)+(.5,0)$) .. controls +(0,1.5) and +(0,-1.5) .. ($(#6)+(.5,0)$); 
		\draw[directedred] ($(#3)+(1,0)$) .. controls +(0,1.5) and +(0,-1.5) .. ($(#6)+(1,0)$); 
		\draw[color=blue!50!black] ($(#1)+(.5,0)$) .. controls +(0,1.5) and +(0,-1.5) .. ($(#5)+(.5,0)$); 
		\draw[redirectedred] (#1) .. controls +(0,1.5) and +(0,-1.5) .. (#5); 
		\draw[directedred] ($(#4)+(.5,0)$) .. controls +(0,1.5) and +(0,-1.5) .. ($(#8)+(.5,0)$); 
		\draw[directedred] ($(#1)+(1,0)$) .. controls +(0,1) and +(0,1) .. (#3);
		\draw[color=white, line width=4pt] ($(#2)+(1,0)$) .. controls +(0,1) and +(0,1) .. (#4);
		\draw[directedred] ($(#2)+(1,0)$) .. controls +(0,1) and +(0,1) .. (#4);
	\draw[directedred] (#8) .. controls +(0,-1) and +(0,-1) .. ($(#7)+(1,0)$)  node [pos=0.65](psi1){} node [pos=0.6, below] {\tiny$\psi^2$};
	\draw[directedred] (#6) .. controls +(0,-1) and +(0,-1) .. ($(#5)+(1,0)$)  node [pos=0.65](psi2){} node [pos=0.6, below] {\tiny$\psi^2$}; 
		\draw[color=white, line width=4pt] (#2) .. controls +(0,1.5) and +(0,-1.5) .. (#7);
		\draw[redirectedred] (#2) .. controls +(0,1.5) and +(0,-1.5) .. (#7); 
		\draw[color=white, line width=4pt] ($(#2)+(.5,0)$) .. controls +(0,1.5) and +(0,-1.5) .. ($(#7)+(.5,0)$); 
		\draw[color=blue!50!black] ($(#2)+(.5,0)$) .. controls +(0,1.5) and +(0,-1.5) .. ($(#7)+(.5,0)$); 
		\fill (psi1) circle (2.9pt);
		\fill (psi2) circle (2.9pt);
}

\newcommand\boxtimestwod[1]{%
    \def\tempa{#1}%
    \boxtimestwodcontinued
}
\newcommand{\boxtimestwodcontinued}[9]{
		\draw[redirectedred] (#1) .. controls +(0,1.5) and +(0,-1.5) .. (#6); 
		\draw[color=blue!50!black] ($(#2)+(.5,0)$) .. controls +(0,1.5) and +(0,-1.5) .. ($(#8)+(.5,0)$);
		\draw[color=blue!50!black] ($(#3)+(.5,0)$) .. controls +(0,1.5) and +(0,-1.5) .. ($(#9)+(.5,0)$);
		\draw[color=blue!50!black] ($(#5)+(.5,0)$) .. controls +(0,1.5) and +(0,-1.5) .. ($(\tempa)+(.5,0)$);
		\draw[directedred] ($(#2)+(1,0)$) .. controls +(0,1.5) and +(0,-1.5) .. ($(#8)+(1,0)$); 
		\draw[directedred] ($(#3)+(1,0)$) .. controls +(0,1.5) and +(0,-1.5) .. ($(#9)+(1,0)$); 
		\draw[directedred] ($(#5)+(1,0)$) .. controls +(0,1.5) and +(0,-1.5) .. ($(\tempa)+(1,0)$); 
		\draw[directedred] ($(#1)+(.5,0)$) .. controls +(0,1) and +(0,1) .. (#2);
		\draw[directedred] ($(#4)+(.5,0)$) .. controls +(0,1) and +(0,1) .. (#5);
		\draw[color=white, line width=4pt] ($(#1)+(1,0)$) .. controls +(0,1) and +(0,1) .. (#3);
		\draw[directedred] ($(#1)+(1,0)$) .. controls +(0,1) and +(0,1) .. (#3);
		\draw[directedred] (#8) .. controls +(0,-1) and +(0,-1) .. ($(#6)+(.5,0)$)  node [pos=0.8](psi1){} node [pos=0.8, below] {\tiny$\psi^2$}; 
		\draw[color=white, line width=4pt] (#9) .. controls +(0,-1) and +(0,-1) .. ($(#6)+(1,0)$); 
		\draw[directedred] (#9) .. controls +(0,-1) and +(0,-1) .. ($(#6)+(1,0)$)  node [pos=0.2](psi2){} node [pos=0.2, above] {\tiny$\psi^2$}; 
		\draw[color=white, line width=4pt] (\tempa) .. controls +(0,-1.5) and +(0,-1.5) .. ($(#7)+(.5,0)$); 
		\draw[directedred] (\tempa) .. controls +(0,-1.5) and +(0,-1.5) .. ($(#7)+(.5,0)$)  node [pos=0.25](psi3){} node [pos=0.2, above] {\tiny$\psi^2$}; 
		\draw[color=white, line width=4pt] (#4) .. controls +(0,1.5) and +(0,-1.5) .. (#7);
		\draw[redirectedred] (#4) .. controls +(0,1.5) and +(0,-1.5) .. (#7); 
		\fill (psi1) circle (2.9pt);
		\fill (psi2) circle (2.9pt);
		\fill (psi3) circle (2.9pt);
}

\newcommand{\unitI}[4]{
		\draw[redirectedred] (#1) .. controls +(0,1) and +(0,-1) .. (#3); 
		\draw[directedred] ($(#2)+(.5,0)$) .. controls +(0,1.5) and +(0,-1.5) .. ($(#4)+(1,0)$);  
		\draw[color=blue!50!black] ($(#1)+(.5,0)$) .. controls +(0,1.5) and +(0,-1.5) .. ($(#4)+(.5,0)$);
		\draw[directedred] (#4) .. controls +(0,-1) and +(0,-1) .. ($(#3)+(.5,0)$)  node [pos=0.65](psi1){} node [pos=0.55, above] {\tiny$\psi^2$}; 
		\draw[directedred] ($(#1)+(1,0)$) .. controls +(0,1) and +(0,1) .. (#2);
		\fill (psi1) circle (2.9pt);
}

\makeatletter
\newif\iftikztransformnodecoordinates
\tikzset{transform node coordinates/.is if=tikztransformnodecoordinates}

\def\tikz@parse@node#1(#2){%
    \pgfutil@in@.{#2}%
    \ifpgfutil@in@
        \tikz@calc@anchor#2\tikz@stop%
    \else%
        \tikz@calc@anchor#2.center\tikz@stop%
        \expandafter\ifx\csname pgf@sh@ns@#2\endcsname\tikz@coordinate@text%
        \else
            \tikz@shapebordertrue%
            \def\tikz@shapeborder@name{#2}%
        \fi%
    \fi%
    \iftikztransformnodecoordinates%
        \pgf@pos@transform{\pgf@x}{\pgf@y}%
    \fi
    \edef\tikz@marshal{\noexpand#1{\noexpand\pgfqpoint{\the\pgf@x}{\the\pgf@y}}}%
    \tikz@marshal%
}

\usepackage{thmtools}
\usepackage{mathtools}
\usepackage{cleveref}

\allowdisplaybreaks

\deffootnote[1em]{1em}{1em}{\textsuperscript{\thefootnotemark}}

\theoremstyle{definition} 
\newtheorem{definition}{Definition}
\newtheorem{proposition}[definition]{Proposition}
\newtheorem{theorem}[definition]{Theorem}
\newtheorem{lemma}[definition]{Lemma}
\newtheorem{corollary}[definition]{Corollary}
\newtheorem{remark}[definition]{Remark}

\newtheorem{example}[definition]{Example}

\numberwithin{equation}{section}
\numberwithin{definition}{section}
\numberwithin{lemma}{section}
\numberwithin{proposition}{section}
\numberwithin{theorem}{section}
\numberwithin{corollary}{section}
\numberwithin{example}{section}
\numberwithin{remark}{section}
\numberwithin{figure}{section}
\numberwithin{table}{section}
\renewcommand{\labelenumi}{(\roman{enumi})}

\DeclareMathOperator{\pr}{pr}
\DeclareMathOperator{\Tr}{Tr}
\DeclareMathOperator{\tr}{tr}
\DeclareMathOperator{\Fun}{Fun}

\DeclareMathOperator{\TQFT}{\mathcal{Z}}

\newcommand*{\longhookrightarrow}{\ensuremath{\lhook\joinrel\relbar\joinrel\rightarrow}}
\newcommand*{\twoheadlongrightarrow}{\ensuremath{\relbar\joinrel\twoheadrightarrow}}
\newcommand{\longdoubleleftrightarrow}{\mathrel{\substack{\xrightarrow{\rule{0.3cm}{0cm}} \\[-.9ex] \xleftarrow{\rule{0.3cm}{0cm}}}}}
\newcommand*{\longhookleftarrow}{\ensuremath{\leftarrow\joinrel\relbar\joinrel\rhook}}
\newcommand*{\LongRRightarrow}{\ensuremath{\equiv\joinrel\Rrightarrow}}

\newcommand{\R}{\mathds{R}}
\newcommand{\C}{\mathds{C}}
\newcommand{\N}{\mathds{N}}
\newcommand{\Z}{\mathds{Z}}
\newcommand{\K}{\mathds{K}}
\newcommand{\Hilb}{\mathcal{H}} 
\newcommand{\Vect}{\mathrm{Vect}} 
\newcommand{\Grp}{\mathrm{Grp}}
\newcommand{\cFA}[1][\K]{\mathrm{cFrAlg}_{#1}}
\newcommand{\nTQFT}[2][n]{#1 \mathrm{TQFT}_{#2}} 
\newcommand{\BG}[1][0]{\ifnum #1<1 {\mathrm{B}G} \else {\mathrm{B}^{#1}G}\fi}
\newcommand{\EG}{\mathrm{EG}}
\newcommand{\Mod}{\mathrm{Mod}}
\DeclareMathOperator{\Rep}{Rep}
\newcommand{\Gcbc}{\mathcal{C}^\times_G}
\newcommand{\CAR}{\mathcal{C}_{\mathbb{A}_R}}
\let\to\undefined
\newcommand{\to}{\longrightarrow}
\let\mapsto\undefined
\newcommand{\mapsto}{\longmapsto}

\newcommand{\Homorb}[1]{{#1}_{\textrm{Hom-orb}}}
\usepackage{bbm}

\DeclareMathOperator{\codim}{codim}
\newcommand{\id}{\mathrm{id}}
\DeclareMathOperator{\Aut}{Aut}
\DeclareMathOperator{\End}{End}
\DeclareMathOperator{\Hom}{Hom}
\newcommand{\ev}{\mathrm{ev}}
\newcommand{\coev}{\mathrm{coev}}
\newcommand{\tril}{\vartriangleleft}
\newcommand{\trir}{\vartriangleright}

\title{
2-Group Symmetries of 3-dimensional Defect TQFTs and Their Gauging}

\author{Nils Carqueville$^*$ \quad Benjamin Haake$^\#$
	\\[0.5cm]
	\normalsize{\texttt{\href{mailto:nils.carqueville@univie.ac.at}{nils.carqueville@univie.ac.at}}} 
	\quad  
	\normalsize{\texttt{\href{mailto:B.Haake@sms.ed.ac.uk}{B.Haake@sms.ed.ac.uk}}}
	\\[0.5cm]  %
	{\normalsize\slshape $^*$Universit\"at Wien, Fakult\"at f\"ur Physik, Boltzmanngasse 5, 1090 Wien, \"{O}sterreich}\\
	{\normalsize\slshape $^\#$University of Edinburgh, School of Mathematics, James Clerk Maxwell Building,\vspace{-.5em}}\\{\normalsize\slshape Peter Guthrie Tait Road, Edinburgh EH9\,3FD, United Kingdom}\\
	{\normalsize\slshape $^\#$Maxwell Institute for Mathematical Sciences, The Bayes Centre,\vspace{-.5em}}\\{\normalsize\slshape 47 Potterrow, Edinburgh EH8 9BT, United Kingdom\vspace{-1em}}
}
\date{}

\begin{document}

\maketitle


\begin{abstract}
	{\small A large class of symmetries of topological quantum field theories is naturally described by functors into higher categories of topological defects. 
	Here we study 2-group symmetries of 3-dimensional TQFTs. 
	We explain that these symmetries can be gauged to produce new TQFTs iff certain defects satisfy the axioms of orbifold data.
	In the special case of Reshetikhin--Turaev theories coming from $G$-crossed braided fusion categories~$\mathcal C^\times_G$, we show that there are 0- and 1-form symmetries which have no obstructions to gauging. 
	We prove that gauging the 0-form $G$-symmetry on the neutral component~$\mathcal C_e$ of~$\mathcal C^\times_G$ produces its equivariantisation~$(\mathcal C^\times_G)^G$, which in turn features a generalised symmetry whose gauging recovers~$\mathcal C_e$.
	If~$G$ is commutative, the latter symmetry reduces to a 1-form symmetry involving the Pontryagin dual group.}
\end{abstract}

\newpage

\tableofcontents


\section{Introduction and Summary}
\label{sec:Introduction}

A standard example of a symmetry is a representation of a group~$G$ on a vector space~$V$, i.e.\ a group homomorphism $G\longrightarrow \textrm{Aut}(V)$. 
Equivalently, this is a functor 
\begin{align}
	\textrm{B}G & \longrightarrow \Vect
	\nonumber
	\\ 
	* & \longmapsto V
\end{align} 
from the delooping $\textrm{B}G$ (the category which has only a single object~$*$ whose endomorphisms are given by the group~$G$) to the category of vector spaces. 

In the context of quantum field theory, $V$ may be a state space on which the (global) symmetry group~$G$ acts as above. 
However, such a $G$-representation is an induced structure. It is natural to consider, at least hypothetically, an $n$-category~$\mathcal D$ whose objects are $n$-dimensional QFTs, and whose $k$-morphisms are \textsl{topological} defects of codimension~$k$. 
A \textsl{$G$-symmetry of a QFT~$Q$} is then an $n$-functor 
\begin{align}
	R\colon 
	\textrm{B}\underline{G} & \longrightarrow \mathcal D 
	\nonumber
	\\
	* & \longmapsto Q
	\label{eq:Gsym1}
\end{align} 
from the delooping of~$\underline{G}$ to~$\mathcal D$, where~$\underline{G}$ is~$G$ viewed as a monoidal $(n-1)$-category with only identity morphisms. 
Such an $n$-functor sends group elements~$g$ to defects $R(g)$ of codimension~1, whose fusion is described by a family of weakly invertible defects $R_{g,h}\colon R(g)R(h) \longrightarrow R(gh)$ of codimension~2, which in turn compose up to defects of codimension~3, and so on. 
All this is contained in the structure of~$R$. 
Returning to the introductory example, if we write $V_Q := \mathrm{End}_{\mathcal D} (\id_{{\id_{\dots}}_{\id_{Q}}})$ for the local state space of~$Q$, then we obtain an ordinary $G$-representation $G\longrightarrow \textrm{Aut}(V_Q)$ by wrapping the defects $R(g)$ as $(n-1)$-spheres around point insertions $\varphi \in V_Q$. 
In dimension $n=2$ we can use the graphical calculus in the pivotal 2-category~$\mathcal D$ to write this $G$-action as
\begin{equation} 
\label{eq:Gaction via defects}
	V_Q \;\ni \; \varphi \longmapsto R(g)(\varphi) \;=\;
	\begin{tikzpicture}[very thick,scale=0.4,color=blue!50!black, baseline]
		\nicepalecolourscheme (-2.8,-1.7) rectangle (2.8,1.7);
		\nicecolourscheme (0,0) circle (1.25);
		\fill (0,0) circle (3pt) node {};
		\fill (0.5,0) circle (0pt) node {{\scriptsize$\varphi$}};
		\draw (0,0) circle (1.25);
		\fill (-45:1.2) circle (0pt) node[right] {{\scriptsize$R(g)$}};
		\draw[<-,  thick] (0.100,-1.25) -- (-0.101,-1.25) node[above] {}; 
		\draw[<-,  thick] (-0.100,1.25) -- (0.101,1.25) node[below] {}; 
	\end{tikzpicture} 
	. 
\end{equation}

The notion of a symmetry of a QFT is naturally generalised by replacing~$\underline{G}$ in the domain of the functor~$R$ in~\eqref{eq:Gsym1} with a representative of a broader class of monoidal $(n-1)$-categories. 
For example, a \textsl{$p$-form symmetry} of~$Q$ is an $n$-functor $\textrm{B}^{p+1}\underline{H}\to \mathcal D$ with $*\longmapsto Q$, where $p\in\{0,1,\dots,n-2\}$, and~$H$ is a group (which is necessarily commutative if $p>0$ for $\textrm{B}^{p+1}\underline{H}$ to be well-defined). 
Such an $n$-functor sends group elements to defects of codimension $p+1$, while otherwise the discussion is analogous to the case of 0-form symmetries as in~\eqref{eq:Gsym1}. 

A slightly broader class of examples of symmetries is given by $n$-functors $\textrm{B}\underline{\mathcal G} \to \mathcal D$, where~$\mathcal G$ is an $(n-1)$-group, i.e.\ a monoidal $(n-2)$-category in which all objects and all morphisms are (weakly) invertible. 
One may think of this as gluing $p$-form symmetries together for all $p\leqslant n-2$, where for $n\geqslant 3$ there is more than one way of doing so. 
This is made explicit e.g.\ by the equivalence of 2-groups and crossed modules (i.e.\ two groups that ``interact'' via two group homomorphisms), as we recall in \Cref{sec:2-groups} below. 

The above symmetries are ``invertible'' in the sense that their domains $\textrm{B}\underline{G}$, $\textrm{B}^{p+1}\underline{H}$ and $\textrm{B}\underline{\mathcal G}$ contain only invertible (higher) morphisms. 
If the domain has non-invertible morphisms, one speaks of \textsl{non-invertible symmetries}. 
For example, symmetries of 2-dimensional rational conformal field theories~$Q$ may come as 2-functors $\textrm{B}\mathcal S \to \mathcal D$ with $*\longmapsto Q$, where~$\mathcal S$ is a spherical fusion category, and analogously for higher fusion categories. 
There is however no general need to restrict to semisimple domains. 

\medskip 

Given a (global) symmetry $R\colon \textrm{B}\mathcal A \to \mathcal D$ of a QFT $Q\in\mathcal D$ for some monoidal $(n-1)$-category~$\mathcal A$, a natural question is whether~$R$ can be gauged to produce a new QFT, or whether there is an obstruction to such a procedure. 
In the context of oriented topological quantum field theories, this is precisely the question of whether~$R$ gives rise to an ``orbifold datum''~$\mathbb{A}_R$. 
We briefly outline the relevant theory below, for details we refer to Sections~\ref{sec:Background} and~\ref{subsec:TQFTPerspective} as well as \cite{Carqueville2017} or the survey \cite{carqueville2023orbifoldstopologicalquantumfield}. 

The starting point is a defect TQFT by which we mean a symmetric monoidal functor~$\mathcal Z$ from a category of stratified and decorated bordisms to $\Vect$. 
The labels for $n$-dimensional strata can be identified with bulk theories (like~$Q$ in the above discussion) while the lower-dimensional decorated strata are topological defects. 
Given one label~$A_j$ for $j$-dimensional defects for all $j\geqslant 1$ and two labels $A_0^+, A_0^-$ for point defects (having to do with opposite orientations), one may try to ``gauge'' such defect data to produce a new closed TQFT~$\mathcal Z_{\mathbb{A}}$ as follows: 
for every bordism, choose an oriented triangulation, label its Poincar\'e dual stratification with $A_n,\dots,A_1,A_0^+,A_0^-$, evaluate the thus-obtained defect bordism with~$\mathcal Z$, and take the colimit over all possible triangulations. 
The last step can be thought of as an attempt to consistently average over all choices. 

For a generic list of defects $\mathbb{A} := (A_n,\dots,A_1,A_0^+,A_0^-)$ the above procedure will not work as no consistent averaging is possible. 
Then~$\mathbb{A}$ cannot be gauged. 
The obstruction comes precisely from the non-invariance of the construction under the choice of triangulations. 
On the other hand, since any two triangulations of a compact manifold can be obtained from one another by a finite number of bistellar (or: Pachner) moves of which there are only finitely many, we can deduce a finite number of constraint equations on~$\mathbb{A}$. 
If they are satisfied, there is no obstruction to gauging in the above sense, and we call~$\mathbb{A}$ an \textsl{orbifold datum} for~$\mathcal Z$, and the gauged theory~$\mathcal Z_{\mathbb{A}}$ is the $\mathbb{A}$-orbifold of~$\mathcal Z$. 

For a wide class of symmetries $R\colon \textrm{B}\mathcal A \to \mathcal D$ of a closed TQFT~$Q$ one can construct associated \textsl{candidate} orbifold data~$\mathbb{A}_R$. 
For example, if $\mathcal A = \underline{G}$ for a finite group~$G$ and $n=2$, this is described in detail in \cite{Brunner_2015}. 
Then whether or not the symmetry can be gauged is the same as the question of whether~$\mathbb{A}_R$ is an actual orbifold datum. 
This question can be phrased and answered in entirely algebraic language, namely internal to the $n$-category $\mathcal D$ that one naturally extracts from the ambient TQFT -- at least for low dimensions, see \cite{Davydov2011, Carqueville2016}. 
Then we also refer to~$\mathbb{A}_R$ as an \textsl{orbifold datum in~$\mathcal D$}. 

\medskip 

In the present paper we carry out the above programme for 0- and 1-form symmetries as well as 2-group symmetries for general 3-dimensional defect TQFTs. 
Moreover, in the special case of symmetries of Reshetikhin--Turaev theories coming from $G$-crossed braided fusion categories, we show that there is no obstruction to gauging, and we make precise contact to the algebraic notions of equivariantisation and de-equivariantisation. 

In slightly more detail, we first consider symmetries $R\colon \textrm{B}\underline{\mathcal G} \to \mathcal D$ where~$\mathcal G$ is a 2-group and~$\mathcal D$ is a $\Vect$-enriched Gray category with duals and a mild completeness property for its Hom 2-categories. 
We refer to \cite{Barrett2012} for Gray categories with duals in general, and note that~$\mathcal D$ can e.g.\ be taken to be the 3-category associated to a 3-dimensional defect TQFT in \cite{Carqueville2016}. 
From the data of~$R$ we then construct a candidate orbifold datum~$\mathbb{A}_R$ and prove (in \Cref{prop:2group orbifold datum}) that it is always a weakly associative algebra internal to~$\mathcal D$. 
This symmetry cannot be gauged in general (but note the discussion in \Cref{rem:HigherGaugableGroupSymmetries}), and neither can the special case of 0-form symmetries $\mathcal G = \underline{G}$ for some finite group~$G$, as we discuss in \Cref{sec:0-form specialisation}. 

The special case of a 1-form symmetry $R\colon \textrm{B}^2\underline{H}\to \mathcal D$ for a finite commutative group~$H$ is analysed in \Cref{sec:1-form specialisation}. 
Such a 3-functor is equivalently described by a braided monoidal 1-functor $\overline{R}\colon \underline{H} \to \textrm{End}_{\mathcal D}(R(\id_*))$. 
We prove (in \Cref{prop:Orbifold datum from 1form symmetry as composite}) that if~$\overline{R}$ is pivotal, then~$\mathbb{A}_R$ is an orbifold datum. 
These are precisely the gaugeable 1-form symmetries. 

\medskip 

After the general discussion in \Cref{sec:OrbDatFrom2Groups}, the remainder of the paper covers specific examples of 0- and 1-form symmetries in Reshetikhin--Turaev TQFTs coming from $G$-crossed braided fusion categories~$\mathcal C^\times_G$, where we find that they are always gaugeable. 
Recall that for a finite group~$G$, a $G$-crossed braided fusion category $\mathcal C^\times_G = \bigoplus_{g\in G} \mathcal C_g$ is a $G$-graded fusion category together with a $G$-action $\rho\colon \underline{G}^{\textrm{rev}} \to \textrm{Aut}^\otimes(\mathcal C^\times_G)$ and a ``crossed'' braiding with components $X\otimes Y \to Y\otimes \rho_h(X)$ for $X\in\mathcal C_h$. 
Assuming that the neutral component $\mathcal C_e \subset \mathcal C^\times_G$ is modular, one obtains a second modular fusion category, namely the equivariantisation $(\mathcal C^\times_G)^G$ with respect to the $G$-action~$\rho$. 
We review the definitions of~$\mathcal C^\times_G$ and $(\mathcal C^\times_G)^G$ in \Cref{sec:General theory of (De)Equivariantisation}. 

From any modular fusion category~$\mathcal C$ one can construct a Reshetikhin--Turaev defect TQFT as well as its associated 3-category~$\mathcal D_{\mathcal C}$. As reviewed in \Cref{exa:3dCategoriesWithAdjoints} this is the delooping of the 2-category of $\Delta$-separable symmetric Frobenius algebras in~$\mathcal C$. 
Given an orbifold datum in~$\mathcal D_{\mathcal C}$, one obtains a new TQFT which (as shown in \cite{CMRSS2021, Carqueville2021}) is another Reshetikhin--Turaev theory based on the modular fusion category~$\mathcal C_{\mathbb{A}}$ constructed in \cite{Mulevicius2022} and reviewed in \Cref{exa:3dOrbifoldData} below. 

One of the key technical steps towards our main results is that from any $G$-crossed braided fusion category~$\mathcal C^\times_G$ we obtain a 0-form symmetry 
\be 
	R\colon \textrm{B}\underline{G} \longrightarrow \mathcal D_{\mathcal C_e}
\ee 
which in turn gives rise to an orbifold datum~$\mathbb{A}_R$ in $\mathcal D_{\mathcal C_e}$, see \Cref{lem:0form orb dat from Gcbrc}. 
Hence we can gauge the 0-form symmetry~$R$. 
In \Cref{sec:equivalence of equivariantisation and orbifolding}, we prove that this gives the TQFT described by the equivariantisation $(\mathcal C^\times_G)^G$, as conjectured in \cite{Carqueville2020} (see also \cite{schweigert2019extendedhomotopyquantumfield} for a related result in a more geometric setting): 

\medskip 

\noindent
\textbf{\Cref{thm:eq is orbifolding}.}
Let~$\mathcal C^\times_G$ be a $G$-crossed braided fusion category, and let~$\mathbb{A}_R$ be the associated orbifold datum of \Cref{lem:0form orb dat from Gcbrc}. 
Then 
\be 
	\left(\mathcal{C}^\times_G\right)^G \; \cong \; \big(\mathcal{C}_e\big){}_{\mathbb{A}_R}
\ee 
as modular fusion categories. 

\medskip

We may now ask for symmetries of the theory described by $(\mathcal C^\times_G)^G$. 
As we explain in \Cref{sec:Deeq is 1-form orbifolding}, it has a canonical symmetry 
\be 
	\mathrm{B}^2\Rep(G)
	\longrightarrow 
	\mathcal{D}_{(\mathcal{C}^\times_G)^G}
	\, . 
\ee 
We show (in \Cref{prop:Deequivariantisation as orbifold datum1}) that this leads to an orbifold datum~$\mathbb{A}_B$ from which one recovers the neutral component: 
\be 
	\left((\mathcal{C}^\times_G)^G\right)_{\mathbb{A}_B}
	\;\cong\; 
	\mathcal{C}_e 
	 \, . 
\ee 

In the special case where~$G$ is commutative, we can identify $\Rep(G)$ with the category $\widehat{G}\text{-}\Vect$ of $\widehat{G}$-graded vector spaces, where $\widehat{G} = \textrm{Hom}_{\textrm{Grp}}(G,\textrm{U}(1))$ is the Pontryagin dual. 
Then we have a 1-form symmetry
\be 
	\widehat{R}\colon \mathrm{B}^2(\widehat{G}\text{-}\Vect)
	\longrightarrow 
	\mathcal{D}_{(\mathcal{C}^\times_G)^G} \, .
\ee 
Note that the domain of~$\widehat{R}$ is not $\textrm{B}^2\widehat{\underline{G}}$ but the double-delooping of the linearisation $\widehat{G}\text{-}\Vect$ of~$\widehat{\underline{G}}$. 
These two domains can be exchanged without loss of generality, however, since the codomain is $\Vect$-enriched.
We show (in \Cref{lem:ARhat} by checking that~$\mathbb{A}_{\widehat{R}}$ is indeed an orbifold datum) that this 1-form symmetry can always be gauged, and that it still recovers the neutral component~$\mathcal C_e$:
\medskip 

\noindent
\textbf{\Cref{thm:inversion of 0form by 1form}.}
Let~$\mathcal C^\times_G$ be a $G$-crossed braided fusion category for a commutative group~$G$, and let~$\mathbb{A}_{\widehat{R}}$ be the associated orbifold datum of \Cref{lem:ARhat}. 
Then 
\be 
	\left((\mathcal{C}^\times_G)^G\right)_{\mathbb{A}_{\widehat{R}}}
	\;\cong\; 
	\mathcal{C}_e 
\ee 
as modular fusion categories. 

\medskip 

In summary, gauging dual 0- and 1-form symmetries of Reshetikhin--Turaev TQFTs coming from $G$-crossed braided fusion categories is mutually inverse. 
This is as expected from the physics literature, which contains the more general statement that gauging a $p$-form symmetry gives rise to a theory with an $(n-p-2)$-form symmetry (given by ``Wilson defects'' \cite{Gaiotto_2015,BBFP2024}) whose gauging recovers the original, see e.g.\ \cite{bhardwaj2023lecturesgeneralizedsymmetries}. 
In our setting we have $p=0$ in $n=3$ dimensions and the resulting 1-form symmetry is given by Wilson lines, i.e.\ $\Rep(G)$.

\medskip 

\noindent
\textbf{Acknowledgements. } 
We are grateful to 
	Tudor Dimofte, 
	Christopher Lieberum, 
	Vincentas Mulevi{\v{c}}ius, 
	Lukas Müller, 
	Iordanis Romaidis 
		and 
	Ingo Runkel 
for insightful discussions and/or helpful comments on an earlier version of the manuscript. 
We also thank the authors of \cite{HPRW}, where they independently prove \Cref{thm:eq is orbifolding}, for coordinating their arxiv submission with ours. 
We acknowledge partial support from the German Science Fund DFG (Heisenberg Programme) and the Austrian Science Fund FWF (project no.\ P\,37046).

\section{Background and Reminders}
\label{sec:Background}

In this section we collect some basic algebraic notions that play central roles subsequently. 
In Section~\ref{subsec:2categories} we present our conventions for 2-functors and their higher maps, and we recall pivotal 2-categories. 
Section~\ref{subsubsec:2dOrbDat} reviews 2-dimensional orbifold data as certain Frobenius algebras, as well as the associated completion operation on pivotal 2-categories. 
In Section~\ref{subsubsec:3cats} we give our conventions for 3-functors, and we recall Gray categories with duals and their graphical calculus. 
Section~\ref{sec:2-groups} is concerned with the 3-category of 2-groups, and how it is equivalent to that of crossed modules. 
Finally, in Section~\ref{subsubsec:3dOrbDat} we describe 3-dimensional orbifold data in general as well as in the special case of 3-categories that one obtains from modular fusion categories. 

\medskip 

Throughout we assume all 1-categories to be idempotent complete. 
This is in particular the case for the category of finite-dimensional $\C$-vector spaces, which we denote by $\Vect$. 

Our conventions for $n$-categories with $n\in\{1,2,3\}$ are as follows. 
We usually denote the composition of $n$-morphisms by~$\circ$, that of $(n-1)$-morphisms by~$\otimes$, and that of $(n-2)$-morphisms by~$\boxtimes$. 
Note that~$\otimes$ is a family of $1$-functors, and~$\boxtimes$ is comprised of $2$-functors. 
The $(n-1)$-functors exhibiting identity morphisms on an object~$X$ are denoted by $I_X\colon \mathds{T}\longrightarrow \End(X)$ where~$\mathds{T}$ is the $(n-1)$-category with a single object and only identity morphisms. 
We denote associators by~$a$, unitors by~$r$ and~$l$, and pentagonators by~$\pi$.

For a monoidal $(n-1)$-category~$\mathcal{C}$ we denote compositions and the monoidal product in such a way that its delooping, i.e.\ the $n$-category $\mathrm{B}\mathcal{C}$ with precisely one object~$*$ with $\textrm{End}_{\mathrm{B}\mathcal{C}}(*) = \mathcal C$, respects the above conventions.

\subsection{2-categories}
\label{subsec:2categories}

In this section we briefly recall the main ingredients of the 3-category of 2-categories, as well as rigid and pivotal structures for 2-categories. 
For details and background we refer to the textbooks \cite{Benabou1967, Johnson2020} and the more casual lectures notes \cite[Sect.\,2.2]{Carqueville_2018}. 

\medskip 

We begin with a brief reminder on (higher) maps between 2-categories; see e.g.\ \cite[Chap.\,4]{Johnson2020} for details. 
Let $\mathcal B, \mathcal B'$ be (not necessarily strict) 2-categories. 
Recall that a \textsl{$2$-functor} $F \equiv (F,F^2,F^0)\colon \mathcal{B}\longrightarrow\mathcal{B}^\prime$ consists of a function~$F$ on objects, a $1$-functor 
	\begin{equation}
		F_{X,Y}\colon \Hom_\mathcal{B}(X,Y)
			\longrightarrow \Hom_{\mathcal{B}^\prime}\big(F(X),F(Y) \big)
	\end{equation}
for each pair of objects $X,Y \in \mathcal B$, a natural isomorphism 
\begin{equation}
	\label{eq:coherence structure for 2-functors}
	F^2_{X,Y,Z}\colon \otimes^\prime\circ\, (F_{Y,Z}\times F_{X,Y})\Longrightarrow F_{X,Z}\circ \otimes \, ,
\end{equation}
for each triple of objects, and a natural isomorphism $F^0_X\colon I^\prime_{F(X)}\Longrightarrow F_{X,X}\circ I_X$ for each $X\in\mathcal B$.
These data are subject to coherence conditions analogous to those of monoidal functors (which are equivalent to the special case of 2-functors between deloopings of monoidal categories). 
We usually suppress the subscripts for $F_{X,Y}$.

A \textsl{pseudonatural transformation} $\beta\colon F\Longrightarrow F^\prime$ between parallel $2$-functors $F,F^\prime\colon \mathcal{B}\longrightarrow\mathcal{B}^\prime$ consists of a $1$-morphism 
\begin{equation}
	\beta_X\colon F(X)\longrightarrow F^\prime(X)
\end{equation}
for each object $X\in\mathcal{B}$, and an invertible natural transformation $(\beta_Y)_* \circ F_{X,Y} \Longrightarrow F^\prime_{X,Y}\circ (\beta_X)^*$ for each pair of objects $X,Y$, whose components for each $1$-morphism $f\colon X\longrightarrow Y$ are $2$-morphisms~$\beta_f$ which fill the diagram 
	\begin{equation}
		\begin{tikzcd}[column sep= 50,row sep= 35]
			F(X)
			\arrow[d, "{\beta_X}"{left}]
			\arrow[r, "{F(f)}"]
			& F(Y)
			\arrow[dl, Rightarrow,shorten <= 2ex, shorten >= 2ex, "{\beta_f}"]
			\arrow[d, "{\beta_Y}"]\\
			F^\prime(X)
			\arrow[r, "{F^\prime(f)}"{below}]
			& F^\prime(Y)
		\end{tikzcd}.
	\end{equation}
These data are subject to unitality and naturality coherence conditions. 

A \textsl{modification} $\Gamma\colon \beta\LongRRightarrow \beta^\prime$ between parallel pseudonatural transformations $\beta, \beta^\prime\colon F\Longrightarrow F^\prime$ is given by a $2$-morphism
\begin{equation}
	\Gamma_X\colon \beta_X \longrightarrow \beta^\prime_X
\end{equation}
for each object~$X$, subject to a coherence axiom. 

\medskip 

A \textsl{right adjoint} of a 1-morphism $f\colon X\to Y$ in a 2-category is a 1-morphism $f^* \colon Y\to X$ such that there exist adjunction maps 
\begin{equation}
	\widetilde{\ev}_f = 
		\begin{tikzpicture}[very thick,scale=1.0,color=blue!50!black, baseline=.4cm]
			\draw[line width=0pt] 
			(3,0) node[line width=0pt] (D) {{\small$f^*$}}
			(2,0) node[line width=0pt] (s) {{\small$f\vphantom{f^*}$}}; 
			\draw[redirected] (D) .. controls +(0,1) and +(0,1) .. (s);
		\end{tikzpicture}
	\colon f \otimes f^* \longrightarrow  \id_Y
	\, , \quad 
	\widetilde{\coev}_f = 
		\begin{tikzpicture}[very thick,scale=1.0,color=blue!50!black,baseline=-.4cm,rotate=180]
			\draw[line width=0pt] 
			(3,0) node[line width=0pt] (D) {{\small$f^*$}}
			(2,0) node[line width=0pt] (s) {{\small$f\vphantom{f^*}$}}; 
			\draw[directed] (D) .. controls +(0,1) and +(0,1) .. (s);
		\end{tikzpicture}
	\colon \id_X \longrightarrow f^* \otimes f
\end{equation}
subject to the usual zigzag identities. 
Here we use the standard graphical calculus where our convention is that string diagrams are read from bottom to top and from right to left. 
If every 1-morphism in a 2-category~$\mathcal B$ has a right adjoint, then a choice of adjunction maps for all 1-morphisms gives~$\mathcal B$ a \textsl{(right) rigid} structure. 
Left adjoints~${}^*\!f$ are defined analogously, and we write $\ev_f, \coev_f$ for their adjunction maps. 

From a rigid structure on a 2-category we obtain a 2-functor $(-)^* \colon \mathcal B \to \mathcal B^{\textrm{op}}$ which is the identity on objects, and whose codomain has the 1-morphisms reversed. 
In this case a \textsl{pivotal} structure on~$\mathcal B$ is a pseudonatural transformation $\id_{\mathcal B} \Longrightarrow (-)^{**}$ whose 1-morphism components are identities. 
Then left and right adjoints are necessarily isomorphic, and one finds that a pivotal structure is equivalent to ${}^*\!f = f^*$ for all 1-morphisms~$f$, as well as 
\begin{equation}

\end{equation}
for all 2-morphisms $\Phi\colon h\to f$ and composable 1-morphisms $g,f$, see e.g.\ \cite[Sect.\,2.3]{Carqueville_2012}. 

Pivotal 2-categories naturally arise from 2-dimensional defect TQFTs, see \cite{Davydov2011} and \Cref{subsec:TQFTPerspective}. 

\begin{example}
	The 2-category $\textrm{B}\Vect$ of finite-dimensional vector spaces has a canonical pivotal structure given by the natural isomorphisms $V\cong V^{**}$ for all $V\in\Vect$. 
\end{example}

\subsection{2-dimensional Orbifold Data}
\label{subsubsec:2dOrbDat}

In this section we review basics about $\Delta$-separable symmetric Frobenius algebras in a given ambient pivotal 2-category. 
Together with their bimodules and bimodule morphisms they form a pivotal 2-category called orbifold completion, which we describe below. 

\medskip 

Let~$\mathcal{B}$ be a 2-category.
Recall that a \textsl{Frobenius algebra} is a 1-endomorphism $A\colon a\to a$ together with 2-morphisms 
\begin{equation}
	\mu=

	\colon A \longrightarrow \id_a
\end{equation}
subject to the conditions
\begin{equation}
	\label{eq:FrobeniusAlgebra}
	\tikzzbox{%
		%
	}%
	\, .
\end{equation}

An algebra is \textsl{separable} if there exists a section for the multiplication as a bimodule morphism, and \textsl{$\Delta$-separable} if this section is the comultiplication. As explained in \cite[Prop.\,3.1]{Mulevicius2022a}, a Frobenius algebra~$A$ is separable iff there exists a 2-morphism $\psi^2 \colon \id_a\to A$ such that $(\id_a \otimes [\mu \circ(\psi^2\otimes\id_a)])\circ\Delta$ is a section for~$\mu$. 
Hence we sometimes denote a separable Frobenius algebra as $(A,\psi)$. 
If there is a braiding, we can speak of commutative algebras, and if~$\mathcal{B}$ is pivotal, we consider symmetric Frobenius algebras. 
The defining conditions of these properties are as follows: 
\begin{equation}
\begin{array}{ccccccc}
	\text{separable}&&\Delta\text{-separable}&&\text{commutative}&&\text{symmetric }\\
\end{array}
\end{equation}
The \textsl{opposite algebra} of an algebra $(A,\mu,\dots)$ is given by $A^{\textrm{op}} := (A,\mu\circ c_{A,A},\dots)$, where~$c$ is the braiding. 
Hence commutative algebras are their own opposites. 

Certain Frobenius algebras play a central role in our constructions below, for which we adopt the nomenclature of \cite{Carqueville2012, Carqueville2017}: 

\begin{definition}
	Let~$\mathcal B$ be a pivotal 2-category. 
	A \textsl{(2-dimensional) orbifold datum} in~$\mathcal{B}$ is a $\Delta$-separable symmetric Frobenius algebra in~$\mathcal{B}$. 
\end{definition}

\begin{remark}
	\label{rem:2dOrbifoldData}
	In the context of 2-dimensional defect TQFTs, $\Delta$-separable symmetric Frobenius algebras arise as the algebraic input data for an internal state sum construction, see \cite{Fukuma1994, LAUDA_2007, Carqueville2012} for details and \cite{carqueville2023orbifoldstopologicalquantumfield} for a review. 
	We refer to this construction as a (generalised) \textsl{orbifold construction} (as it also describes topologically twisted sigma models whose targets are quotients of manifolds by group actions as another special case).
	
	The basic idea is to label the Poincar\'e dual of an oriented triangulation of a surface with the data of the algebra, namely 1-strata with~$A$ and 0-strata with~$\mu$ or~$\Delta$. 
	Then the condition of triangulation invariance is guaranteed by the defining conditions of $\Delta$-separable symmetric Frobenius algebras.
\end{remark}

A \textsl{left module} over an algebra~$A$ is a composable 1-morphism~$X$ together with a compatible left $A$-action, i.e.\ a 2-morphism
\begin{equation}
	\label{eq:Frob module structure}
 .
\end{equation}
If~$A$ is Frobenius, then~$X$ has an induced \textsl{comodule structure} given by
\begin{equation}
	\label{eq:Frob induced comodule structure}
	%
	\, . 
\end{equation}
Right modules are defined analogously, and an \textsl{$(A^\prime,A)$-bimodule} is a 1-morphism equipped with a left $A^\prime$-module structure and a right $A$-module structure which commute with each other.
Morphisms of (bi)modules are morphisms between the underlying objects that commute with the algebra actions.

Let~$A$ be a $\Delta$-separable Frobenius algebra, let~$Y$ be a left $A$-module and let~$X$ be a right $A$-module. 
The \textsl{relative tensor product} $Y\otimes_AX$ is the coequaliser of the two induced maps $Y\otimes A \otimes X \to Y\otimes X$. 
Since we assume all our categories to be idempotent complete, the relative tensor product can be computed by splitting the idempotent
\begin{equation}
	\label{eq:relative tensor idempotent}
	e \equiv e_{Y,A,X} := 

	\, .
\end{equation}
Indeed, as spelled out e.g.\ in \cite[Sect.\,2.3]{Carqueville2012}, there is a retraction and a section denoted by
\begin{equation}
	\label{eq:RetractionSection}
	r \equiv r_{Y,A,X} = 
	%
	\, , 
\end{equation}
respectively, such that $s\circ r = e$ and $r \circ s = \id_{Y\otimes_A X}$. 

The retraction~$r$ is \textsl{balanced} in the sense that
\begin{equation}
	%
\end{equation}
and~$s$ satisfies a similar cobalancing condition involving the induced comodule structures as in~\eqref{eq:Frob induced comodule structure}.
Maps out of the relative tensor product are isomorphic to balanced maps out of the ordinary horizontal composite, and maps into the relative tensor product are isomorphic to cobalanced maps: 
\begin{align}
	\Hom(Y\otimes_AX,K)& \stackrel{\cong}{\longrightarrow} \Hom^\textrm{bal}(Y\otimes X,K) \, ,
	&
	\Hom(K, Y\otimes_AX)& \stackrel{\cong}{\longrightarrow} \Hom^\textrm{cobal}(K, Y\otimes X) \, .
	\nonumber 
	\\
	f & \longmapsto 
	\begin{tikzpicture}[very thick,scale=0.5,color=blue!50!black, baseline]
		\pgfmathsetmacro{\yy}{0.4}
		\draw[line width=0pt] (1,-1.5) node[below] (X) {{\scriptsize$X$}} 
		(0,-1.5) node[below] (Y) {{\scriptsize$Y$}}
		(0.5,1.5) node[above] (XY) {{\scriptsize$K$}};
		\coordinate (Xm) at (1,0);
		\coordinate (Ym) at (0,0);
		\fill[color=blue!50!black] (0.5,0.85) circle (4pt) node[right] (meet) {{\scriptsize$f$}};
		\draw
		(X) -- (Xm)
		(Y) -- (Ym)
		(0.5,0) -- (XY);
		\draw[color=black] ($(Xm)+(\yy,0)$) -- ($(Ym)+(-\yy,0)$); 
	\end{tikzpicture}
	&
	g & \longmapsto 
	\begin{tikzpicture}[very thick,scale=0.5,color=blue!50!black, baseline, rotate=180]
		\pgfmathsetmacro{\yy}{0.4}
		\draw[line width=0pt] (1,-1.5) node[above] (X) {{\scriptsize$Y$}} 
		(0,-1.5) node[above] (Y) {{\scriptsize$X$}}
		(0.5,1.5) node[below] (XY) {{\scriptsize$K$}};
		\coordinate (Xm) at (1,0);
		\coordinate (Ym) at (0,0);
		\fill[color=blue!50!black] (0.5,0.85) circle (4pt) node[right] (meet) {{\scriptsize$g$}};
		\draw
		(X) -- (Xm)
		(Y) -- (Ym)
		(0.5,0) -- (XY); 
		\draw[color=black] ($(Xm)+(\yy,0)$) -- ($(Ym)+(-\yy,0)$);
	\end{tikzpicture}
	\label{eq:BalancedCobalanced}
\end{align}
Composition of (co)balanced morphisms is composition in the underlying category, and the identification \eqref{eq:BalancedCobalanced} preserves composition of a balanced morphism with a cobalanced morphism since $r\circ s=\id_{Y\otimes_A X}$.

There is a 2-category whose objects are $\Delta$-separable symmetric Frobenius algebras which naturally appears in a generalisation of state sum constructions for arbitrary 2-dimensional topological quantum field theories to include defects, introduced in \cite{Carqueville2012}: 

\begin{definition}[label=def:2orbcompletion]
	Let $\mathcal{B}$ be a pivotal $2$-category. 
	Its \textsl{orbifold completion} is the 2-category $\mathcal{B}_{\textrm{orb}}$ where
	\begin{itemize}
		\item 
		objects are 2-dimensional orbifold data ($\Delta$-separable symmetric Frobenius algebras), 
		\item 
		1-morphisms $A\to B$ are $(B,A)$-bimodules, 
		\item 
		2-morphisms are bimodule morphisms, 
		\item 
		horizontal composition is the relative tensor product, 
		\item 
		vertical composition is induced from~$\mathcal B$, 
		\item 
		the unit 1-morphism of $(a,A) \in \mathcal{B}_{\textrm{orb}}$ is~$A$ viewed as a bimodule over itself. 
	\end{itemize}
\end{definition}

This is indeed a completion, as originally shown in \cite[Prop.\,4.2\,\&\,4.7]{Carqueville2012} and described in terms of a universal property (essentially as the left adjoint of a forgetful 3-functor) in \cite[Sect.\,4.2.2]{Carqueville2023}: 

\begin{theorem}
	Let $\mathcal{B}$ be a pivotal $2$-category. 
	\begin{enumerate}[label={(\roman*)}]
		\item 
		\label{item:BorbPivotal}
		$\mathcal{B}_{\textrm{orb}}$ has a pivotal structure with adjunction data
		\begin{equation}
			\label{eq:adjunction data for modules}
			\ev_X = \!
			\begin{tikzpicture}[very thick,scale=1.0,color=blue!50!black, baseline=.8cm]
				\draw[line width=0pt] 
				(1.75,1.75) node[line width=0pt, color=green!50!black] (A) {{\small$A\vphantom{X^*}$}}
				(1,0) node[line width=0pt] (D) {{\small$X\vphantom{X^*}$}}
				(0,0) node[line width=0pt] (s) {{\small$X^*$}}; 
				\draw[directed] (D) .. controls +(0,1.5) and +(0,1.5) .. (s);
				
				\draw[color=green!50!black] (1.25,0.55) .. controls +(0.0,0.25) and +(0.25,-0.15) .. (0.86,0.95);
				\draw[-dot-, color=green!50!black] (1.25,0.55) .. controls +(0,-0.5) and +(0,-0.5) .. (1.75,0.55);
				\draw[color=green!50!black] (1.5,-0.1) node[Odot] (unit) {}; 
				\draw[color=green!50!black] (1.5,0.15) -- (unit);
				\fill (0.86,0.95) circle (2pt) node (meet) {};
				\draw[color=green!50!black] (1.75,0.55) -- (A);
			\end{tikzpicture}
			\circ s_{X^*,B,X}
			\, , \quad
			\coev_X =  r_{X,A,X^*} \circ 
			\begin{tikzpicture}[very thick,scale=1.0,color=blue!50!black, baseline=-.8cm,rotate=180]
				\draw[line width=0pt] 
				(3.21,1.85) node[line width=0pt, color=green!50!black] (B) {{\small$B\vphantom{X^*}$}}
				(3,0) node[line width=0pt] (D) {{\small$X\vphantom{X^*}$}}
				(2,0) node[line width=0pt] (s) {{\small$X^*$}}; 
				\draw[redirected] (D) .. controls +(0,1.5) and +(0,1.5) .. (s);
				\draw[color=green!50!black] (2.91,0.85) .. controls +(0.2,0.25) and +(0,-0.75) .. (B);
				\fill (2.91,0.85) circle (2pt) node (meet) {};
			\end{tikzpicture}
		\end{equation}
		for all $X \in \mathcal{B}_{\textrm{orb}}((a,A),(b,B))$, where we use the notation of~\eqref{eq:RetractionSection}. 
		\item 
		\label{item:2dCompletion}
		There is a pivotal equivalence $(\mathcal{B}_{\textrm{orb}})_{\textrm{orb}} \cong \mathcal{B}_{\textrm{orb}}$. 
	\end{enumerate}
	\label{thm:2dOrbifoldCompletion}
\end{theorem}

From now on we will usually not display the splitting data as in~\eqref{eq:adjunction data for modules}, but implicitly use the description of maps out of or into relative tensor products in terms of (co)balanced maps of~\eqref{eq:BalancedCobalanced}. 

\begin{example}
	\label{exa:2dOrbifoldCompletion}
	The orbifold completion of the pivotal 2-category $\textrm{B}\Vect$ is the pivotal 2-category of $\Delta$-separable symmetric Frobenius $\C$-algebras, their bimodules and bimodule maps, 
	\begin{equation}
			\big( \textrm{B}\Vect \big)_{\textrm{orb}} 
				\;\cong\;
			\textrm{$\Delta$ssFrobAlg}(\Vect) \, . 
	\end{equation}
	More generally, given any pivotal monoidal category~$\mathcal C$, we denote the pivotal 2-category $(\textrm{B}\mathcal C)_{\textrm{orb}}$ as $\textrm{$\Delta$ssFrobAlg}(\mathcal C)$, which is a full sub-2-category of the 2-category $\textrm{ssFrobAlg}(\mathcal C)$ of all separable symmetric Frobenius algebras in~$\mathcal C$. 
\end{example}

\subsection{3-categories}
\label{subsubsec:3cats}

In this section we briefly recall 3-functors between 3-categories, as well as coherent adjoints for 1- and 2-morphisms. 
For the general theory of 3-categories we refer to \cite{Gurski2013}, while our discussion of adjoints is based on the work \cite{Barrett2012} in the semi-strict setting. 

\medskip 

Let $\mathcal T, \mathcal T'$ be (not necessarily strict) 3-categories. 
Recall that a \textsl{3-functor} 
\begin{equation}
	\label{eq:3functor}
	R \equiv (R,\chi,\iota,\omega,\gamma,\delta)\colon \mathcal{T}\longrightarrow \mathcal{T}^\prime
\end{equation} 
consists of a function~$R$ on objects, 
a $2$-functor 
\begin{equation} 
	R_{a,b}\colon \Hom_\mathcal{T}(a,b)\longrightarrow \Hom_{\mathcal{T}^\prime}\big(R(a),R(b)\big)
\end{equation}
for each pair of objects $a,b\in\mathcal T$, 
a pseudonatural transformation 
\begin{equation} 
	\chi\colon \boxtimes^\prime\circ\, (R\times R)\Longrightarrow R\circ \boxtimes 
\end{equation}
for each triple of objects (which we suppress in the subscripts, and similarly in the other subscripts below), a pseudonatural transformation 
$
	\iota\colon I^\prime_{R(a)}\Longrightarrow R_{a,a}\circ I_a 
$
for each object, 
a modification
$
	\omega\colon \chi\circ (\chi\times \id)\LongRRightarrow \chi\circ (\id\times\chi)
$
for each quadruple of objects (suppressing associators), 
and two modifications
$
	\gamma\colon \chi\circ (\iota\times \id)\LongRRightarrow \id
$
and 
$\delta\colon \id\LongRRightarrow\chi\circ (\id\times \iota)$
for each pair of objects (suppressing unitors). 
These data are subject to various coherence conditions, see \cite[Def.\,4.10]{Gurski2013} for details.

\medskip 

Every 3-category is equivalent to a strict-up-to-tensorators \textsl{Gray category}, which is a category enriched over the category of strict 2-categories and strict 2-functors with the Gray tensor product, see e.g.\ \cite{Gurski2013}. 
Hence the Hom 2-categories of a Gray category~$\mathcal T$ are all strict, and the only not-necessarily-strict aspect of~$\mathcal T$ is the interchange law for the $\boxtimes$-composition of 1-morphisms. 

There is a 3-dimensional graphical calculus for Gray categories, developed in \cite{toddtrimblesurface, Barrett2012}. 
We use the conventions of \cite{Carqueville2016}, reading $\circ$-composition from bottom to top, $\otimes$-composition from right to left, and $\boxtimes$-composition from front to back. 
Hence for example 
\begin{equation}
	\tikzzbox{
}
\end{equation}
represents a 3-morphisms $\varphi \colon U \otimes (\id_r \boxtimes V) \to W \otimes (X \boxtimes \id_u)$, involving in particular the 1-morphism $u\colon a\to b$ and the 2-morphism $U \colon r \boxtimes v \to w$. 

The graphical calculus also works for \textsl{Gray categories with duals}. 
These were introduced in \cite{Barrett2012} as Gray categories whose Hom 2-categories come with pivotal structures $\id \Longrightarrow (-)^{**}$, and where every 1-morphism $X\colon a\to b$ has an adjoint $X^\#\colon b\to a$ with $X^{\#\#} = X$ and chosen adjunction 2-morphisms $\coev_X\colon \id_b \to X \boxtimes X^\#$ and $\ev_X \colon X^\# \boxtimes X \to \id_a$ such that $\ev_X = (\coev_{X^\#})^*$, and their zigzag identities hold up to coherent 3-isomorphisms. 
Graphically, we have 
\begin{equation}
	\id_{\ev_X} = 

\end{equation}
so that an $X$-labelled 2-sphere evaluates to the 3-morphism $\widetilde{\ev}_{\ev_X} \circ \coev_{\coev_{X^\#}}$. 

Gray categories with duals naturally arise from 3-dimensional defect TQFTs, see \cite{Carqueville2016} and \Cref{subsec:TQFTPerspective}. 
We will identify 3-categories with adjoints for 1- and 2-morphisms which are coherently equivalent to Gray categories with duals with the latter. 

\begin{example}
	\label{exa:3dCategoriesWithAdjoints}
	For any modular fusion category~$\mathcal C$, there is an associated Gray category with duals~$\mathcal D_{\mathcal C}$. 
	We review its construction below, for details we refer to \cite{kapustin2010surfaceoperators3dtopological, Koppen2021, Carqueville2023}. 
	
	Recall that a modular fusion category~$\mathcal C$ is a braided spherical fusion category whose braiding is non-degenerate, see e.g.\ \cite{EGNO2015, TuraevVirelizier}. 
	This means that~$\mathcal C$ is a semisimple $\C$-linear monoidal category with only finitely many isomorphism classes of simple objects. 
	A spherical structure is a pivotal structure with the property that left and right traces are equal, and a non-degenerate braiding $\{c_{X,Y} \colon X\otimes Y \to Y \otimes X\}_{X,Y\in\mathcal C}$ is such that if $c_{X,Y} \circ c_{Y,X} = \id_{Y\otimes X}$ for all $X\in\mathcal C$, then~$Y$ is isomorphic to ${\mathbbm{1}}^{\oplus s}$ for some $s\in\N$. 
	
	Using the notation of Example~\ref{exa:2dOrbifoldCompletion}, we have a pivotal 2-category $\textrm{ssFrob}(\mathcal C)$ for every modular fusion category~$\mathcal C$. 
	Thanks to the tensor product and braiding of~$\mathcal C$, this 2-category has a natural monoidal structure so that we may consider its delooping
	\begin{equation}
		\label{eq:defect3cat associated to an MFC}
		\mathcal D_{\mathcal C} := \textrm{B}\,\textrm{ssFrob}(\mathcal C) \, . 
	\end{equation}
	Unpacking this, we have that in~$\mathcal D_{\mathcal C}$ 
	\begin{itemize}
		\item 
		there is only a single object, 
		\item 
		1-morphisms are separable symmetric Frobenius algebras in~$\mathcal C$, 
		\item 
		2-morphisms are bimodules in~$\mathcal C$, 
		\item 
		3-morphisms are bimodule morphisms in~$\mathcal C$, 
		\item 
		Hom 2-categories come with a pivotal structure (induced via Theorem~\ref{thm:2dOrbifoldCompletion}\ref{item:BorbPivotal}), 
		\item 
		adjoints of 1-morphisms are given by opposite algebras, 
	\end{itemize}
	and~$\mathcal D_{\mathcal C}$ is (equivalent to) a Gray category with duals. 
\end{example}

Note that the construction of~$\mathcal D_{\mathcal C}$ is well-defined for arbitrary ribbon categories~$\mathcal C$, modularity is not required.

\subsection{2-groups}
\label{sec:2-groups}

A \textsl{$2$-group} is a monoidal 1-category where all objects and morphisms are invertible. 
When working with explicit examples, it is often useful to equivalently describe 2-groups in terms of purely group theoretic data. 
Here we recall this relationship. 

\medskip 

A \textsl{crossed module} $(G,H,t,\varphi)$ consists of two groups~$G$ and~$H$, as well as group homomorphisms $t\colon H\longrightarrow G$ and $\phi\colon G^{\textrm{rev}}\longrightarrow \Aut(H)$,\footnote{By the superscript ``rev'' we indicate the reversal of the group multiplication or the monoidal product. 
This means that we use right $G$-actions in accordance with later sections and \cite[Sect.\,5]{Carqueville2020}.} and we write $\phi_g:=\phi(g)$. 
These data are subject to the compatibility conditions 
\begin{align}
	\label{eq:Peiffer identity}
	\phi_{t(h)}(h^\prime) &= h^{-1}h^\prime h 
		\, , \\
		\label{eq:equivariance of t}
	t(\phi_g(h))& = g^{-1}t(h)g
\end{align}
for all $g\in G$ and all $h,h^\prime\in H$. 

There is a $3$-category 2Grp whose objects are deloopings of 2-groups, and whose 1-, 2-, and 3-morphisms are 2-functors, pseudonatural transformations, and modifications, respectively. 
There is also a $3$-category $\chi$Mod of crossed modules, whose 1-morphisms are so-called butterflies (because of the shape of their defining diagram), 2-morphisms are butterfly-compatible group homomorphisms $G\longrightarrow H^\prime$, and 3-morphisms are elements of~$H^\prime$. 
See \cite{Aldrovandi2009} for butterflies and \cite[Sect.\,3.3\,\&\,Def.\,8.4--8.6]{Noohi2008} for details on the equivalence; below we are exclusively interested in mapping objects from $\chi$Mod to 2Grp. 

There is an equivalence of $3$-categories 
\begin{align}
	\chi\textrm{Mod} & \stackrel{\cong}{\longrightarrow} \textrm{2Grp} 
	\nonumber
	\\
	(G,H,t,\phi) & \longmapsto \mathrm{B}\mathcal G_{(G,H,t,\phi)}
\end{align}
where
\begin{align}
	\mathrm{Ob}\,\mathcal G_{(G,H,t,\phi)} &:= G 
		\, , 
	\\
	\Hom_{\mathcal G_{(G,H,t,\phi)}}(g,g^\prime) &:= \big\lbrace (g,h) \;\big|\; t(h)g=g^\prime\,\big\rbrace\, .
	\label{eq:Hom2group}
\end{align}
Identity morphisms in $\mathcal G_{(G,H,t,\phi)}$ are given by $(g,e_H)$, and composition is induced from multiplication in~$H$, i.e.\ $(t(h)g,h^\prime)\circ (g,h) := (g,h^\prime h)$. 
In the 3-category~$\mathrm{B}\underline{\mathcal{G}}{}_{(G,H,t,\phi)}$, where $\underline{\mathcal{G}}_{(G,H,t,\phi)}$ denotes the monoidal 2-category given by~$\mathcal{G}_{(G,H,t,\phi)}$ with additional identity 2-morphisms, this is presented graphically as 
\begin{equation}
	\label{eq:vertical composition in 2group}

\end{equation} 
where here and below we label lines that represent a 2-morphism $(g,h)$ (which is a morphism in $\mathcal{G}_{(G,H,t,\phi)}$) simply by~$h$. 
The strict monoidal structure of~$\mathcal G_{(G,H,t,\phi)}$ is given by group multiplication in~$G$ on objects and by the semidirect product 
$
(g,h)\otimes(g^\prime,h^\prime):=(gg^\prime,h\phi_{g^{-1}}(h^\prime)) \in\Hom(gg^\prime, t(h)g\,t(h^\prime)g^\prime)
$ 
on morphisms, hence in $\mathrm{B}\underline{\mathcal{G}}{}_{(G,H,t,\phi)}$ it reads 
\begin{equation}
	\label{eq:semidirect product of 2group}

\end{equation}
and the unit object is the unit element in~$G$. Conversely, starting from the 2-group, this serves as the definition of~$\phi$ by setting $h=e_H$.
From \eqref{eq:vertical composition in 2group} we can read off that~$t$ has to be a morphism of groups, and from \eqref{eq:semidirect product of 2group} we can similarly deduce that it is $G$-equivariant, since $t(h)g\, t(h')g'=t(h\phi_{g^{-1}}(h'))gg'$ implies \eqref{eq:equivariance of t}.
The Peiffer identity \eqref{eq:Peiffer identity} follows from the two ways to compute the central diagram
\begin{equation}

\end{equation}
by composing the lines on the back plane first (right-hand side) or performing the composition of surfaces first (left-hand side). This uses the fact that all coherence data in $\mathcal{G}_{(G,H,t,\phi)}$ is given by identities.

Unless specified otherwise, from now on by a 2-group we mean one of type $\mathcal G_{(G,H,t,\phi)}$.

\subsection{3-dimensional Orbifold Data}
\label{subsubsec:3dOrbDat}

Recall from \Cref{subsubsec:2dOrbDat} that 2-dimensional orbifold data are $\Delta$-separable symmetric Frobenius algebras in a given pivotal 2-category. 
As indicated in \Cref{rem:2dOrbifoldData} they arise by demanding invariance under 2-dimensional Pachner moves. 
Here we review the analogous discussion in three dimensions, following \cite{Carqueville2012, Carqueville2017, Carqueville2023}. 

\begin{definition}
	\label{def:3dOrbifoldDatum}
	Let~$\mathcal D$ be a Gray category with duals. 
	A \textsl{(3-dimensional, special) orbifold datum\footnote{What we call an orbifold datum in~$\mathcal D$ here is what is equivalent to what is called an orbifold datum in the ``Euler completion'' of~$\mathcal D$ in \cite{Carqueville2017}. Since here we are exclusively interested in the Euler completion (which amounts to allowing non-trivial ``point insertions''~$\psi$ and~$\phi$ in 2-strata and in 3-strata, respectively), we permit ourselves the shorthand.}} in~$\mathcal D$ is a tuple 
	\begin{equation}
		\mathbb{A} := \left(a, \;\;
			\tikzzbox{

			\, , \;\;
			\psi 
			\, , \;\;
			\phi 
		\right) 
	\end{equation}
	where $a\in\mathcal D$ is an object, $A\colon a \to a$ is a 1-morphism, $T \colon A\boxtimes A \to A$ is a 2-morphism, $\alpha \colon T \otimes (\id_A \boxtimes T) 
	\to T \otimes (T \boxtimes \id_A)$, $\psi\in\textrm{End}(\id_A)$, and $\phi\in\textrm{End}(\id_{\id_a})$ are 3-isomorphisms, and they are subject to the relations in Figure~\ref{fig:OrbifoldDatumAxioms}.
\end{definition}

\begin{figure}[!ht]
	\captionsetup[subfigure]{labelformat=empty}
	\centering
	\vspace{-60pt}
	\begin{subfigure}[b]{0.6\textwidth}
		\centering
		\begin{equation*}
		\tag{O1}
		\label{eq:O1}
		\tikzzbox{
}
		\end{equation*}
	\end{subfigure}
	\vspace{0pt}
	\caption{Defining conditions on 3-dimensional orbifold data $\mathbb A = (a,A,T,\alpha,\overline{\alpha},\psi,\phi)$; striped surfaces correspond to the adjoint of~$A$ and the (double) primed versions of $\alpha,\overline\alpha$ are obtained by pre- and post-composition with adjunction maps for~$T$, e.g.\ $\alpha^\prime=${\protect\tikz[semithick,scale=.3,color=green!30!black,baseline=-.1cm]{
	\protect\draw[color=green!30!black, semithick, rounded corners=0.5pt, postaction={decorate}, decoration={markings,mark=at position .8 with {\protect\arrow[draw=green!30!black]{>}}}] 
		(0,0) .. controls +(-0.25,0.25) and +(0,-0.25) .. (-0.5,1);
		\protect\draw[color=green!30!black, semithick, rounded corners=0.5pt, postaction={decorate}, decoration={markings,mark=at position .2 with {\protect\arrow[draw=green!30!black]{>}}}] 
		(1.25,1) .. controls +(-0.25,-1.5) and +(0.25,-0.75) .. (0,0);
		\protect\draw[color=white, line width=3pt] 	(1.25,-1) .. controls +(-0.25,1.5) and +(0.25,0.75) .. (0,0);
		\protect\draw[color=green!30!black, semithick, rounded corners=0.5pt, postaction={decorate}, decoration={markings,mark=at position .2 with {\protect\arrow[draw=green!30!black]{<}}}] 
		(1.25,-1) .. controls +(-0.25,1.5) and +(0.25,0.75) .. (0,0);
		\protect\draw[color=green!30!black, semithick, rounded corners=0.5pt, postaction={decorate}, decoration={markings,mark=at position .5 with {\protect\arrow[draw=green!30!black]{>}}}] 
		(-0.5,-1) .. controls +(0,0.25) and +(-0.25,-0.25) .. (0,0);
		\protect\fill (0,0) circle (3.8pt) node[black, opacity=0.6, left] {{\tiny$\overline{\alpha}$}};
		}}, see \cite[Def.\,4.2]{Carqueville2017} for details.}
	\label{fig:OrbifoldDatumAxioms}
\end{figure}

\begin{remark}
	\begin{enumerate}
		\item 
		In the context of 3-dimensional defect TQFTs, orbifold data~$\mathbb A$ as above can be used to label the Poincar\'e duals of oriented triangulations of 3-manifolds: $(a,\phi)$ then labels 3-strata, $(A,\psi)$ labels 2-strata, $T$ labels 1-strata, and $\alpha,\overline{\alpha}$ label positively/negatively oriented 0-strata. 
		As explained in \cite{Carqueville2017}, the defining conditions of orbifold data in Figure~\ref{fig:OrbifoldDatumAxioms} ensure independence of the choice of triangulation in the construction. 
		In the special cases relevant for our applications in subsequent sections, $(A,\psi)$ turns out to be a separable symmetric Frobenius algebra as in \Cref{subsubsec:2dOrbDat}. 
		\item 
		Just like every 2-dimensional orbifold datum is in particular an associative algebra in its ambient 2-category, every 3-dimensional orbifold datum~$\mathbb{A}$ is in particular a (weakly) associative algebra in its ambient 3-category. 
		Indeed, it follows from the conditions~\eqref{eq:O2} and~\eqref{eq:O3} that $\alpha,\overline\alpha$ are invertible, and hence the ``multiplication''~$T$ is weakly associative. 
		Moreover, condition~\eqref{eq:O1} says that~$\alpha$ and~$\overline{\alpha}$ satisfy the pentagon identity up to the isomorphism~$\psi^2$. 
	\end{enumerate}
\end{remark}

\begin{example}
	\label{exa:3dOrbifoldData}
	Let~$\mathcal C$ be a ribbon category, and let~$\mathcal D_{\mathcal C}$ be the Gray category with duals as in \Cref{exa:3dCategoriesWithAdjoints}. 
	As shown in \cite[Sect.\,3.2]{Carqueville2017} and \cite[Sect.\,5.2]{Mulevicius2022a}, an orbifold datum~$\mathbb A$ in~$\mathcal D_{\mathcal C}$ is equivalently described by
	\begin{itemize}
		\item 
		a separable symmetric Frobenius algebra $(A,\psi)$ in~$\mathcal C$, 
		\item 
		an $(A,A\otimes A)$-bimodule~$T$ in~$\mathcal C$, 
		\item 
		$(A,A\otimes A\otimes A)$-bimodule morphisms $\alpha\colon T\otimes_2 T \to T\otimes_1 T$ and $\overline\alpha\colon T\otimes_1 T\to T\otimes_2 T$, 
		\item 
		an invertible scalar $\phi \in \C^\times$, 
	\end{itemize}
	subject to the relations in Figure~\ref{fig:OrbifoldDatumAxiomsInDC}. 
	Here for $i\in\{1,2\}$ by the notation $T \otimes_i T$ we mean the relative tensor product over the $i$-th factor~$A$ acting from the right on the left~$T$ (and from the left on the right~$T$), and we have already evaluated the tensor product $T\otimes_1 A\cong T$ which appears in the domain of~$\alpha$ in \Cref{def:3dOrbifoldDatum}. 
	We set 
	\begin{equation}
		\label{eq:PsiMaps}
 
		\, . 
	\end{align}
	\vspace{0pt}
	\caption{Defining conditions on 3-dimensional orbifold data in $\mathcal D_{\mathcal C}$.}
	\label{fig:OrbifoldDatumAxiomsInDC}
\end{figure}

\begin{figure}
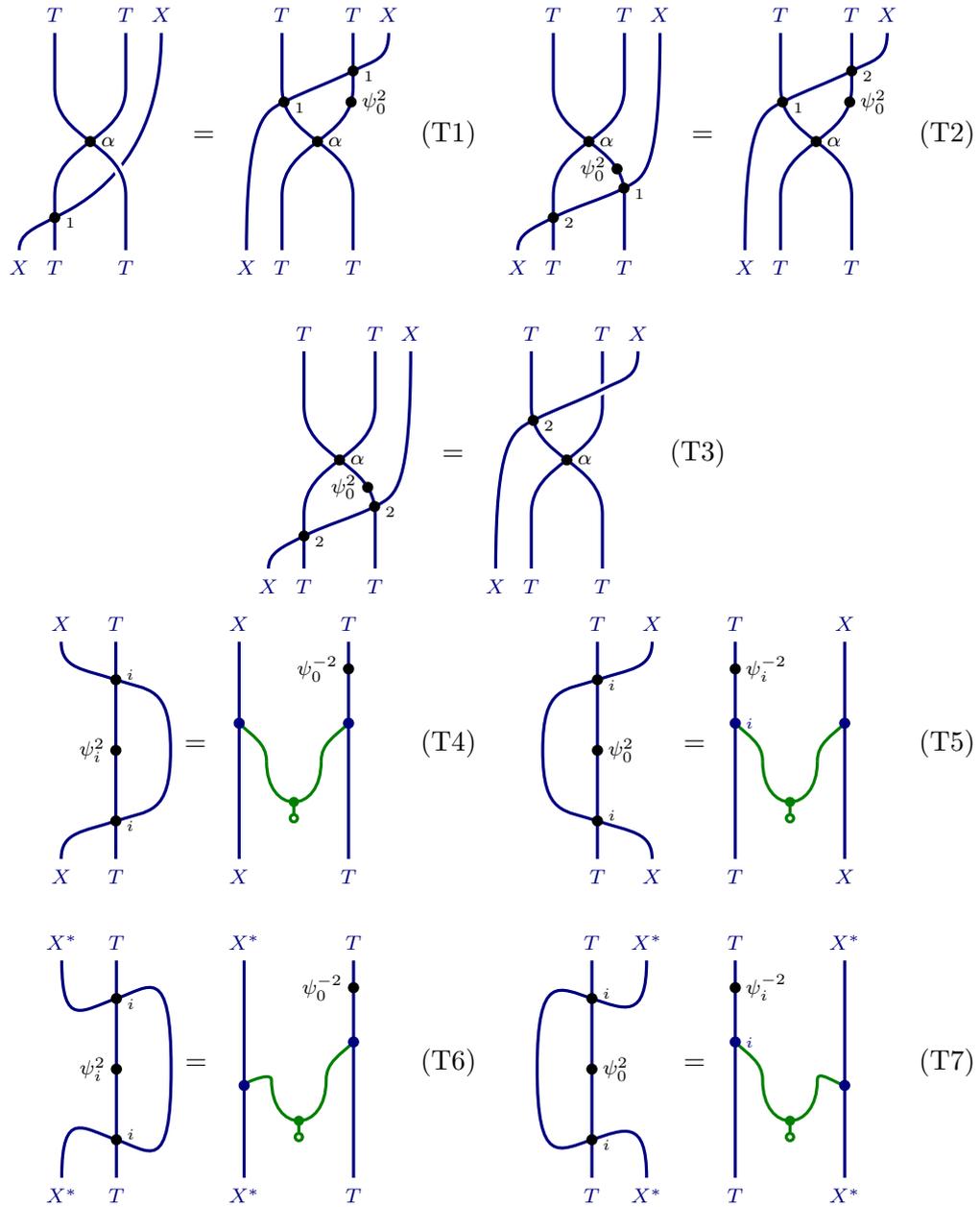
 
	\captionsetup[subfigure]{labelformat=empty}
	\centering
	\begin{subfigure}{0.45\textwidth}
	\begin{equation}
		\tag{T1}
		\label{eq:T1}
		%
 
	\end{equation}
	\end{subfigure}
	\caption{Defining conditions for $\tau$-structures.}
	\label{fig:tau-structure axioms}
\end{figure}

\begin{remark}
	\label{rem:MFCCA}
	Given an orbifold datum $\mathbb{A} = (*,A,T,\alpha,\overline\alpha,\psi,\phi)$ for a ribbon category~$\mathcal C$ as in \Cref{exa:3dOrbifoldData}, one can construct a new, ``orbifolded'' ribbon category~$\mathcal C_{\mathbb{A}}$ as first carried out in \cite{Mulevicius2022}. 
	Objects of~$\mathcal C_{\mathbb{A}}$ are $(A,A)$-bimodules~$X$ together with $(A,A^{\otimes 3})$-bimodule morphisms $\tau^X_i\colon X\otimes_A T \to T\otimes_i X$ and $\overline\tau^X_i\colon T\otimes_i X\to X\otimes_A T$ for $i\in\{1,2\}$, subject to the conditions in \Cref{fig:tau-structure axioms}, cf.\ \cite[Fig.\,3.1]{Mulevicius2022}. 
	The maps are called \textsl{$T$-crossings}, but we may also refer to them as a \textsl{$\tau$-structure on~$X$}. Morphisms $f\colon (X,\tau^X)\to (Y,\tau^Y)$ are $(A,A)$-bimodule morphisms such that 
	\begin{equation}
	\label{eq:morphisms in CAA}
	\tau^Y_i \circ (f\otimes_A \id_T) = (\id_T \otimes_i f)\circ\tau^X_i.
	\end{equation}		
	Using the description of maps of relative tensor products in terms of (co)balanced maps as in~\eqref{eq:BalancedCobalanced}, we represent the maps $\tau^X_i, \overline\tau^X_i$ in string diagrams as 
	\begin{equation} 
		\tau_i^X \;=\; 
 
		\, . 
	\end{equation}
	The monoidal structure on~$\mathcal C_{\mathbb A}$ is given by the relative tensor product $X \otimes_A Y$ together with the $T$-crossings 
	\begin{equation}
	\label{eq:CAA monoidal product tau}
	\tau^{X,Y}_i:=
		%
 
		\, , 
	\end{equation}
	using the notation introduced in~\eqref{eq:PsiMaps}. The monoidal unit is~$A$ as an $(A,A)$-bimodule whose $T$-crossings~$\tau_i^\mathds{1}$ are given by the composition of module and $i$-th comodule structure of~$T$ and lastly $\omega_i^{-1}$.
	The braided structure of~$\mathcal C_{\mathbb A}$ has components 
	\begin{equation}
	\label{eq:CAA braiding}
		c^{\mathcal{C}_\mathbb{A}}_{X,Y}:=\phi^2 \cdot 
		%
 
	\end{equation}
	viewed as maps $X\otimes_A Y \to Y\otimes_A X$ and we used the endomorphisms~$\omega$ for $(A,A)$-bimodules in a direct generalisation of \eqref{eq:PsiMaps}. 
	Finally, the rigid structure is such that the dual of $(X,\tau_i^X, \overline\tau_i^X) \in \mathcal C_{\mathbb A}$ is the dual $(A,A)$-bimodule~$X^*$ in~$\mathcal C$ together with the $T$-crossings 
	\begin{equation}
	\label{eq:dual tau}
	\tau^{X^*}_i:=
		%
	\end{equation}
	and the adjunction data in \cite[Eq.\,(3.7)]{Mulevicius2022}. 
	The morphisms underlying the adjunction data are those of \eqref{eq:adjunction data for modules}, but with additional $\psi$-insertions.
	For further details we refer to \cite[Sect.\,3]{Mulevicius2022}.
\end{remark}

Recall the notion of orbifold completion of a pivotal 2-category~$\mathcal B$ of \Cref{def:2orbcompletion} which produces a new 2-category $\mathcal B_{\textrm{orb}}$ whose objects are 2-dimensional orbifold data.
Analogously, there is a notion of orbifold completion for a Gray category with duals~$\mathcal D$ which produces a new 3-category~$\mathcal D_{\textrm{orb}}$ whose objects are 3-dimensional orbifold data.
This completion is constructed in \cite{Carqueville2023} with morphisms in~$\mathcal D_{\textrm{orb}}$ given by (weak) bimodules, bimodule morphisms and their modifications, where the bimodules and bimodule morphisms satisfy additional conditions, and it is shown that its Hom 2-categories are pivotal and orbifold complete, and every 1-morphism in~$\mathcal D_{\textrm{orb}}$ has equivalent left and right adjoints. 
Generalising the case of 1- and 2-categories, $(-)_{\textrm{orb}}$ is expected to be left adjoint to a forgetful 4-functor. 

We will only need the full sub-3-category 
\begin{equation}
	\label{eq:DHomorb}
	\Homorb{\mathcal D} \;\subset\; \mathcal D_{\textrm{orb}}
\end{equation}
whose objects are the trivial orbifold data $a \equiv (a,\id_a,l_{\id_a},\dots)$ whose structure morphisms are made up of the structure morphisms of~$\mathcal D$. 
Hence we can identify the objects of $\Homorb{\mathcal D}$ with those of~$\mathcal D$, and the Hom 2-categories of the former are the 2-dimensional orbifold completions of the Hom 2-categories of the latter, 
\begin{equation}
	\label{eq:Hom-orbData}
	\textrm{Ob}\big( \Homorb{\mathcal D} \big) 
		\cong 
		\textrm{Ob}\big( \mathcal D \big)
	\, , \quad 
	\Hom_{\Homorb{\mathcal D}}(a,b) 
		= 
		\big(\Hom_{\mathcal D}(a,b)\big)_{\textrm{orb}} \, . 
\end{equation}

\section{Orbifold Data from 2-Group Symmetries}
\label{sec:OrbDatFrom2Groups}

In this section we consider 2-group symmetries~$R$ (see \Cref{def:3d 2-group symmetry}) internal to certain 3-categories, and we associate candidate orbifold data~$\mathbb{A}_R$ to them. 
In \Cref{subsec:GeneralCase} we prove that~$\mathbb{A}_R$ is a weakly associative algebra and interpret its potential failure to be an orbifold datum as an obstruction to gauging the symmetry. 
In \Cref{sec:0-form specialisation} we briefly specialise to 0-form symmetries, while in \Cref{sec:1-form specialisation} we study the other special case of 1-form symmetries from braided fusion categories. 
We show that in the latter case the symmetry can always be gauged if the symmetry functor is pivotal, i.e.\ then~$\mathbb{A}_R$ is necessarily an orbifold datum, and we describe how it behaves under transport along ribbon functors.

\subsection{General Case}
\label{subsec:GeneralCase}

From now on we assume that $\mathcal{D}$ is a Gray category with duals whose $\Hom$ 2-categories contain direct sums and are $\Vect$-enriched, or~$\mathcal D$ is equivalent to such a 3-category. 
Let~$\mathcal{G}$ be a 2-group described by a crossed module $(G,H,t,\phi)$ as in \Cref{sec:2-groups}. 
By~$\underline{\mathcal{G}}$ we denote the monoidal 2-category given by~$\mathcal{G}$ with additional identity 2-morphisms.

\begin{definition}[label=def:3d 2-group symmetry]
	A \textsl{$\mathcal{G}$-symmetry in~$\mathcal{D}$} is a 3-functor
	\begin{equation}
		\label{eq:3d 2-group symmetry}
		R\colon \mathrm{B}\underline{\mathcal{G}}\longrightarrow \mathcal{D}\, .
	\end{equation}
	We also refer to such an~$R$ as a \textsl{$\mathcal{G}$-symmetry on (the object) $R(\ast)\in\mathcal D$}.
\end{definition}

In this section we work with general 2-group symmetries and construct candidate orbifold data from them. 
Then in \Cref{sec:0-form specialisation,sec:1-form specialisation} we specialise explicitly to 0-form and 1-form symmetries, i.e.\ trivial~$H$ and trivial~$G$, respectively. 

\medskip 

A 3-functor~$R$ as in~\eqref{eq:3d 2-group symmetry} picks out certain 1-endomorphisms $R(\ast)\longrightarrow R(\ast)$ and 2-morphisms between these both of which are associated to elements of the 2-group. 
In this way, the 3-functor~$R$ describes an action of~$\mathcal{G}$ on an object in~$\mathcal{D}$. 
This is part of the reason for calling it a ``symmetry'', analogous to the 1-categorical case of functors $\textrm{B}G \to \Vect$. 
We discuss the physical interpretation of this functor in \Cref{subsec:TQFTPerspective}. 

Before turning to the algebraic constructions of this section, we briefly motivate them from a TQFT perspective. 
To a first approximation, the goal is to construct an orbifold datum from every $2$-group $\mathcal{G}$-symmetry. 
This would mean that we could gauge the symmetry.
However, since there are anomalous symmetries (which cannot be gauged), we cannot expect to achieve this in general. 
We can, however, construct a weakly associative algebra, whose potential failure to be an orbifold datum describes an anomaly. 

\medskip

As discussed in \Cref{rem:2dOrbifoldData}, the orbifold construction involves labelling strata in the dual of a triangulation with objects and morphisms from the category of defects, namely those given by the orbifold datum. 
In our present situation, we want to enhance this by also choosing triangulations on the 2-strata, which can be thought of as a 2-step orbifolding process. 
The first step concerns the objects of~$\mathcal{G}$ which make up the elements of the group~$G$, while the second step is about the morphisms of~$\mathcal G$, i.e.\ the elements of~$H$. 
In this interpretation, it is natural to work with 1-strata on the 2-strata of the first orbifold datum, as these are the line defects in the orbifolded theory after the first step \cite{CMRSS2021}. 
As a result, the line defects on the 2-strata must be such that the construction is independent of the choice of triangulations of the 2-strata, hence they must be 2-dimensional orbifold data. 
Therefore, the orbifold datum candidate that we construct lives in the 3-category $\Homorb{\mathcal{D}}$ discussed at the end of \Cref{subsubsec:3dOrbDat}.
In the following, we set the TQFT perspective aside and focus purely on the algebraic constructions that arise from these considerations.

\medskip 

We now turn to construct an algebra $\mathbb{A}_R = (a,A,T,\alpha,\overline{\alpha})$ in $\Homorb{\mathcal{D}}$ (recall~\eqref{eq:Hom-orbData}), for which we define 
\begin{equation}
	\label{eq:object of general 2group orb dat}
	a:=R(\ast)
\end{equation}
to be the underlying object. 
The 1-endomorphism~$A=(A_G,A_H)$ of~$a$ in $\Homorb{\mathcal D}$ consists of a 1-morphism $A_G\colon a\to a$ in~$\mathcal{D}$, and a Frobenius algebra~$A_H$ (consisting of 2- and 3-morphisms in~$\mathcal D$) on~$A_G$. 
The 1-morphism~$A_G$ and the 2-morphism underlying~$A_H$ are given by 
\begin{equation}
	\label{eq:1-morphism of general 2group orb dat}
	A:= (A_G,A_H):=\Big(\bigoplus_{g\in G}R(g),\bigoplus_{\substack{g\in G\\h\in H}}R(g,h)\Big) , 
\end{equation}
where we recall from~\eqref{eq:Hom2group} that morphisms in our 2-group~$\mathcal G$ are pairs $(g,h)$. 
The Frobenius algebra structure of~$A_H$ is given by the multiplicative coherence data of~$R$ with respect to $\otimes$-composition. 
Recall from \eqref{eq:coherence structure for 2-functors} that these data consist of 3-isomorphisms
\begin{align}
	R^2_{g,h,h^\prime}:=R^2_{(t(h)g,h^\prime),(g,h)}\colon R(t(h)g,h^\prime)\otimes R(g,h)&\longrightarrow R(g,h^\prime h)\, ,
	\\
	R^0_g\colon \id_{R(g)}&\longrightarrow R(g,e_H)\, .
\end{align}
Using these we can define the Frobenius algebra structure on~$A_H$ to be 
\begin{align}
	\!\!\!\mu&:=
\!
	:=|H|\!\Big(\bigoplus_{g\in G} \left(R^0_g\right)^{-1}\Big)\circ \mathrm{pr}\, ,
	\label{eq:Frob structure on AH}
\end{align}
where the projection and inclusion are those of $\bigoplus_{g\in G}R(g,e_H)$ in~$A_H$.

\begin{lemma}[label=lem:Frobenius algebra structure on AH]
Let~$\mathcal{G}$ be a $2$-group associated to a crossed module $(G,H,t,\phi)$, let~$\mathcal{D}$ be a $3$-category whose $\Hom$ 2-categories are pivotal, contain direct sums and are $\Vect$-enriched, and let $R\colon \textrm{B}\underline{\mathcal{G}}\longrightarrow \mathcal{D}$ be a $3$-functor whose 2-functor components on $\Hom$ 2-categories are pivotal. 
Then $(A_H,\mu,\Delta,\eta,\varepsilon)$ as in \eqref{eq:1-morphism of general 2group orb dat}, \eqref{eq:AHMuDelta} and \eqref{eq:Frob structure on AH} is a $\Delta$-separable symmetric Frobenius algebra.
\end{lemma}

\begin{proof}
	The algebra is associative and coassociative because of the hexagon (i.e.\ associativity) axiom for~$R^2$, and it is unital and counital because of the square (i.e.\ unitality) axiom for~$R^0$ and~$R^2$. 
	
	Since~$\mu$ and~$\Delta$ are sums whose terms are inverses of each other, their composition $\mu\circ\Delta$ is a sum of identities. 
	The multiplicity is cancelled by the normalisation of~$\Delta$, which implies $\Delta$-separability. 
	Symmetry follows from the fact that the $2$-functors on $\Hom$ 2-categories are pivotal, i.e.\ the adjunction data is given by the composition of counit with multiplication and the composition of comultiplication with unit (up to a unique isomorphism which is identical for left and right adjoints). 
	The zigzag identity then implies the claim. 
	
	To prove the Frobenius property, we introduce a simplified notation. 
	The relevant diagrams are ones of 2-endomorphisms of~$A_G$, so we can use 2-dimensional diagrams. 
	We denote $R(g)$ by~$g$ and only label the 2-strata once, the same labels apply for each diagram. 
	We denote a line $R(g,h)$ simply by~$h$ since the element~$g\in G$ is uniquely determined, namely it coincides with the label of the plane to the right of the line. 
	The Frobenius property is an identity of 3-morphisms $A_H\otimes A_H\longrightarrow A_H\otimes A_H$, whose components are of the form $h_1^\prime\otimes h_1\longrightarrow h_2^\prime\otimes h_2$ in our simplified notation. 
	We prove the Frobenius identity for each component by first applying the zigzag identity on~$h^\prime_2$ and unitality on~$h_2$, then using invertibility of~$R^0$, associativity, and lastly invertibility of~$R^2$:
	\begin{equation}
		\label{eq:LongGreenIdentities}
		\!\!\!

	\end{equation}
	Note that, because we are looking at components, the inclusion and projection in the definition of~$\eta$ and~$\varepsilon$ are omitted, so the unit and counit reduce to~$R^0$ and its inverse, thus becoming invertible. 
	Similarly, since we are not summing over~$H$ in these diagrams, no multiplicities need to be cancelled when applying $\Delta$-separability, and thus we omit factors of~$|H|$. 
	
	Continuing the calculation of~\eqref{eq:LongGreenIdentities}, we introduce another comultiplication and counit on the~$h^\prime_1$-strand and then essentially repeat the previous steps to obtain the desired result:
	\begin{equation}

		\, . 
	\end{equation}
	The mirrored version follows analogously.
\end{proof}

Next we turn to the component $T\colon A\otimes A\longrightarrow A$ of~$\mathbb{A}_R$. 
As a 2-morphism  $A_G\boxtimes A_G\longrightarrow A_G$ in~$\mathcal{D}$ we define it to be 
\begin{equation}
	\label{eq:2-morphism of general 2group orb dat}
T:=A_H\otimes \chi:=\bigoplus_{g_1,g_2\in G}\Big(\bigoplus_{h\in H} R(g_1g_2,h)\Big)\otimes \chi_{g_1,g_2} =	
		%
\end{equation}
where~$\chi$ is part of the structure of the 3-functor~$R$, as explained in the discussion around~\eqref{eq:3functor}. 
The $(A_H,A_H\boxtimes A_H)$-bimodule structure of~$T$ is given by the multiplication~$\mu$ in~\eqref{eq:AHMuDelta} and the 3-morphism component of~$\chi$ (denoted as a black vertex) as follows:
\begin{equation}
\label{eq:general T module structure}
	%
\, .
\end{equation}
More precisely, the black vertex on the right is given by the sum over the 3-morphism components of~$\chi$, i.e.\ by the sum over
\begin{equation}
	\label{eq:chighgshs}
	\chi_{(g,h),(g^\prime,h^\prime)}\colon \chi_{t(h)g,t(h^\prime)g^\prime}\otimes \big(R(g,h)\boxtimes R(g^\prime,h^\prime)\big) 
		\longrightarrow 
		R\big(gg^\prime, h\phi_{g^{-1}}(h^\prime)\big)\otimes \chi_{g,g^\prime}\, .
\end{equation}

\begin{lemma}[label=lem:T-module structure]
	The 2-morphism~$T$ in \eqref{eq:2-morphism of general 2group orb dat} equipped with the structure in \eqref{eq:general T module structure} is an $(A_H,A_H\boxtimes A_H)$-bimodule.
	\begin{proof}
		The left $A_H$-module structure is given by multiplication of the algebra and therefore well-defined. 
		The right diagram in \eqref{eq:general T module structure} defines a right $(A_H\boxtimes A_H)$-module structure on~$T$, which follows from the naturality and unitality of $\chi$ as well as associativity of~$A_H$. 
		Naturality and unitality of~$\chi$ are given by the following identities:
	\begin{equation}
		\label{eq:naturality and unitality of chi}
\, .
	\end{equation}
	From these, the right $(A_H\boxtimes A_H)$-module axioms follow by inserting another $A_H$-line to the left and composing with~$\mu$, then applying the algebra axioms of~$A_H$.
	\end{proof}
\end{lemma}

Finally we turn to the components $\alpha,\overline\alpha$ of~$\mathbb{A}_R$, which we define as
\begin{equation}
\label{eq:3-morphism of general 2group orb dat}
	\alpha:= 
	%
}
		\, ,
\end{equation}
where the black vertex above the~$\omega^{-1}$ (or~$\omega$) denotes the sum over the inverses $\chi_{(g,h),(g^\prime,h^\prime)}^{-1}$ weighted by $1/|H|$, i.e.\ it is given by
\begin{equation}
	\bullet^{-1} = \frac{1}{|H|}\chi_{(g,h),(g^\prime,h^\prime)}^{-1} \, .
\end{equation} 
Similar to the role of the coefficient $1/|H|$ in~$\Delta$ \eqref{eq:Frob structure on AH} in showing $\Delta$-separability, the factor $1/|H|$ ensures that when a purple $A_H$-line to the left of a $T$-line passes the latter and thereby ``splits'' into two lines via~$\chi^{-1}$ and then passes back, recombining into single line via~$\chi$, this is equal to the identity 3-morphism:
\begin{equation}
\label{eq:cancelling of chi and its inverse}				
\, .
\end{equation}
The factor $1/|H|$ also ensures that when we turn the diagrams in \eqref{eq:naturality and unitality of chi} upside down, the identities still hold, as~$\varepsilon$ comes with a factor of $|H|$.

\begin{lemma}
	The 3-morphisms $\alpha,\overline{\alpha}$ in~\eqref{eq:3-morphism of general 2group orb dat} are $(A_H,A_H^{\boxtimes 3})$-bimodule morphisms in~$\mathcal{D}$. 
\end{lemma}

\begin{proof}
Note that the module structure on the domain of~$\alpha$, i.e.\ on $T\otimes_1 T=(A_H\otimes \chi)\otimes_{A_H} ((A_H\otimes \chi)\boxtimes\id_{A_G})$, is induced by \eqref{eq:general T module structure}. 
As a string diagram, we write the former as parallel lines connected by the idempotent:
 \begin{equation}
 	\label{eq:TT-idempotent}
		\tikzzbox{
}
 \end{equation}
Then the right $A_H^{\boxtimes 3}$-module structure is given by either diagram in \eqref{eq:right module structure of TT} where the equality is due to the naturality of~$\chi$ \eqref{eq:naturality and unitality of chi} and the properties of the Frobenius algebra~$A_H$:
\begin{equation}
	\label{eq:right module structure of TT}
		\tikzzbox{%
}\, .
\end{equation}
The module structure on $T\otimes_2 T$ is given similarly, exchanging front and back. 

In addition to the naturality and unitality axioms \eqref{eq:naturality and unitality of chi}, we also require the modification axiom for~$\omega$, which in the graphical calculus reads 
\begin{equation}
\label{eq:modification axiom for omega}
		\tikzzbox{%
}\, .
\end{equation}
To show that~$\alpha$ is a right $A_H^{\boxtimes 3}$-module morphism, we need to prove that the two ways of stacking the diagrams \eqref{eq:right module structure of TT} and \eqref{eq:3-morphism of general 2group orb dat} are equal. Starting with~$\alpha$ on top, the multiplications involved in the module structure \eqref{eq:right module structure of TT} can be moved from the bottom of the diagram to the top via \eqref{eq:modification axiom for omega}, and naturality of~$\chi$ and the properties of the Frobenius algebra~$A_H$ allow us to rearrange them into the desired shape, recovering \eqref{eq:right module structure of TT} on top of~$\alpha$. For the left $A_H$-module property, note that~$\mu$ and~$\Delta$ are module morphisms, hence~$\alpha$ is, too.
The analogous argument holds for~$\overline{\alpha}$.
\end{proof}

Recall again the 3-category $\Homorb{\mathcal{D}}$ from~\eqref{eq:Hom-orbData}. 

\begin{proposition}[label=prop:2group orbifold datum]
	Let~$\mathcal{G}$ be a $2$-group associated to a crossed module $(G,H,t,\phi)$, let~$\mathcal{D}$ be a Gray category with duals whose $\Hom$ 2-categories have direct sums and are $\Vect$-enriched, and let $R\colon \underline{\mathcal{G}}\longrightarrow \mathcal{D}$ be a $3$-functor whose 2-functor components on $\Hom$ 2-categories are pivotal. 
	Then $\mathbb{A}_R:= (a,A,T,\alpha,\overline{\alpha})$
	of \eqref{eq:object of general 2group orb dat}, \eqref{eq:1-morphism of general 2group orb dat}, \eqref{eq:2-morphism of general 2group orb dat}, and \eqref{eq:3-morphism of general 2group orb dat} is a weakly associative algebra in $\Homorb{\mathcal D}$.
\end{proposition}

\begin{proof}
	We have to show that the conditions~\eqref{eq:O1}, \eqref{eq:O2}, and~\eqref{eq:O3} are satisfied. 
	In each case our basic strategy is to move the purple $A_H$-lines all the way to the left using the 3-morphism components of~$\chi$ and the modification axiom for~$\omega$, before applying the appropriate identity for the $T$-lines. 
	Then, we rearrange the $A_H$-lines to fit the desired result.
	
	We start with \eqref{eq:O2}. 
	Note that since the $\otimes$-composition in the $\Hom$ 2-categories is replaced by relative tensor products in $\Homorb{\mathcal{D}}$, idempotents \eqref{eq:relative tensor idempotent} are introduced, similar to what happens in $\mathcal{D}_\mathcal{C}$, cf.\  \Cref{fig:OrbifoldDatumAxiomsInDC}, and also \Cref{eq:TT-idempotent}. 
	Thus the identity to check for~\eqref{eq:O2} is 
	\begin{equation}
		\label{eq:O2 for 2group symm}
		\tikzzbox{
}
		\, . 
	\end{equation}
	Starting on the left, note that the composition of~$\eta$ with the black unlabelled vertex (the sum over 3-morphisms components of~$\chi$, \eqref{eq:chighgshs}) reduces the sum to those terms including the unit $e_H\in H$ on the side of~$\eta$, i.e.\ $\chi_{(g,e_H),(g^\prime,h^\prime)}$ or $\chi_{(g,h),(g^\prime,e_H)}$ if $\eta$ is behind or infront of~$\chi$, respectively.
	This holds similarly for~$\varepsilon$. 
	As in the proof of \Cref{lem:Frobenius algebra structure on AH}, we can then use the invertibility of $R^0_g$ to connect the strings in the back, and the restriction of \eqref{eq:cancelling of chi and its inverse} to only those terms that include $e_H$ on one side lets us push the lines to the left. These steps are given by the following diagrams:
	\begin{equation}		
	\label{eq:cancelling of chi and its inverse with identity}
\, .
	\end{equation}
	This leaves us with a loop $\mu\circ\Delta$ on the leftmost plane which can be removed using $\Delta$-separability. 
	Cancelling~$\omega^{-1}$ and~$\omega$, we are left with two parallel $T$-lines, and we use the Frobenius property to rewrite the remaining $\Delta\circ\mu$ as $(\mu\otimes\id_{A_H})\circ(\id_{A_H}\otimes \Delta)$. 
	Next we insert two black vertices on the upper part of the left line by pushing the string connecting~$\mu$ and~$\Delta$ through~$T$ above the already existing vertices (going from right to left in \eqref{eq:cancelling of chi and its inverse}). 
	Then we use the naturality of~$\chi^{-1}$, i.e.\ \eqref{eq:naturality and unitality of chi} turned upside-down, to move the comultiplication to the right side of the $T$-line. 
	One of the comultiplications cancels immediately, since it ends in a counit, the other connects the upper right purple $A_H$-line to the two lines along~$T$. 
	We can now cancel the bottom two black vertices, leaving us with a diagram that is like the one on the right-hand side of~\eqref{eq:O2 for 2group symm}, but without the full idempotent. 
	To remedy this, observe that $\Delta=(\id\otimes \mu)\circ ((\Delta\circ\eta)\otimes \id)$ by the Frobenius property and unitality. 
	This finishes the proof of~\eqref{eq:O2}, and~\eqref{eq:O3} is checked similarly. 
	
	We now turn to~\eqref{eq:O1}, so we have to check the identity 
	\begin{equation}
		\label{eq:O1ForCrazyPeople}
		\tikzzbox{
}
		\, . 
	\end{equation}
	We start on the right-hand side with the immediate goal of moving the purple $A_H$-line which connects the two black vertices across the central $\omega^{-1}$. 
	Locally, we take the following steps:
	\begin{equation}
		%
	 \stackrel{\eqref{eq:modification axiom for omega}}{=} |H| \cdot 
		%
	 \stackrel{\eqref{eq:cancelling of chi and its inverse with identity}}{=}
		%
	\, .
	\end{equation}	
	Here, we extend the lines connecting the two black vertices to units and counits downwards and upwards, respectively.
	We then push the $A_H$-line on the front plane through the upper $\chi$-line using the 3-morphism component of~$\chi$ and its inverse. Next, we apply the modification axiom~\eqref{eq:modification axiom for omega} in order to move the lines coming from the lower black vertex and the unit to the bottom of the central $\omega^{-1}$. Lastly, we cancel two of the black vertices on the bottom as in \eqref{eq:cancelling of chi and its inverse with identity}.
	On the right-hand side of \Cref{eq:O1ForCrazyPeople}, this leaves us with the three copies of $\omega^{-1}$ connected by $\chi$-lines without additional $\chi$-insertions. 
	We can then apply the pentagon axiom, which $\omega^{-1}$ satisfies since it is part of a 3-functor, resulting in the configuration of $\chi$-lines of the left-hand side of~\eqref{eq:O1ForCrazyPeople}. 
	The resulting diagram is almost exactly the one on the left-hand side, apart from the fact that the $A_H$-lines coming from the very right are now passing through the $\chi$-lines 
	below and above the $\omega^{-1}$, thus without intersecting their connection:	
	\begin{equation}
		\label{eq:IntermediateO1}
				\tikzzbox{
}
	\end{equation}
	
	Now the bottom of the diagram in~\eqref{eq:IntermediateO1} bundles $A_H$-lines to the left by pairing them with units and passing them through $T$-lines via 3-components of~$\chi$, and the top does the same in reverse, using the inverse 3-morphism components of~$\chi$ and counits. 
	By naturality of~$\chi$ and the Frobenius algebra properties of~$A_H$, we can rearrange the remaining morphisms of $A_H$-lines on the left-hand side of~\eqref{eq:O1ForCrazyPeople} into a sequence of multiplications followed by a sequence of comultiplications. 
	Using the modification axiom again, we can pass the $A_H$-lines coming from the very right back across the two $\omega^{-1}$, and the unitality condition~\eqref{eq:naturality and unitality of chi} for~$\chi$ lets us reduce the units and counits. 
	By similar methods as before, i.e.\ pushing $A_H$-lines across $\chi$-lines and using naturality and unitality of~$\chi$ to move the two 3-morphism components of~$\chi$ and~$\chi^{-1}$ past each other, we obtain the desired result.
\end{proof}

We remark that the algebra~$\mathbb{A}_R$ can be enhanced by a unital structure constructed from the unital coherence data of~$R$. 
This turns~$\mathbb{A}_R$ into an $E_1$-algebra, as discussed and used in this setting in \cite{Carqueville2023}.

Conceptually, $\mathbb{A}_R$ is the ``2-group algebra'' associated to the $2$-group~$\mathcal{G}$. 
If~$\mathcal{G}^\oplus$ denotes the completion of~$\mathcal G$ with respect to direct sums, and $\textrm{B}_{\C}\underline{\mathcal{G}^\oplus}$ denotes the $\Vect$-enrichment of $\textrm{B}\underline{\mathcal{G}^\oplus}$ at the level of its categories of 2- and 3-morphisms, then the 2-group algebra is an orbifold datum inside $\Homorb{\textrm{B}_{\C}\underline{\mathcal{G}^\oplus}}$ itself. 
Our construction above then corresponds to the transport of this orbifold datum along~$R$ to $\Homorb{\mathcal{D}}$. 
We expect this transport to provide an orbifold datum iff~$R$ is pivotal, in an appropriate 3-dimensional sense. 

In two dimensions, we would start with a mere group~$G$ instead of a 2-group~$\mathcal G$. 
Then the analogous completion is to use the category $G$-$\Vect$ of $G$-graded vector spaces instead of~$\underline{G}$. 
In this case it is known \cite[Ex.\,3.3]{Turaev2010} that the group algebra (understood as the sum over $g\in G$ of copies of~$\C$ in degree~$g$) carries the structure of a $\Delta$-separable symmetric Frobenius algebra. 
We describe this structure in detail in \Cref{sec:Deeq is 1-form orbifolding}. 
Its transport along a pivotal monoidal functor gives an orbifold datum in the codomain (cf.\ \cite[Sect.\,4.2]{Mulevicius2022a}):

\begin{proposition}
	Let $F\colon \mathcal{B}\longrightarrow\mathcal{B}^\prime$ be a pivotal (strong) monoidal 2-functor between pivotal 2-categories. 
	Then~$F$ preserves orbifold data.
\end{proposition}

\begin{proof}
	Let $(a,A)$ be an orbifold datum in~$\mathcal{B}$, i.e.\  a $\Delta$-separable symmetric Frobenius algebra which consists of an object~$a$ in~$\mathcal{B}$ and 1- and 2-morphisms $(A,\mu,\Delta,\eta,\varepsilon)$ in $\End_\mathcal{B}(a)$ as in \Cref{subsubsec:2dOrbDat}. 
	Recall from \cite{DayPastro2008,Mulevicius2022a} that its transport along~$F$ is $(F(A),\mu_F,\Delta_F,\eta_F,\varepsilon_F)$, given by
	\begin{equation}
		\label{eq:Frobenius transport}
		\mu_F:= F(\mu)\circ F^2_{A,A}\,,\quad  \Delta_F:=(F^2_{A,A})^{-1}\circ F(\Delta)\,,\quad \eta_F:= F(\eta)\circ F^0_A\,,\quad \varepsilon_F:=(F^0_A)^{-1}\circ F(\varepsilon)\,.
	\end{equation}
	The proof that this is a Frobenius algebra is analogous to that of \Cref{lem:Frobenius algebra structure on AH}, additionally using the Frobenius algebra axioms of $(A,\mu,\Delta,\eta,\varepsilon)$. 
	$\Delta$-separability follows from the fact that~$F$ is strong monoidal, symmetry follows from the fact that~$F$ is pivotal.
\end{proof}

\begin{remark}
	\label{rem:HigherGaugableGroupSymmetries}
	In our construction in \Cref{prop:2group orbifold datum} of the 3-dimensional candidate orbifold datum~$\mathbb{A}_R$ we used that the 1-form symmetry ``within~$R$'' is gaugeable (because of pivotality of the 2-functor components, see \Cref{prop:Orbifold datum from 1form symmetry as composite} below) and then produced the surface defect of~$\mathbb{A}_R$ by orbifolding the 1-form symmetry on the surface defect~$A_G$ where the latter is given by the sum over 0-form symmetry defects. 
	We expect that an analogous construction should provide a 4-dimensional candidate orbifold datum from a finite 3-group symmetry $\mathcal R\colon \mathrm{B}\mathfrak{G}\longrightarrow \mathcal{D}$, where~$\mathcal D$ is an appropriate 4-category, and the 3-group~$\mathfrak G$ in particular involves three groups $G_0,G_1,G_2$ (analogous to a 2-group involving two groups $G,H$). 
	
	Indeed, assuming that the 3-functor components~$R$ of the 4-functor~$\mathcal R$ are pivotal, we expect the 3-dimensional candidate orbifold datum $\mathbb{A}_R$ built from the 1- and 2-form symmetries ``within~$\mathcal R$'' to be a genuine one (as discussed above), so we can gauge it on the 3-dimensional defect~$A_{G_0}$ of the 0-form symmetry within~$\mathcal R$. 
	The result is a 1-endomorphism in the 4-category $\Homorb{\mathcal D}$ which is defined in analogy to the 3-dimensional case (and using the results of \cite{Carqueville2023} to orbifold the Hom 3-categories) and we use it as the 3-dimensional defect of the 4-dimensional candidate orbifold datum $\mathbb{A}_\mathcal{R}$.
	The surface defect connecting three copies of this 3-dimensional defect should come from the coherence data of $\mathcal{R}$ and the appropriate product with the surface defect $(A_{G_1},A_{G_2})$ of~$\mathbb{A}_{\mathcal R}$, similar to~\eqref{eq:2-morphism of general 2group orb dat}. In this setting, the surface defect is the gauging of the 2-form symmetry $G_2$ on the sum of 1-form symmetry defects $A_{G_1}$.
	The line and point defects of~$\mathbb{A}_{\mathcal R}$ should follow similarly in principle.
	We expect that this produces a genuine orbifold datum if the 4-functor~$\mathcal{R}$ is pivotal. 
	More generally, we expect that $n$-dimensional orbifold data can be iteratively constructed from pivotal $(n-1)$-group symmetries for any $n\geqslant 2$.
\end{remark}

\subsection{0-form Symmetries}
\label{sec:0-form specialisation}

For any group~$G$, there is a crossed module $(G,H,t,\phi)$ where~$H,t,\phi$ are each trivial. 
Specialising \Cref{def:3d 2-group symmetry} to this case, we have a \textsl{0-form $G$-symmetry}. As a 3-functor, this takes the form 
\begin{equation}
	\label{eq:0form symm in 3d}
	R\colon \mathrm{B}\underline{G}\longrightarrow \mathcal{D}
\end{equation}
where the underline now denotes the addition of identity 2- and 3-morphisms.
Note that in this case, the algebra~$A_H$ defined in~\eqref{eq:1-morphism of general 2group orb dat}, \eqref{eq:AHMuDelta} and~\eqref{eq:Frob structure on AH} is a sum of identity 2-morphisms, with every summand $R(g,e_H)$ carrying its own Frobenius algebra structure via $\mu_g:=R^2_{(g,e_H),(g,e_H)}$ and $\eta:=R^0_g$, and~$\Delta$ and~$\varepsilon$ given by the respective inverses. 

It is a straightforward calculation to check that 
\begin{equation}
	R^0_g\colon \id_{R(g)}\stackrel{\cong}{\longrightarrow} R(g,e_H)
\end{equation}
is an isomorphism of Frobenius algebras for any 3-functor~$R$, where the Frobenius structure on $\id_{R(g)}$ is given by unitors and identities. 
We apply this observation to the associative algebra constructed in \Cref{prop:2group orbifold datum}, concluding that we can identify the 1-morphism~$A$ in $\Homorb{\mathcal D}$ with $(A_G,\bigoplus_{g\in G}\id_{R(g)})$, which is the trivial Frobenius algebra on~$A_G$. Similarly, $R^0$ provides a 3-isomorphism between~$T$ and $\id_{R(g)}\otimes\chi$ which has an $\id_{A_G}$-module structure. 
Thus, the algebra~$\mathbb{A}_R$ of \Cref{prop:2group orbifold datum} lies in the (essential) image of the inclusion functor $\mathcal{D}\longrightarrow \Homorb{\mathcal D}$, and we can pull it back along this functor, forgetting the Frobenius algebra and the module structure entirely. 
The resulting associative algebra is in~$\mathcal{D}$ (as opposed to $\Homorb{\mathcal D}$) and we denote it by~$\mathbb{A}_R$ as well, given as follows:
\begin{gather}
	a:= R(\ast)
	\, , \quad 
	A:=\bigoplus_{g\in G}R(g)
	\, , \quad 
	T:=\bigoplus_{g_1,g_2\in G}\chi_{g_1,g_2} 
	\, , 
	\nonumber 
	\\ 
	\alpha:=
}
			\, .
	\label{eq:0form orbdat in 3d}
\end{gather}

\subsection{1-form Symmetries}
\label{sec:1-form specialisation}

We now consider the other edge case where in a crossed module $(G,H,t,\phi)$ all of $G,t,\phi$ are trivial, and we are left only with the group~$H$. 
Due to~\eqref{eq:Peiffer identity}, $H$ is necessarily commutative. 
Put differently, we consider a pure \textsl{1-form symmetry} based on~$H$. 
The functor \eqref{eq:3d 2-group symmetry} then takes the form 
\begin{equation}
	R\colon \mathrm{B}^2\underline{H}\longrightarrow \mathcal{D}\,.
\end{equation} 
For the pivotal braided category of $2$-endomorphisms of $R(\id_\ast)$ we write  
\begin{equation}
	\mathcal{C} := \textrm{End}_{\mathcal D}\big(R(\id_*)\big) 
\end{equation} 
where~$*$ is the unique object in $\mathrm{B}^2\underline{H}$. 
The category~$\mathcal C$ inherits its pivotal braided monoidal structure from~$\mathcal D$ (cf.\ \cite[Lem.\,2.18]{Barrett2012}), while~$\underline{H}$ is naturally equipped with such a structure because~$H$ is commutative. 
Using this, one finds that (up to strictification)~$R$ descends to and is equivalently described by a braided monoidal functor 
\begin{equation}
	\label{eq:Ras1functor}
	\overline{R}\colon \underline{H} \longrightarrow \mathcal{C} \, .
\end{equation}

In this section, we construct a new orbifold datum $\mathbb{A}_{A_{\overline{R}}}$ from the functor~$\overline{R}$ using methods similar to those of \Cref{subsec:GeneralCase} and \cite{Carqueville2020}. 
We show in \Cref{prop:Orbifold datum from 1form symmetry as composite} that the construction reproduces the orbifold datum~$\mathbb{A}_R$ from the general construction in \Cref{prop:2group orbifold datum}.
The new presentation of the orbifold datum allows for various new tools in calculating orbifold data, such as \Cref{prop:orbifold data for compositions of functors}.

\medskip 

\begin{proposition}[label=prop:FrobeniusAlgebra from 1form symmetry]
	Given a braided monoidal functor $\overline{R}\colon \underline{H}\longrightarrow \mathcal{C}$, we set
	\begin{gather}
	\begin{alignat}{2}
		A_{\overline{R}}&:= \;\,\bigoplus_{h\in H}\;\, \overline{R}(h)\,,&&
		\\
		\mu&:=\bigoplus_{h,h^\prime\in H}\overline{R}^2_{h,h^\prime} \colon & &A_{\overline{R}}\otimes A_{\overline{R}}\longrightarrow A_{\overline{R}}\,,
		\\
		\Delta&:=\bigoplus_{h,h^\prime\in H}\frac{1}{|H|}\Big(\overline{R}^2_{h,h^\prime}\Big)^{-1}\!\colon  &&A_{\overline{R}}\longrightarrow A_{\overline{R}}\otimes A_{\overline{R}}\,,
		\\
		\eta&:=\mathrm{incl}_{\overline{R}(e)}\circ \overline{R}^0\colon &&\mathds{1}\longrightarrow \overline{R}(e)\longrightarrow A_{\overline{R}}\,,
		\\
		\varepsilon&:=|H| \cdot (\overline{R}^0)^{-1}\circ \pr_{\overline{R}(e)}\colon &&A_{\overline{R}}\longrightarrow \overline{R}(e)\longrightarrow\mathds{1}\,.
	\end{alignat}
	\end{gather}
	Then $(A_{\overline{R}},\mu,\Delta,\eta,\varepsilon)$ is a commutative $\Delta$-separable Frobenius algebra. 
	If~$\overline{R}$ is pivotal, then~$A_{\overline{R}}$ is symmetric.
\end{proposition}

\begin{proof}
	With the exception of commutativity, the proof is the same as that of \Cref{lem:Frobenius algebra structure on AH}.
	To show commutativity, we note that the braiding on~$\underline{H}$ is trivial since~$H$ is commutative. 
	Therefore, since~$\overline{R}$ is braided, the braiding on $A_{\overline{R}}\otimes A_{\overline{R}}$ (in~$\mathcal{C}$) is given by 
	\begin{equation}
		\label{eq:braiding on A}
		c_{A_{\overline{R}},A_{\overline{R}}}=\bigoplus _{h,h^\prime\in H}\Big(\overline{R}^2_{h,h^\prime}\Big)^{-1}\circ \overline{R}^2_{h^\prime,h}\,. 
	\end{equation}
	Note that this is different from $\Delta\circ\mu$ which includes additional terms $(\overline{R}^2_{h_1,h_1^\prime})^{-1} \circ \overline{R}^2_{h_2^\prime,h_2}$ for any two pairs of group elements satisfying $h_1h_1^\prime=h_2 h^\prime_2$. 
	When we compose $c_{A_{\overline{R}},A_{\overline{R}}}$ with~$\mu$, terms cancel in such a way that the result is again~$\mu$. 
	Hence~$A_{\overline{R}}$ is commutative.
\end{proof}

Recall the 3-category~$\mathcal D_{\mathcal C}$ that we associated to a ribbon category~$\mathcal C$ in \Cref{exa:3dCategoriesWithAdjoints}. 
Every commutative $\Delta$-separable symmetric Frobenius algebra in~$\mathcal C$ can be turned into an orbifold datum: 

\begin{proposition}[name={\cite[Prop.~3.15]{Carqueville2020}}, label=prop:orbifold datum from Frobenius algebra] 
	Let~$\mathcal{C}$ be a braided fusion category and let $(A,\mu,\Delta,\eta,\varepsilon)$ be a commutative $\Delta$-separable symmetric Frobenius algebra in~$\mathcal{C}$. 
	Then 
	\begin{equation}
	\label{eq:orbdat from cssFrob}
		\mathbb{A}_A 
			:=
			\big(\mathcal{C}, \, A, \, T=A, \, \alpha=\overline{\alpha}=\Delta\circ\mu, \,  \psi=\eta, \, \phi=1\big)
	\end{equation}
	is a (special) orbifold datum in~$\mathcal{D}_\mathcal{C}$, where the $(A,A\otimes A)$-bimodule structure of $T$ is given by various applications of $\mu$.
\end{proposition}

Since $\Delta$-separable Frobenius algebras are in particular separable, we have $\mathrm{B}(\mathrm{B} \mathcal C)_\text{orb}\subset D_{\mathcal C}$. 
We want to generalise \Cref{prop:orbifold datum from Frobenius algebra} to any 3-category which is equivalent to a Gray category with duals, using the category of 2-endomorphisms of the identity 1-morphism on some object~$a$ as~$\mathcal{C}$, i.e.\ $\mathcal{C}=\End_\mathcal{D}(\id_a)$. 
In this setting, we have (recall the discussion around~\eqref{eq:DHomorb})
\begin{equation}
	\mathrm{B}(\mathrm{B} \mathcal{C})_\text{orb}\longhookrightarrow \Homorb{\mathcal D}\,,
\end{equation} 
and thus the orbifold datum constructed in \eqref{eq:orbdat from cssFrob} is an orbifold datum in $\Homorb{\mathcal D}$.

\medskip

\begin{proposition}[label=prop:Orbifold datum from 1form symmetry as composite]
	The orbifold datum~$\mathbb{A}_R$ associated to a $1$-form symmetry $R\colon \mathrm{B}^2\underline{H} \longrightarrow \mathcal{D}$ coincides with the one constructed by combining \Cref{prop:FrobeniusAlgebra from 1form symmetry} and \Cref{prop:orbifold datum from Frobenius algebra},
	\begin{equation}
		\mathbb{A}_R=\mathbb{A}_{A_{\overline{R}}}\,.
	\end{equation}
	In particular, if a 1-form symmetry~$R$ has pivotal 2-functor components on $\Hom$ 2-categories, then \Cref{prop:2group orbifold datum} yields an orbifold datum (and not just an algebra).
\end{proposition}

\begin{proof}
	We assume that $R=(R,\chi,\iota,\omega,\gamma,\delta)$ is sufficiently strictified, such that $R(\id_\ast)=\id_{R(\ast)}$, and the unital data $\iota,\gamma,\delta$ from \Cref{subsubsec:3cats} are identities. 
	Then the 2-morphism components of~$\chi$ are given by unitors in~$\mathcal{D}_\mathcal{C}$. 
	From the fact that $\overline{R}{}^2=R^2$ and $\overline{R}{}^0=R^0$ where applicable, it follows that the algebras~$A_{\overline{R}}$ from \Cref{prop:FrobeniusAlgebra from 1form symmetry} and~$A_H$ defined in~\eqref{eq:1-morphism of general 2group orb dat} and~\eqref{eq:Frob structure on AH} are identical.
	
	Next we show that the orbifold datum of \Cref{prop:2group orbifold datum} reduces to the one built from~$A_H$ via \Cref{prop:orbifold datum from Frobenius algebra}.
	To this end, we can identify $\boxtimes$-products and $\otimes$-products between $2$-morphisms in $\End_{\mathcal D_{\mathcal C}}(\id_a)$ by an Eckmann--Hilton argument. 
	More precisely, we have a 3-isomorphism
	\begin{equation}
		r_{\id_a}\otimes (f\boxtimes g)\cong (f\otimes g)\otimes r_{\id_a}\,,
	\end{equation}
	given by unitors and the coherence morphism $\boxtimes^2$. 
	In a Gray category, the right unitors in this equation are identities and the 3-isomorphism reduces to the tensorator. 
	Under this identification, we can show that the 3-morphism components of~$\chi$ are given by~$R^2$, by starting with~\eqref{eq:naturality and unitality of chi} and applying invertibility of~$R^2$. 
	Using these identifications, it follows that $T=A_H=A_{\overline{R}}$ as objects and the module structure on~$T$ is given by the multiplication~$\mu$ on~$A_H$. 
	Therefore, the identity also holds on the level of modules. 
	Similarly, the bimodule morphisms~$\alpha$ and~$\overline{\alpha}$ of~\eqref{eq:3-morphism of general 2group orb dat} reduce to $\Delta\circ \mu$ by the (co)unitality axioms of $A_H$.
\end{proof}

\medskip

Recall from \cite[Def.\,5.12]{Mulevicius2022a} the notion of transport of an orbifold datum~$\mathbb{A}$ in a Gray category with duals of the type~$\mathcal{D}_{\mathcal{C}_1}$ along a braided strong monoidal functor~$F\colon\mathcal{C}_1\longrightarrow\mathcal{C}_2$. 
This produces an orbifold datum $F(\mathbb{A})$ in~$\mathcal{D}_{\mathcal{C}_2}$, given by
\begin{equation}
\label{eq:orbifold transport in DC}
 	F(\mathbb{A}):= \big(F(A),F(T),\alpha_F,\overline{\alpha}_F, \psi_F,\phi_F\big)
\end{equation}
where $F(A)$ is given by \eqref{eq:Frobenius transport}, $F(T)$ is similarly given by the underlying object with module structure composed with (the sum over)~$F^2$, and~$\alpha$ and~$\overline{\alpha}$ are also composed with~$F^2$ and its inverse. 
There is an additional assumption of compatibility between~$F$ and~$\mathbb{A}$ which ensures that suitable~$\psi_F$ and~$\phi_F$ exist, see \cite[Def.\,5.10]{Mulevicius2022a} for details.
 
As with \Cref{prop:orbifold datum from Frobenius algebra}, there is a direct generalisation of this construction to the setting of the present section. 
Orbifold data from 1-form symmetries are compatible with transport along such functors in the following sense.

\begin{proposition}[label=prop:orbifold data for compositions of functors]
	Let $\overline{R}\colon \underline{H}\to\mathcal{C}$ and $F\colon\mathcal{C}\to\mathcal C'$ be ribbon functors.
	\begin{enumerate}
	\item
	Given any commutative $\Delta$-separable symmetric Frobenius algebra~$A$ in~$\mathcal{C}$, the transport of the associated orbifold datum \eqref{eq:orbdat from cssFrob} (with $\psi_F=\eta_F$ and $\phi_F=1$) is the orbifold datum associated to $F(A)$:
	\begin{equation}
		F(\mathbb{A}_A)=\mathbb{A}_{F(A)}.
	\end{equation}
	\item 
	The orbifold data associated to 1-form symmetries~$\overline{R}$ and $F\circ \overline{R}$ are related by transport along~$F$:
	\begin{equation}
		F(\mathbb{A}_{A_{\overline{R}}})=\mathbb{A}_{F(A_{\overline{R}})}=\mathbb{A}_{A_{F\circ \overline{R}}}\,.
	\end{equation}
	\end{enumerate}
\end{proposition}

\begin{proof}
	The first statement is obvious by comparing definitions \eqref{eq:Frobenius transport}, \eqref{eq:orbifold transport in DC}, and \Cref{prop:orbifold datum from Frobenius algebra}: the underlying algebra is the same by definition, the module structure of $F(T)=F(A)$ is given by $F(\mu)\circ F^2_{A,A}$ in both cases, and similarly $\alpha_F=\Delta_F\circ\mu_F$. 
	The second statement then follows from the observation that $F(A_{\overline{R}})=A_{F\circ \overline{R}}$ which follows from unravelling the definition of composition of monoidal functors.
\end{proof}

\section{Symmetries from $\boldsymbol{G}$-crossed Braided Fusion Categories}
\label{sec:SymmetriesFromGCrossedBraidedCategories}

In this section we discuss a specific class of categories (of line defects) where we can gauge a finite group-like symmetry. 
To this end, we consider $G$-crossed braided fusion categories which in particular carry a $G$-action for some finite group~$G$.
We gauge this $G$-action in two ways: 
On the one hand, given a $G$-crossed braided fusion category we can directly apply the standard equivariantisation construction (which we recall in \Cref{sec:General theory of (De)Equivariantisation}), while on the other hand we can also produce a 3-functor (a $G$-symmetry) and its associated orbifold datum from it (\Cref{sec:equivariantisation is orbifolding}) via the methods of Section~\ref{sec:0-form specialisation}. 
The latter recovers the orbifold datum of \cite[Thm.\,5.1]{Carqueville2020}.
In \Cref{sec:equivalence of equivariantisation and orbifolding} we prove that both paths lead to equivalent results. 

\Cref{sec:Deeq is 1-form orbifolding} is dedicated to undoing the gauging of the $G$-symmetry by gauging an emergent dual symmetry. 
This is called de-equivariantisation, and can also be understood as gauging a $1$-form symmetry. 
When~$G$ is non-commutative, this symmetry is non-invertible and thus falls outside our previous discussion. 
However, we identify an orbifold datum that achieves the de-equivariantisation in either case, leading to a clear and general picture.

\subsection{Equivariantisation and De-equivariantisation}
\label{sec:General theory of (De)Equivariantisation}

In this section we review group-like actions on (braided) monoidal categories and the related constructions of equivariantisation and de-equivariantisation. 
Our main application is in terms of $G$-crossed braided fusion categories, for which we also collect the relevant details. 
Our main references for this section are \cite[Sect.\,4]{Drinfeld2009} and \cite[Sect.\,7]{etingof2009fusioncategorieshomotopytheory}. 

\medskip 

Given monoidal categories $V$ and $\mathcal{C}$, recall that an \textsl{action of~$V$ on~$\mathcal{C}$} is a monoidal functor $V\longrightarrow \End^\otimes(\mathcal{C})$ into the category of monoidal endofunctors. 
Equivalently, this is a functor 
\begin{equation}
	\label{eq:module category action}
	V\times \mathcal{C}\longrightarrow\mathcal{C} 
\end{equation}
equipped with additional associativity and unitality coherence isomorphisms, subject to the usual pentagon and triangle conditions, which involve the associators and unitors of both the action as well as those in~$V$; see e.g.\ \cite{Ostrik2006} for details.

A special case is a \textsl{$G$-action on a monoidal category~$\mathcal C$} for some group~$G$, by which we mean a monoidal functor\footnote{By the superscript ``rev'' we indicate the reversal of the monoidal product. This means that we use right $G$-actions in accordance with \cite[Sect.\,5]{Carqueville2020}.}
\begin{equation}
	\rho\colon \underline{G}^\text{rev}\longrightarrow \Aut^{\otimes}(\mathcal{C}) \, . 
\end{equation}
Here~$\underline{G}$ denotes the monoidal category whose objects are group elements $g\in G$, the tensor product is group multiplication, and~$\underline{G}$ has only identity morphisms. The codomain $\Aut^{\otimes}(\mathcal{C})$ is the category of monoidal autoequivalences.
For every $g\in G$, we abbreviate $\rho_g := \rho(g)$ and we denote its monoidal structure isomorphisms as
\begin{equation}
	\rho_g^2\colon\! \otimes\circ\, (\rho_g\times\rho_g) \Longrightarrow \rho_g\circ \otimes
		\, , \quad 
	\rho^0_g\colon \mathds{1}\longrightarrow \rho_g(\mathds{1}) \, .
\end{equation}
Since~$\rho$ itself is also monoidal, it also comes with structure isomorphisms given by invertible natural transformations 
\begin{equation}
	\rho^2\colon \rho_{(-)} \circ \rho_{(-)} \Longrightarrow\rho_{(-)\cdot(-)}
		\, , \quad 
	\rho^0\colon \id_\mathcal{C}\Longrightarrow \rho_e \, ,
\end{equation}
where~$\rho^0$ is monoidal and~$\rho^2$ has components which are themselves monoidal natural transformations, 
\begin{equation}
	\rho^2_{g,h}\colon \rho_g\circ\rho_h\Longrightarrow\rho_{hg} \, .
\end{equation}
Note that all of (the component morphisms of) these natural transformations can be distinguished by their sub- and superscripts.

From now on we assume that the group~$G$ is finite. 

\begin{definition}[label=def:equivariantisation]
	Let~$\mathcal{C}$ be a fusion category, and let $\rho\colon \underline{G}^\text{rev}\longrightarrow \Aut^\otimes(\mathcal{C})$ be a monoidal functor.
	\begin{enumerate}
		\item 
		A \textsl{$G$-equivariant object} in~$\mathcal{C}$ is a pair $(X, u:= \lbrace u_g\rbrace_{g\in G})$ with $X \in \mathcal{C}$ and a family of isomorphisms 
		\begin{equation}
			u_g\colon \rho_g(X)\stackrel{\cong}{\longrightarrow} X
		\end{equation}
		such that 
		the following diagram commutes:
		\begin{equation}
		\label{eq:defining property of equivariant structure}
			\begin{tikzcd}[column sep= 40,row sep= 30, ampersand replacement=\&]
			\rho_g(\rho_h(X))
			\arrow[d, "{(\rho^2_{g,h})_X}"{left}]
			\arrow[r, "{\rho_g(u_h)}"]
			\& \rho_g(X)
			\arrow[d, "{u_g}"]\\
			\rho_{hg}(X)
			\arrow[r, "u_{hg}"{below}]
			\& X			
			\end{tikzcd}
		\end{equation}
		\item 
		The \textsl{equivariantisation} of~$\mathcal{C}$ (with respect to~$\rho$) is the category~$\mathcal{C}^G$ whose objects are $G$-equivariant objects in $\mathcal{C}$ and whose morphisms $(X,u)\longrightarrow (Y,v)$ are morphisms $f\colon X\longrightarrow Y$ in~$\mathcal{C}$ such that
		\begin{equation}
		\label{eq:defining property of morphisms in equivariantisation}
		\begin{tikzcd}[column sep= 40,row sep= 30, ampersand replacement=\&]
			\rho_g(X)
			\arrow[d, "{u_g}"{left}]
			\arrow[r, "{\rho_g(f)}"]
			\& \rho_g(Y)
			\arrow[d, "{v_g}"]\\
			X
			\arrow[r, "f"{below}]
			\& Y
		\end{tikzcd} 
		\end{equation}
		commutes for all $g\in G$. 
	\end{enumerate}
\end{definition}

Note that the above definition of~$\mathcal C^G$ as a mere category does not use the fact that~$G$ acts through \textsl{monoidal} automorphisms $(\rho_g,\rho_g^2,\rho_g^0)$. 
However, we can use this additional structure to give a monoidal structure on the equivariantisation $\mathcal{C}^G$. 
For further details see e.g.\ \cite[Sect.\,4]{Drinfeld2009} or \cite[App.\,5]{Turaev2010}:

\begin{lemma}
	The equivariantisation~$\mathcal{C}^G$ has a monoidal structure with $(X,u)\otimes (Y,v):=(X\otimes Y, w)$, where
	\begin{equation}
		w_g := \Big( \!
		\begin{tikzcd}[column sep= 40,row sep= 30, ampersand replacement=\&]
			\rho_g(X\otimes Y)
			\arrow[r, "(\rho_g^2)^{-1}_{X,Y}"]
			\& 
			\rho_g(X)\otimes\rho_g(Y)
			\arrow[r, "u_g\otimes v_g"]
			\&
			X\otimes Y 
		\end{tikzcd} 
		\!\Big) , 
	\end{equation}
	and the monoidal unit is $(\mathds{1},\lbrace (\rho^0_g)^{-1}\colon \rho_g(\mathds{1})\longrightarrow \mathds{1}\rbrace_{g\in G})$.
\end{lemma} 

We briefly pause to outline the interpretation of~$\mathcal C^G$, how it fits into the structure of the remainder of this paper, and how it connects with our motivation from topological quantum field theory. 
In the context of Reshetikhin--Turaev theory, the underlying category~$\mathcal{C}$ is interpreted as the category of line defects. 
While local operators in the gauged theory have to be invariant under the $G$-action, line defects may have non-trivial junctions between them. 
Therefore, instead of demanding $\rho_g(X)=X$ for a line operator $X\in\mathcal{C}$, we can impose $\rho_g(X)\cong X$ and choose a set of such isomorphisms (junctions) compatible with group multiplication.
Hence, the equivariantisation is the category of line defects in the gauged theory. 
One should also include twisted sectors, which we do in \Cref{def:G-crossed braided category}.
Assuming $\rho_g(\mathds{1})=\mathds{1}$ and $\rho^0_g=\id_\mathds{1}$ for the moment, we also recover the invariant local operators as endomorphisms of the ``transparent'' line defect $\mathds{1}$, i.e.\ $(\End_\mathcal{C}(\mathds{1}))^G=\End_{\mathcal{C}^G}(\mathds{1},\id_\mathds{1})$.
The physics literature then provides arguments that this new theory features ``Wilson lines'' that implement a potentially non-invertible symmetry, see e.g.\ \cite{Gaiotto_2015, furlan2024finitegroupglobalgauged}, and whose gauging undoes the gauging of the 0-form symmetry~$R$.
We turn to describe it as a $\Rep(G)$-action on~$\mathcal C^G$, and we will discuss its gauging in \Cref{sec:Deeq is 1-form orbifolding}. 

\medskip 
 
We start by noting that if we equip $\Vect$ with the trivial $G$-action, then its equivariantisation is simply the category of finite-dimensional $G$-representations, $\Vect^G = \Rep(G)$. 
Next, observe that any fusion category~$\mathcal C$ has an action \eqref{eq:module category action} of $\Vect$, where for a vector space~$V$ and an object $X\in \mathcal C$, the object $V\otimes X \in \mathcal C$ is the one which represents the functor $V\otimes_\C \Hom_\mathcal{C}(X,-)$, and we take tensor products on morphisms.
A choice of basis in~$V$ induces the identification $V\otimes X\cong X^{\oplus \mathrm{dim}(V)}$. If we restrict to the action on the unit object in~$\mathcal{C}$, this induces an inclusion $\Vect\longrightarrow \mathcal{C}$ as a subcategory (which is full if the unit object is simple). We can equivariantise this action to obtain the functor 
\begin{align}
	\Rep(G)\times \mathcal{C}^G &\longrightarrow \mathcal{C}^G
	\nonumber 
	\\
	\big((V,\phi),(X,u)\big)&\longmapsto \big(V\otimes X, \lbrace \phi(g)\otimes u_g\rbrace_{g\in G}\big) \,,
	\label{eq:RepGCGCG}
\end{align}
which describes precisely a $\Rep(G)$-symmetry. Here, $(V,\phi)$ denotes a representation with underlying vector space~$V$ and right $G$-action $\phi\colon G^\text{rev}\longrightarrow \Aut(V)$.

By restricting the functor~\eqref{eq:RepGCGCG} to the unit in the second argument, we obtain  an inclusion 
\begin{align}
	\label{def:inclusion of RepG}
	I\colon \Rep(G)\longrightarrow \mathcal{C}^G
\end{align}
of $\Rep(G)$ into~$\mathcal C^G$ as a full monoidal subcategory.
This functor allows us to undo the equivariantisation, as we recall in~\eqref{eq:eq-deeq equivalence} below. 
For this, first note that the $\C$-algebra $\C(G)$ of functions on~$G$ can be naturally lifted to an algebra in $\Rep(G)$ through the action 
\begin{equation}
	\label{eq:right action on alg of functions}
	\big(\phi^\textrm{r}_g(f)\big)(h):=f(hg^{-1}) \, ,
\end{equation} 
where $g,h\in G$ and $f\in \C(G)$. Recall that an algebra in $\Rep(G)$ is a $\C$-algebra with $G$-action whose multiplication is $G$-equivariant. 
This holds for $\phi^\textrm{r}$ because 
\begin{equation}
	\phi^\textrm{r}_g(f_1f_2)(h)=(f_1f_2)(hg^{-1})=\phi^\textrm{r}_g(f_1)(h)\phi^\textrm{r}_g(f_2)(h) \, .
\end{equation}

Now let~$\mathcal B$ be a monoidal category, and let $J\colon\Rep(G)\longrightarrow Z(\mathcal{B})$ be a braided monoidal functor. 
Recall that the Drinfeld centre $Z(\mathcal{B})$ is the category whose objects are pairs $(X,\theta)$ consisting of an object $X\in\mathcal B$ and a half-braiding $\theta\colon (-)\otimes X \longrightarrow X \otimes (-)$, and morphisms in~$Z(\mathcal{B})$ are morphisms of the underlying objects compatible with the half-braidings. 
We use the standard symmetric braiding on $\Vect$ and the canonical braiding~$c$ on $Z(\mathcal{B})$, with components $c_{(X,\theta),(X^\prime,\theta^\prime)}=\theta^\prime_X$. 
Then the image of $\C(G)$ is 
\begin{equation}
	(A,\theta) := J(\C(G)) \in Z(\mathcal{B})
\end{equation}
which inherits an algebra structure via transport along~$J$ as described in~\eqref{eq:Frobenius transport}. 
In this setting we have (see \cite[Thm.\,8.23.3]{EGNO2015} and \cite[Sect.\,4.2.4]{Drinfeld2009} for details): 

\begin{definition}
	Let~$\mathcal B$ be a monoidal category, and let $J\colon\Rep(G)\longrightarrow Z(\mathcal{B})$ be a braided monoidal functor. 
	The \textsl{de-equivariantisation} of~$\mathcal B$ (with respect to~$J$) is the category~$\mathcal{B}_G$ of $A$-modules in~$\mathcal B$, 
	\begin{equation}
		\mathcal{B}_G \;:=\; A\text{-}\Mod(\mathcal B) \, . 
	\end{equation}
\end{definition}
We note that there is a monoidal structure on~$\mathcal{B}_G$ which requires the half-braidings~$\theta$, and that the left action of~$G$ on $\C(G)$ induces a right action on~$\mathcal{B}_G$.

\medskip 

Equivariantisation and de-equivariantisation induce equivalences between the 2-categories of monoidal categories~$\mathcal C$ with $G$-action and of \textsl{monoidal categories over $\Rep(G)$}, i.e.\ monoidal categories~$\mathcal B$ equipped with a braided functor $\Rep(G)\longrightarrow {Z}(\mathcal{B})$ as above, see \cite[Prop.\,4.19]{Drinfeld2009}.
In particular, we have
\begin{align}
	\mathcal{C} \cong (\mathcal{C}^G)_G \label{eq:eq-deeq equivalence}
		\, , \quad 
	\mathcal{B} \cong (\mathcal{B}_G)^G
\end{align}
as equivalences in the respective 2-categories. 
On the left-hand side, we used a functor $J\colon \Rep(G)\longrightarrow Z(\mathcal{C}^G)$  extracted from~$I$ in~\eqref{def:inclusion of RepG}. 
The image of~$I$ is essentially direct sums of the monoidal unit, and thus it carries canonical half-braidings via unitors, and the functor is braided.  
Explicitly, the first equivalence in~\eqref{eq:eq-deeq equivalence} is given by
\begin{align}
	\mathcal{C}& \stackrel{\cong}{\longrightarrow} (\mathcal{C}^G)_G 
	\nonumber 
	\\
	X & \longmapsto \mathrm{Ind}(X):=\bigoplus_{g\in G} \rho_g(X)
	\label{eq:eq-deeq-eq}
\end{align}
with an $I(\C(G))$-module structure as follows. 
Choosing the basis $\{\delta_h\}_{h\in G} \subset \C(G)$ of functions supported only on a single element~$h$ (with value~1), we get that 
\begin{equation}
	\C(G)\otimes_\C \mathrm{Ind}(X)
		\cong 
		\bigoplus_{g,h\in G}\C\langle \delta_h\rangle \otimes_\C \rho_g(X)
		\cong 
		\bigoplus_{g,h\in G}\rho_g(X) \, . 
\end{equation}
The module structure is given by the projection onto the diagonal, i.e.\ those summands for which we have $g=h$. 
In other words, we think of a function $f\in\C(G)$ as acting on $\mathrm{Ind}(X)$ by $\bigoplus_{g\in G}f(g)\cdot \id_{\rho_g(X)}$, and so~$\delta_h$ projects out $\rho_h(X)$. 

\medskip 

Our last preparatory step is to enhance the above constructions with braidings, which is required in order to use \Cref{prop:FrobeniusAlgebra from 1form symmetry} in Section~\ref{sec:equivariantisation is orbifolding} below. 
More precisely, what we need is a ``braiding'' that is twisted by the $G$-action (\cite[Sect.\,4.4.3]{Drinfeld2009}): 

\begin{definition}[label=def:G-crossed braided category]
	A \textsl{$G$-crossed braided fusion category} is a fusion category $\mathcal{C}^\times_G$ together with
	\begin{enumerate}
		\item 
		a monoidal functor $\rho\colon \underline{G}^\textrm{rev}\longrightarrow \Aut^\otimes(\mathcal{C}^\times_G)$,
		\item 
		a decomposition $\mathcal{C}^\times_G=\bigoplus_{g\in G}\mathcal{C}_g$, 
		\item 
		a \textsl{$G$-braiding}~$c$ with components 
		\begin{equation}
			c_{X,Y}
			\equiv
			\hphantom{\text{{\scriptsize$\rho_{h}($}}}
			\begin{tikzpicture}[very thick,scale=0.7,color=blue!50!black,baseline]
				\draw (-1,-1) node[below] (X) {{\scriptsize$X$}};
				\draw (1,-1) node[below] (Y) {{\scriptsize$Y$}};
				\draw (1,1) node[above] (Xu) {{\scriptsize$\rho_h(X)$}};
				\draw (-1,1) node[above] (Yu) {{\scriptsize$Y\vphantom{\rho_h(X)}$}};
				\draw (1,-1)  .. controls +(0,1) and +(0,-1) .. (-1,1); 
				\draw[color=white, line width=4pt] (-1,-1)  .. controls +(0,1) and +(0,-1) .. (1,1); 
				\draw (-1,-1)  .. controls +(0,1) and +(0,-1) .. (1,1);
			\end{tikzpicture}
			\colon X\otimes Y\stackrel{\cong}{\longrightarrow} Y\otimes \rho_h(X) 
		\end{equation}
		for all $h\in G$, $X\in\mathcal{C}^\times_G$, and $Y\in \mathcal{C}_h$.
	\end{enumerate}
	These are subject to the compatibility conditions:
	\begin{enumerate}[label=(\alph*)]
		\item
		$\rho_g(\mathcal{C}_h)\subset \mathcal{C}_{g^{-1}hg}$ for all $g,h\in G$,
		\item 
		the isomorphisms $c_{X,Y}$ are natural in~$X$ and~$Y$,
		\item 
		the isomorphisms $c_{X,Y}$ are compatible with the $G$-action in the sense that for all $g\in G$ we have 
		\begin{equation}
			\rho_g(c_{X,Y})=c_{\rho_g(X),\rho_g(Y)} \, ,
		\end{equation}
		\item 
		the following identities hold for all $g,h\in G$, $Y\in \mathcal{C}_h$, and $Z\in\mathcal{C}_k$: 
		\begin{align}
		\label{eq:braiding and monoidal product inGcbc}

		\end{align}
		where we suppress the monoidal structure of~$\rho$ as we do in following.
	\end{enumerate}
\end{definition}

As an additional tool in drawing string diagrams, we define the opposite braiding~$\widetilde{c}$ whose components are
\begin{align}
	\widetilde{c}_{\,Y,X}&:=c^{-1}_{\rho_{h^{-1}}(X),Y}\circ (\id_Y\otimes ((\rho^2_{h,h^{-1}})^{-1}_X)\circ \rho^0_X))
	\equiv
	%
.
\end{align}
The braiding on the equivariantisation is given by composing the braiding in $\mathcal{C}^\times_G$ with the equivariant structure of the braided object \cite[Sect.\,4.4.4]{Drinfeld2009}:

\begin{lemma}[label=lem:eq of Gcbc is braided]
	Let $(\mathcal{C}^\times_G, \rho, c)$ be a $G$-crossed braided fusion category. 
	Its equivariantisation $(\mathcal{C}^\times_G)^G$ is an ordinary braided fusion category with braiding components 
	\begin{align}
		c^G_{(X,u),(Y,v)} 
		\;:= 
		\bigoplus_{h\in G}(\id_{Y_h}\otimes u_h) \circ c_{X,Y} 
		\;= \bigoplus_{h\in G}
		%
	\end{align}
	for $h\in G$, $X\in\mathcal{C}^\times_G$, and $Y=\bigoplus_{h\in G} Y_h$, with $Y_h\in \mathcal{C}_h$. 
\end{lemma}

In this setting, the equivariantisation/de-equivariantisation correspondence~\eqref{eq:eq-deeq-eq} is between $G$-crossed braided fusion categories and braided fusion categories~$\mathcal B$ containing $\Rep(G)$, in the sense that there is a fully faithful braided monoidal functor $\Rep(G)\longrightarrow \mathcal{B}$. 
In particular, the functor~$I$ in~\eqref{def:inclusion of RepG} is fully faithful and {braided} monoidal. 

\medskip 

Finally, recall that for an algebra~$B$ in a braided monoidal category, by using either the braiding or its inverse, we can equip any left $B$-module with two a priori different right $B$-module structures. 
A \textsl{local $B$-module} is a left $B$-module~$X$ for which these two coincide, i.e.
\begin{equation}
	\begin{tikzpicture}[very thick,scale=0.75,color=blue!50!black, baseline]
	\draw (0,-1) node[below] (X) {{\scriptsize$X$}};
	\draw[color=green!50!black] (0.6,-1) node[below] (A1) {{\scriptsize$B$}};
	\draw (0,1) node[above] (Xu) {{\scriptsize$X$}};
	\draw (0,-1) -- (0,1); 
	\draw[color=white, line width=4pt] (A1) .. controls +(0,0.6) and +(0,-0.3) .. (-0.25,0);
	\draw[color=green!50!black] (A1) .. controls +(0,0.8) and +(-0.8,-0.5) .. (0,0.3);
	\fill[color=blue!50!black] (0,0.3) circle (2.9pt) node (meet2) {};
	\end{tikzpicture} 
	=
	\begin{tikzpicture}[very thick,scale=0.75,color=blue!50!black, baseline]
	\draw (0,-1) node[below] (X) {{\scriptsize$X$}};
	\draw[color=green!50!black] (0.6,-1) node[below] (A1) {{\scriptsize$B$}};
	\draw (0,1) node[above] (Xu) {{\scriptsize$X$}};
	\draw[color=green!50!black] (A1) .. controls +(0,0.8) and +(-0.8,-0.5) .. (0,0.3);
	\draw[color=white, line width=4pt] (0,-1) -- (0,0); 
	\draw (0,-1) -- (0,1); 
	\fill[color=blue!50!black] (0,0.3) circle (2.9pt) node (meet2) {};
	\end{tikzpicture} 
\, . 
\end{equation}
Restricting de-equivariantisation to local modules over $B:=I(\C(G))$, we obtain the component $\mathcal{C}_e\subset \mathcal{C}^\times_G$ graded by the unit element $e\in G$, 
\begin{equation}
	\label{eq:deequivariantisation as local modules}
	\mathcal{C}_e
		\;\cong \;
		B\text{-}\Mod^\text{loc}\big((\mathcal{C}^\times_G)^G\big)
		\;\subset \;
		\big((\mathcal{C}^\times_G)^G\big)_G
\end{equation}
as explained in \cite[App.\,5.3.8]{Turaev2010}.

\subsection{0-form Symmetries}
\label{sec:equivariantisation is orbifolding}
In \Cref{sec:General theory of (De)Equivariantisation} we reviewed the equivariantisation $(\mathcal{C}^\times_G)^G$ of a $G$-crossed braided fusion category~$\mathcal{C}^\times_G$. 
In the present section we lay the groundwork to understand this as the gauging of a $0$-form symmetry 
\begin{equation}
	R\colon \mathrm{B}\underline{G}
		\longrightarrow 
		\mathcal{D}_{\mathcal{C}_e}
\end{equation}
that we construct below, whose codomain is the 3-category~$\mathcal{D}_{\mathcal{C}_e}$ of \Cref{exa:3dCategoriesWithAdjoints} associated to the neutral component $\mathcal C_e \subset \mathcal{C}^\times_G$. 
From~$R$ we then extract an orbifold datum~$\mathbb{A}_R$, which is the main result of this section.
In Section~\ref{sec:equivalence of equivariantisation and orbifolding} we prove that the associated ribbon category $\mathcal{C}_{\mathbb{A}_R}$ (recall Remark~\ref{rem:MFCCA}) is equivalent to $(\mathcal{C}^\times_G)^G$. 
Throughout we continue to use the 3-categorical conventions reviewed in Sections~\ref{sec:Background} and~\ref{sec:OrbDatFrom2Groups}. 

\medskip 

We construct the 3-functor $R\colon \mathrm{B}\underline{G}\longrightarrow \mathcal{D}_{\mathcal{C}_e}$ as a composition, going through the intermediate 3-category $\mathcal{G}_{\mathcal{C}^\times_G}$. 
For a given $G$-crossed braided category $(\mathcal{C}^\times_G, \rho, c)$ as in \Cref{def:G-crossed braided category}, we proceed with a description of $\mathcal{G}_{\mathcal{C}^\times_G}$, following \cite{Cui2016, Jones2020}. 
It is given by the delooping of a monoidal $2$-category $G\rtimes \mathcal{C}^\times_G$. 
The objects of this $2$-category are the elements of $G$, and the Hom categories are 
\begin{equation}
	\Hom(g,h) := \lbrace g\rbrace \times \mathcal{C}_{hg^{-1}} \, . 
\end{equation}
The vertical composition of $2$-morphisms (denoted~$\circ$) is the composition of morphisms in $\mathcal{C}_{hg^{-1}}$, while horizontal composition (denoted~$\otimes$) is the monoidal product in~$\mathcal{C}^\times_G$. 
More precisely, given $1$-morphisms $(g,X)\in\Hom(g,h)$ and $(h,Y)\in \Hom(h,k)$, we have 
\begin{equation}
	(h,Y)\otimes (g,X):=(g, Y\otimes X) \, .
\end{equation}
Naturally, the identity $1$-morphisms $\id_g$ are given by $(g,\mathds{1})$ using the monoidal unit  $\mathds{1} \in \mathcal C_e \subset \mathcal{C}^\times_G$ and the unitors and associator also coincide with those in $\Gcbc$. 
Lastly, the monoidal structure~$\boxtimes$ of $G\rtimes \mathcal{C}^\times_G$ is given by the group multiplication on objects and by a semidirect product on $1$- and $2$-morphisms:
\begin{align}
\label{eq:braiding on equivariantisation}
	g\boxtimes g^\prime &:=gg^\prime\,,
	\\
	(g,X)\boxtimes (g^\prime,X^\prime)&:=\big(gg^\prime, X\otimes \rho_{g^{-1}}(X^\prime)\big)\,,
	\\
	(g,f)\boxtimes (g^\prime,f^\prime)&:=\big(gg^\prime, f\otimes \rho_{g^{-1}}(f^\prime)\big)
\end{align}
for objects $g, g^\prime$, $1$-morphisms $(g,X), (g^\prime,X^\prime)$ and $2$-morphisms $(g,f), (g^\prime,f^\prime)$.

To complete~$\boxtimes$ to a $2$-functor, it has to be equipped with pseudonatural transformations $\boxtimes^2$ (interchanger) and $\boxtimes^0$ (unitor). 
This can be done by using the crossed braiding of $\mathcal{C}^\times_G$ together with~$\rho^2$ and~$\rho^0$, respectively. 
Lastly, the unitors, associator, triangulators and pentagonator for $\boxtimes$ are defined canonically up to a normalisation, as explained in \cite[Rem.\,6.1]{Cui2016}.

Now we can simply deloop this monoidal 2-category, resulting in the desired 3-category 
\begin{equation}
	\label{eq:3cat of Gcbc}
	\mathcal{G}_{\mathcal{C}^\times_G} 
		\;:=\; 
		\mathrm{B}(G\rtimes\mathcal{C}^\times_G) \, . 
\end{equation}
The graphical calculus for this 3-category corresponds to the one we used for 2-groups (see \Cref{sec:2-groups}). The only difference is that the line defects are labelled by objects of $\Gcbc$ instead of a group, and the target surface is encoded in the degree of the object, i.e.\ we have $(g,X)\in\Hom_{\mathcal{G}_{\mathcal{C}^\times_G} }(g,t(X)g)$ where $t(X):=k$ for $X\in \mathcal{C}_k$.

\medskip 

Next we turn to symmetries in the 3-category~$\mathcal{D}_{\mathcal{C}_e}$ from Example~\ref{exa:3dCategoriesWithAdjoints}. 
Recall that it is the delooping of the 2-category of $\Delta$-separable symmetric Frobenius algebras, bimodules and bimodule morphisms in the ribbon category~$\mathcal C_e$. 
To construct our $0$-form symmetry as a $3$-functor 
\begin{equation}
	R\equiv(R,\chi,\iota,\omega,\gamma, \delta)\colon 
	\mathrm{B}\underline{G}\longrightarrow \mathcal{D}_{\mathcal{C}_e} \, ,
\end{equation}
we continue to use the notation of \Cref{subsubsec:3cats}.
Since $\mathrm{B}\underline{G}$ has only a single object, there is only a single pair/triple/quadruple of objects which may label the components of the coherence data, hence there is exactly one pseudonatural transformation~$\chi$, and similarly for the remaining coherence data.

We define~$R$ as the composition of two $3$-functors. 
The first piece is the inclusion 
\begin{equation}
	\mathrm{B}\underline{G}\longhookrightarrow \mathcal{G}_{\mathcal{C}^\times_G}\,.
\end{equation} 
In both categories, there is only one object and $1$-morphisms are given by elements of~$G$ with their group multiplication. 
All higher morphisms are identities in the domain $\mathrm{B}\underline{G}$, and the inclusion maps them to the associated identities in $\mathcal{G}_{\mathcal{C}^\times_G}$. 
All coherence morphisms are identities.

The second piece is a $3$-functor 
$
	M_G\colon \mathcal{G}_{\mathcal{C}^\times_G}\longrightarrow \mathcal{D}_{\mathcal{C}_e}
$ 
which we define as follows. 
We arbitrarily pick a simple object $m_g\in\mathcal{C}_g$ for each $g\in G$, and we write $m_g^* \in \mathcal C_{g^{-1}}$ for its dual. 
For objects and morphisms we then set
\begin{align}
	M_G\colon \mathcal{G}_{\mathcal{C}^\times_G} & \longrightarrow \mathcal{D}_{\mathcal{C}_e}
	\nonumber
	\\
	\ast &\longmapsto \ast
	\\
	g &\longmapsto A_g:=m_g^*\otimes m_g
	\\
	\Hom(g,h)\;\ni\; (g,X)&\longmapsto m_h^*\otimes X\otimes m_g
	\\
	\Hom\big((g,X),(g,Y)\big)\;\ni\; f&\longmapsto \id_{m_h^*}\otimes f\otimes \id_{m_g}\,.
\end{align}
Note that once again, all the images are in~$\mathcal{C}_e$. 

In the following, when we draw string diagrams (in~$\mathcal{C}^\times_G$) involving~$m_g$ or its dual, we only label strings by the group element~$g$, with an arrow pointing upwards indicating~$m_g$ whereas an arrow pointing downwards indicates~$m_g^*$. 
In particular, we have
\begin{equation}

	=\id_{m_g^*}\,.
\end{equation}
For later convenience, we also introduce 
\begin{equation}
	\psi:= \sum_{g\in G}(\mathrm{dim}(m_g))^{-\frac{1}{2}}\cdot \id_{A_g}\,.
\end{equation}
By slight abuse of notation, we write $\psi^2$ also when we mean $(\mathrm{dim}(m_g))^{-1}\cdot\id_{A_g}$ in diagrams, instead of the entire sum. 
We do not distinguish between~$\psi$ and its dual~$\psi^*$, as we have $\mathrm{dim}(m_g)=\mathrm{dim}(m_g^*)$ by virtue of~$C_G^\times$ being fusion. 
As an illustration of these conventions, we have that 
\begin{equation}
	%
	=\mathrm{dim}(m_g^*)^{-1}\cdot\id_{m_g^*}\,.
\end{equation}

We endow 
\begin{equation}
	A_g := m_g^*\otimes m_g
\end{equation}
with the structure of a $\Delta$-separable Frobenius algebra in~$\mathcal C_e$ by specifying the following (co)multiplication and (co)unit, adapted from \cite{Carqueville2020}: 

\begin{align}
\label{eq:AgGextension}
&\mu_g:=\,
%
\!\!.
\end{align}
Moreover, for each $g\in G$ we have a left $A_g$-module structure on $m_g^*\in\mathcal{C}_{g^{-1}}$, and a right module structure on $m_g\in\mathcal{C}_g$ given by
\begin{align}
\label{eq:mmodulestructures}
%
 
.
\end{align}
\noindent
Using these, we can also construct $(A_{gh},A_g\otimes A_h)$-bimodules 
\begin{equation}
	\label{eq:chi-def}
	\chi_{g,h} :=m_{gh}^* \otimes m_g \otimes m_h \; \in \mathcal C_e 
\end{equation}
with module structure
\begin{align}
& \sigma^0_{g,h}:=
%
 
\label{eq:Tghactions}
\end{align}%

\noindent{}where here and below the labels~1 and~2 in string diagrams refer to the right actions of $A_g\otimes \mathds{1}$ and $\mathds{1} \otimes A_h$, respectively, coming from the $(A_g\otimes A_h)$-action. 
The bimodules~$\chi_{g,h}$ are denoted $T_{g,h}$ in \cite{Carqueville2020}; here we stick to our notational conventions of 3-functors. 

The 2-morphisms $\chi_{g,h}\colon A_g\otimes A_h\longrightarrow A_{gh}$ in~\eqref{eq:chi-def} are in fact equivalences, and their duals $\chi_{g,h}^*=m_h^*\otimes m_g^*\otimes m_{gh}$ are their weak inverses. 
This follows from
\begin{align}
	\chi_{g,h}\otimes_{A_g\otimes A_h}\chi_{g,h}^* & \cong  A_{gh}\,,
	\\
	\chi_{g,h}^*\otimes_{A_{gh}}\chi_{g,h} & \cong m_h^*\otimes m_g^*\otimes m_g\otimes m_h \cong A_g\otimes A_h \, ,
	\label{eq:chiStarRelativeChi}
\end{align}
where we used $m_g\otimes_{A_g} m^*_g \cong \mathds{1}$ (which follows from standard arguments of computing relative tensor products over $\Delta$-separable Frobenius algebras by splitting idempotents, see e.g.\ \Cref{subsubsec:2dOrbDat} and \cite[Sect.\,3]{Carqueville2020}), and the last isomorphism in~\eqref{eq:chiStarRelativeChi} is the braiding.
The $(A_h,A_g)$-bimodule structures on $2$-morphisms of the form $m_h^*\otimes X\otimes m_g$ are defined in complete analogy to~\eqref{eq:Tghactions}.

We continue by describing the coherence data for the $2$-functor 
\begin{equation}
	\label{eq:MG2functor}
	(M_G)_{\ast,\ast}\colon \End_{\mathcal{G}_{\mathcal{C}^\times_G}}(\ast)\longrightarrow \End_{\mathcal{D}_{\mathcal{C}_e}}(\ast)\,,
\end{equation} 
where we only spelled out the subscript ``$\ast,\ast$'' here for clarity.
Using $m_g\otimes_{A_g} m^*_g \cong \mathds{1}$ as well as unitors, we define $M_G^2$ to be the pseudonatural transformation whose components are given by the canonical isomorphisms (again using $m_h\otimes_{A_h} m^*_h \cong \mathds{1}$)
\begin{equation}
	(M_G^2)_{(h,Y),(g,X)}\colon 
		\underbrace{M_G\big((h,Y)\otimes (g,X)\big)}_{=m_k^*\otimes Y\otimes X\otimes m_g}
			\stackrel{\cong}{\longrightarrow} 
		\underbrace{M_G\big((h,Y)\big)\otimes_{A_h} M_G\big((g,X)\big)}_{=m_k^*\otimes Y\otimes m_h\otimes_{A_h} m^*_h\otimes X\otimes m_g} .
\end{equation}
The identity 2-morphism of the algebra~$A_g$ is given by~$A_g$ viewed as a bimodule over itself. 
This coincides with $M_G((g,\mathds{1}))= A_g$, and we set $M_G^0:=\id_{A_g}$. 
This completes the definition of the $2$-functor~\eqref{eq:MG2functor}, and one readily checks that the associativity and unitality axioms are satisfied.

\medskip 

We proceed to specify the required pseudonatural transformations $\chi, \iota$ and modifications $\omega,\gamma, \delta$ associated to the $3$-functor~$M_G$. 
Since the inclusion $\mathrm{B}\underline{G}\longhookrightarrow \mathcal{G}_{\mathcal{C}^\times_G}$ is a strict $3$-functor, the coherence data of~$R$ will coincide with that of~$M_G$, and we use the same symbols for both. 
For composition of $1$-morphisms, we define the component $2$-equivalences of~$\chi$ as the bimodules given in \eqref{eq:chi-def}, 
\begin{equation}
	\chi_{g,h}\colon \underbrace{M_G(g)\boxtimes M_G(h)}_{=A_g\otimes A_h} 
	\stackrel{\cong}{\longrightarrow} 
	\underbrace{M_G(gh)}_{=A_{gh}}\,.
\end{equation}
Its $3$-morphism components $\chi_{(g,X), (g^\prime, X^\prime)}$ are of the form
\begin{equation}
	\chi_{h,h^\prime}\otimes_{A_h\otimes A_{h^\prime}} \big(M_G(g,X)\otimes M_G(g^\prime,X^\prime)\big)
	\stackrel{\cong}{\longrightarrow} 
	 M_G\big((g,X)\boxtimes(g^\prime,X^\prime)\big)\otimes_{A_{gg^\prime}}\chi_{g,g^\prime}\,,
\end{equation}
that is, they fill the diagrams
\begin{equation}
	\begin{tikzcd}[column sep= 95,row sep= 50, ampersand replacement=\&]
	A_g\otimes A_{g^\prime}
	\arrow[d, "{\chi_{g,g^\prime}}"{left}]
	\arrow[r, "{M_G(g,X)\otimes M_G(g^\prime,X^\prime)}"]
	\& A_h\otimes A_{h^\prime}
	\arrow[dl, Rightarrow,shorten <= 5ex, shorten >= 5ex, "{\chi_{(g,X), (g^\prime, X^\prime)}}"]
	\arrow[d, "{\chi_{h,h^\prime}}"]\\
	A_{gg^\prime}
	\arrow[r, "{M_G((g,X)\boxtimes(g^\prime,X^\prime))}"{below}]
	\& A_{hh^\prime}
\end{tikzcd} 
\end{equation}
where composition is given by relative tensor products in $\mathcal{C}_e$. 
We define these 3-morphisms to be
\begin{equation}
	\chi_{(g,X),(g^\prime,X^\prime)} 
	\;\;:=\quad 

	\, . 
\end{equation}
As before, this string diagram is in~$\mathcal{C}^\times_G$ and gives the balanced version of the morphism (recall the discussion around~\eqref{eq:BalancedCobalanced}).

Lastly, we need to provide an invertible modification~$\omega$ that acts as the associator with respect to~$\chi$. 
Its components take the form
\begin{equation}
	\omega_{g,h,k}\colon \chi_{gh,k}\otimes_{A_{gh}\otimes A_k}(\chi_{g,h}\otimes A_k)
	\stackrel{\cong}{\longrightarrow} 
	\chi_{g,hk}\otimes_{A_g\otimes A_{hk}}(A_g\otimes \chi_{h,k})\label{eq:1-associator}
\end{equation}
and we define them to be the following morphisms in~$\mathcal{C}_e$:
\begin{equation}
	\omega_{g,h,k}:=\quad
\,.
\end{equation}
We suppressed the associators because they are trivial on the left-hand side (of \eqref{eq:1-associator}) since~$G$ is associative, while on the right-hand side the associator of~$\mathcal{C}_e$ is invisible in string diagrams as usual.

Regarding unitors, note first that $M_G(e)=A_e\cong\mathds{1}$ via $\iota:=\mathds{1}$ as a $(\mathds{1},A_e)$-bimodule, and $\chi_{g,e}=A_g\otimes \mathds{1}$ as bimodules. 
Thus, all components of the unitor modifications~$\delta$ and~$\gamma$ are given by unitors in~$\mathcal{C}_e$, evaluations, and coevaluations.

\begin{lemma}
	The ingredients described above assemble into a 3-functor
	\begin{equation}
		\label{eq:0-form symmetry functor}
		R\equiv(R,\chi,\iota,\omega,\gamma, \delta)
			:= 
			\Big( \mathrm{B}\underline{G} 
			\longhookrightarrow 
			\mathcal{G}_{\mathcal{C}^\times_G}
			\stackrel{M_G}{\longrightarrow} 
			\mathcal{D}_{\mathcal{C}_e} \Big) .
	\end{equation}
\end{lemma}

\begin{proof}
	Since the inclusion $\mathrm{B}\underline{G} \longhookrightarrow \mathcal{G}_{\mathcal{C}^\times_G}$ is a strict 3-functor, what has to be checked is that the structure morphisms of~$M_G$ satisfy the coherence axioms of a 3-functor. 
	This amounts to a tedious yet direct computation of standard string diagram manipulations of the type also carried out in \Cref{sec:equivalence of equivariantisation and orbifolding}. 
	We refrain from giving all details here, but to illustrate the general idea, we check the modification axiom of~$\omega$. 
	As an identity of pasting diagrams, it reads 
	\begin{equation}
			\begin{tikzcd}[column sep= 50,row sep= 35]
			F(\xi)
			\arrow[d, bend right=50,"{\alpha_\xi}"{right}]
			\arrow[r, "{F(f)}"]
			& F(\zeta)
			\arrow[d, bend right=50, "{\alpha_\zeta}"{left}]{r}[name=a]{}
			\arrow[d, bend left=50, "{\beta_\zeta}"]{}[name=b]{}
			\arrow[Rightarrow, "{\omega_\zeta}"{above},shorten >=1ex, shorten <=1ex, from=a, to=b]\\
			F^\prime(\xi)
			\arrow[ur, Rightarrow,shorten <= 4ex, shorten >= 4ex, "{\alpha_f^{-1}}", xshift=-3ex]
			\arrow[r, "{F^\prime(f)}"{below}]
			& F^\prime(\zeta)
		\end{tikzcd}
		\;\;=\;\;
			\begin{tikzcd}[column sep= 50,row sep= 35]
			F(\xi)
			\arrow[d, bend right=50,"{\alpha_\xi}"{left}]{r}[name=a]{}
			\arrow[d, bend left=50, "{\beta_\xi}"]{}[name=b]{}
			\arrow[r, "{F(f)}"]
			& F(\zeta)
			\arrow[d, bend left=50, "{\beta_\zeta}"]
			\arrow[Rightarrow, "{\omega_\xi}"{above},shorten >=1ex, shorten <=1ex, from=a, to=b]\\
			F^\prime(\xi)
			\arrow[ur, Rightarrow,shorten <= 4ex, shorten >= 4ex, "{\beta_f^{-1}}", xshift=3ex]
			\arrow[r, "{F^\prime(f)}"{below}]
			& F^\prime(\zeta)
		\end{tikzcd}.
\end{equation}
where $F:=\boxtimes\,\circ (\boxtimes\times \id)\circ (M_G^{\times 3})$, $F^\prime:=M_G\circ \,\boxtimes\,\circ (\boxtimes\times \id)$, $\xi:=(g,h,k)$, $\zeta:=(g^\prime,h^\prime,k^\prime)$, $f:=((g,X),(h,Y),(k,Z))$, $\alpha:=\chi\otimes (\chi\,\boxtimes \,\id_{M_G})$, and lastly $\beta:=\chi\otimes (\id_{M_G}\,\boxtimes \,\chi)$. Since this includes the identity pseudonatural transformation of~$M_G$ and whiskering with~$\boxtimes$, it involves in particular the unitors \eqref{eq:mmodulestructures} (as the 3-morphism components of the identity pseudonatural transformation) and various braidings (since whiskering with the $\boxtimes$-composition introduces $\boxtimes^2$). The corresponding identity of string diagrams is thus given by
\begin{equation}
\label{eq:string diagram for modification axiom}

\end{equation}	
where the horizontal blue dotted line at the bottom left indicates the tensor product $X\otimes \rho_g^{-1}(Y)\otimes\rho^{-1}_{gh}(Z)$ which splits into~$X$ and $\rho_g^{-1}(Y)\otimes\rho^{-1}_{gh}(Z)$ in the left diagram and into $X\otimes \rho_g^{-1}(Y)$ and $\rho^{-1}_{gh}(Z)$ in the right diagram. 
Overall, the diagram is a morphism
\begin{align}
	&\Big(m_{g^\prime h^\prime k^\prime}^*\otimes (X\otimes \rho_g^{-1}(Y)\otimes\rho^{-1}_{gh}(Z)) \otimes m_{ghk}\Big) \otimes T_{gh,k}\otimes T_{g,h}\otimes A_k \nonumber\\
	&\longrightarrow T_{g^\prime, h^\prime k^\prime}\otimes A_{g^\prime}\otimes T_{h^\prime, k^\prime}\otimes \left(m_{g^\prime}^*\otimes X\otimes m_g\right)\otimes \left(m_{h^\prime}^*\otimes Y\otimes m_h\right)\otimes \left(m_{k^\prime}^*\otimes Z\otimes m_k\right).
\end{align}
By carefully removing all the loops (their contributions are cancelled by their $\psi^2$-insertions) and tracing the remaining lines, one observes that the equality holds.
\end{proof}

We note that our construction of the 0-form symmetry~$R$ depends on a choice of simple objects~$m_g$ for each $g\in G$. 
However, any two such choices give rise to Morita equivalent algebras \cite[Cor.\,5.4]{Carqueville2020}, and we expect these choices to produce equivalent $3$-functors.

\medskip 

From~$R$ we now extract the main algebraic entity of this section. 
Recall from Example~\ref{exa:3dOrbifoldData} the notion of an orbifold datum in the 3-category~$\mathcal D_{\mathcal C}$ we associate to a ribbon category~$\mathcal C$. 
The orbifold datum we specify below is essentially the one we constructed for a general 0-form symmetry in \eqref{eq:0form orbdat in 3d}. 
However, it turns out that this is not quite an orbifold datum, therefore we modify it by introducing~$\psi$ and~$\phi$ explicitly and adjusting~$\alpha$ and~$\overline{\alpha}$ by factors of $\psi^{-1}$. 
(Such modifications are familiar e.g.\ from obtaining Turaev--Viro--Barrett--Westbury models as orbifold TQFTs, see \cite[Sect.\,4]{Carqueville2020}.)

\begin{lemma}[label=lem:0form orb dat from Gcbrc]
	The data of the 3-functor~$R$ in~\eqref{eq:0-form symmetry functor} gives rise to an orbifold datum 
	\begin{equation}
		\mathbb{A}_R := \big(A,T,\alpha,\overline{\alpha}, \psi,\phi \big)
	\end{equation}
	in~$\mathcal D_{\mathcal C_e}$ as follows: 
	\begin{align}
		A&:=\bigoplus_{g\in G}R(g)=\bigoplus_{g\in G}A_g\,,\\
		T&:=\bigoplus_{g,h\in G}\chi_{g,h}\,,\\
		\label{eq:alphaG}
		\alpha&:=\bigoplus_{g,h,k\in G} \mathrm{dim}(m_{gh}) (\sigma_{gh,k}^2\otimes \id_{\chi_{g,h}})\circ (\id_{\chi_{gh,k}}\otimes c_{\chi_{g,h},A_k})\circ\omega_{g,h,k}^{-1}\circ  (\id_{\chi}\otimes \eta_k\otimes \id_{\chi})\,,\\
		\label{eq:alphabG}
		\overline{\alpha}&:= \bigoplus_{g,h,k\in G} \mathrm{dim}(m_{hk})(\sigma_{g,hk}^1\otimes \id_{\chi_{h,k}})\circ \omega_{g,h,k}\circ (\id_{\chi\otimes \chi}\otimes \eta_k)\,,\\
		\psi&:=\bigoplus_{g\in G} \left(\mathrm{dim}(m_g)\right)^{-\frac{1}{2}}\cdot \id_{A_g}\,, \label{eq:psiInAR}
		\\
		\phi&:=\frac{1}{\sqrt{|G|}}\,.
		\label{eq:phiInAR}
	\end{align}
\end{lemma}

\begin{proof}
	As noted in~\Cref{exa:3dOrbifoldData}, orbifold data in~$\mathcal{D}_{\mathcal{C}_e}$ involve partially evaluating the relative tensor products in the domain and codomain of~$\alpha$. 
	We make this explicit by composing $\omega_{g,h,k}$ and its inverse with the module action~\eqref{eq:mmodulestructures} which serves as the retraction onto the relative tensor product, and with the unit $\eta_g$ from \eqref{eq:AgGextension} which is its section. 
	In the case of~$\alpha$, the full module action includes the braiding $c_{\,\chi_{g,h},A_k}$.
	
	These compositions are written out in \eqref{eq:alphaG} and \eqref{eq:alphabG}. They yield the $3$-morphisms~$\alpha$ and~$\overline{\alpha}$, respectively, which in terms of string diagrams are given by 
	\begin{align}
	\label{eq:string alphaG}
		\alpha&=\bigoplus_{g,h,k\in G}
.
	\end{align}
	As a consequence, for the functor~$R$ in~\eqref{eq:0-form symmetry functor}, we recover our $\mathbb{A}_R=(A, T, \alpha,\overline{\alpha}, \psi,\phi)$ as the special orbifold datum constructed in \cite[Thm.~5.1]{Carqueville2020}. 
\end{proof}

\subsection{1-form Symmetries}
\label{sec:Deeq is 1-form orbifolding}

In \Cref{sec:equivariantisation is orbifolding} we produced an orbifold datum~$\mathbb{A}_R$ from a 0-form symmetry $R\colon \mathrm{B}\underline{G}\longrightarrow \mathcal{D}_{\mathcal{C}_e}$, such that $\mathbb{A}_R$-orbifolding is the same as equivariantisation (proven in \Cref{sec:equivalence of equivariantisation and orbifolding}). 
In the present section we construct an orbifold datum~$\mathbb{A}_{\widehat{R}}$ from a 1-form symmetry~$\widehat{R}$, 
and we show in which sense $\mathbb{A}_{\widehat{R}}$-orbifolding undoes equivariantisation. 

\medskip 

Let~$G$ be any finite group, and let~$\mathcal C_G^\times$ be a $G$-crossed braided fusion category. 
Recall the functor 
\begin{equation}
	\label{eq:IAgain}
	I\colon\Rep(G)\longrightarrow (\mathcal{C}^\times_G)^G	
\end{equation}
from~\eqref{def:inclusion of RepG} (with $\mathcal C = \mathcal{C}^\times_G$). 
We first show how to obtain an algebra~$B$ and an orbifold datum~$\mathbb{A}_B$ in $(\mathcal{C}^\times_G)^G$ from~$I$. 
Then we show that~$\mathbb{A}_B$ in fact arises from a 1-form symmetry~$\widehat{R}$. 

\medskip 

We start by noting that the function algebra $\C(G)$ has the structure of a Frobenius algebra in $\Rep(G)$. 
Indeed, we already saw that together with the $G$-action~\eqref{eq:right action on alg of functions}, $\C(G)$ is an algebra in $\Rep(G)$. 
Further equipping it with the counit $\varepsilon\colon \delta_g \longmapsto 1$ and comultiplication $\Delta\colon \delta_g \longmapsto \delta_g\otimes \delta_g$ for all $g\in G$ turns $\C(G)$ into a commutative $\Delta$-separable Frobenius algebra in $\Rep(G)$, as a simple calculation shows.
Here as before $\{\delta_g\}_{g\in G} \subset \C(G)$ is the basis of functions with value 1 supported on only one element. 

The functor~$I$ allows us to associate to the Frobenius algebra $\C(G)$ an orbifold datum in $(\mathcal{C}^\times_G)^G$ such that gauging it amounts to extracting the neutral component~$\mathcal C_e$ of the $G$-crossed braided fusion category: 

\begin{proposition}[label=prop:Deequivariantisation as orbifold datum1]
	Let~$\mathcal C_G^\times$ be a $G$-crossed braided fusion category. 
	Then 
	\begin{equation}
		B := I\big(\C(G)\big)
	\end{equation}
	is a commutative $\Delta$-separable symmetric Frobenius algebra in $(\mathcal{C}^\times_G)^G$. 
	Moreover, the associated orbifold datum~$\mathbb{A}_B$ (via \Cref{prop:orbifold datum from Frobenius algebra}) is such that
	\begin{equation}
	\label{eq:non-commutative deequivariantisation via orbifold}
		\left((\mathcal{C}^\times_G)^G\right)_{\mathbb{A}_B}
		\; \cong \; 
		\mathcal{C}_e \, .
	\end{equation}
\end{proposition}

\begin{proof}
	Recall that the functor~$I$ is fully faithful braided monoidal and factors through the symmetric centre. 
	Therefore, the image~$B$ of $\C(G)$ carries the structure of a Frobenius algebra and it is commutative and $\Delta$-separable. 
	Since~$I$ also preserves duals, the resulting algebra is symmetric. 
	This proves the first claim.
	
	Applying \Cref{prop:orbifold datum from Frobenius algebra} to~$B$ we obtain an orbifold datum~$\mathbb{A}_B$ in $(\mathcal{C}^\times_G)^G$. 
	In fact~$B$ is also haploid (meaning $\mathrm{dim}_\C(\Hom(\mathds{1},B))=1$), because~$I$ is fully faithful. 
	Thus we can apply \cite[Thm.~4.1]{Mulevicius2022} to find that orbifolding with~$\mathbb{A}_B$ gives the ribbon category of local $B$-modules, 
	\begin{equation}
		\label{eq:orbifold category of a condensable algebra}
		\big((\mathcal{C}^\times_G)^G\big)_{\mathbb{A}_B}
		\;\cong\; 
		B\text{-}\Mod^\textrm{loc}\big((\mathcal{C}^\times_G)^G\big)\,.
	\end{equation}
	Combining this with~\eqref{eq:deequivariantisation as local modules} implies the second claim.
\end{proof}

\medskip 

Next we give a ``dual'' description of the above results. 
For any group, we can consider the \textsl{Pontryagin dual group} 
\begin{equation}
	\widehat{G} := \Hom_{\Grp}\big(G,\textrm{U}(1)\big)
\end{equation}
with pointwise multiplication. 
Its \textsl{group algebra} $\C[\widehat{G}]$ is the vector space spanned by the elements of~$\widehat{G}$ together with multiplication induced by multiplication in~$\textrm{U}(1)$. 
Its unit is the constant function with value~$1$, comultiplication is given by the sum over tensor products of functions whose product is the original function (all divided by~$|G|$), and the counit maps the unit to~$|G|$ and every other function to~$0$.
This turns $\C[\widehat{G}]$ into a Frobenius algebra over~$\C$, which lifts to a Frobenius algebra in $\Rep(G)$ (see the proof of \Cref{lem:equivalence of group algebras} below for details).

As the function algebra $\C(G)$ of the original group~$G$ is also a Frobenius algebra in $\Rep(G)$, a natural question is how it compares to $\C[\widehat{G}]$. 
Here we note that the functor~$I$ in~\eqref{eq:IAgain} pertains to a non-invertible symmetry (its domain $\Rep(G)$ is not generated by invertible objects) unless~$G$ is commutative. 
In the context of equivariantisation and de-equivariantisation, this is one motivation to assume that~$G$ is commutative.
Then both Frobenius algebras are isomorphic: 

\begin{lemma}[label=lem:equivalence of group algebras]
	If $G$ is commutative, then $\C[\widehat{G}]\cong \C(G)$ as Frobenius algebras in $\Rep(G)$.
\end{lemma}

\begin{proof}
	Finite commutative groups can be decomposed into direct sums of cyclic groups, i.e.\ given such a group~$G$, there is a list of natural numbers $(n_1,\ldots,n_k)$ such that
	$
	G \;\cong\; \bigoplus_{i=1}^k\Z_{n_i}
	$. 
	This decomposition carries over to the Pontryagin group, 
	\begin{equation}
		\widehat{G} \cong \bigoplus_{i=1}^k\Hom_{\Grp}\big(\Z_{n_i},\textrm{U}(1)\big)\,.
	\end{equation}
	Thus, if we show the claim for cyclic groups, then the general statement follows.
	
	The dual group~$\widehat{\Z}_n$ consists of~$n$ elements that can be identified with their image of $1\in\Z_n$, 
	\begin{align}
		\phi_j\colon \Z_n & \longrightarrow \textrm{U}(1)
		\\
		1 & \longmapsto \textrm{e}^{2\pi\textrm{i}j/n} 
	\end{align}
	for $0\leqslant j\leqslant n-1$. 
	The group algebra $\C[\widehat{\Z}_n]$ is the $\C$-linear hull of these functions equipped with multiplication. 
	The comparison with $\C(\Z_n)$ is immediate, as it also consists of such functions. 
	Setting $\xi:=\phi_1(1)=\textrm{e}^{2\pi\textrm{i}/n}$, the change of basis between $\{\delta_l\}$ and $\{\phi_j\}$ is 
	\begin{align}
		\phi_j =\sum_{l=0}^{n-1}\xi^{lj}\delta_l 
		\, ,\quad 
		\delta_l=\frac{1}{n}\sum_{j=0}^{n-1}\xi^{-lj}\phi_j \, ,
	\end{align}
	which is a finite Fourier transformation.
	This shows that $\C[\widehat{\Z}_n]=\C(\Z_n)$ as vector spaces. 
	
	To show equality also as Frobenius algebras and $G$-representations, we summarise both sets of structures in Table~\ref{tab:Structures on group algebras}, where $\delta_{l,l^\prime}$ is the Kronecker delta and all indices are mod$\,n$. 
	It is elementary to show that these prescriptions are compatible with the change of basis given above, i.e.\ they define the same structures on the underlying vector space. 
	\begin{table}[t]
		\begin{center}
			\begin{tabular}{l|c|c}
				Structure&$\C(\Z_n)$&$\C[\widehat{\Z}_n]$\\\hline
				\alignCenterstack{&G\text{-action}\\&\text{multiplication}\\&\text{unit}\\&\text{comultiplication}\\&\text{counit}}&
				\alignCenterstack{
				\delta_l\tril g&\;:=\delta_{g+l}\\
				\delta_l\cdot\delta_{l^\prime}&\;:=\delta_{l,l^\prime}\delta_l\\
				\eta(1)&\;:=\sum_{l=0}^{n-1}\delta_l\\
				\Delta(\delta_l) &\;:= \delta_l\otimes\delta_l\\
				\varepsilon (\delta_l)&\;:=1}&
				\alignCenterstack{
					\phi_j\tril g&\;:=\xi^{-gj}\phi_j\\
					\phi_j\cdot\phi_{j^\prime}&\;:=\phi_{j+j^\prime}\\
					\eta(1)&\;:= \phi_0\\
					\Delta(\phi_j) &\;:= \frac{1}{n}\sum_{k=0}^{n-1}\phi_k\otimes \phi_{j-k}\\
					\varepsilon (\phi_j)&\;:=n\delta_{j,0}}
			\end{tabular}
			\caption{Frobenius algebra and representation structures on $\C[\widehat{\Z}_n]$ and $\C(\Z_n)$.}
			\label{tab:Structures on group algebras}
		\end{center}
	\end{table}
	
	We conclude that $\C[\widehat{\Z}_n]=\C(\Z_n)$ as Frobenius algebras in $\Rep(\Z_n)$, and thus $\C[\widehat{G}]\cong \C(G)$ for any finite commutative group~$G$. 
\end{proof}

Next we construct the $1$-form symmetry
\begin{equation}
	\label{eq:1FormSymmetryHat}
	\widehat{R}\colon \mathrm{B}^2\widehat{G}\text{-}\Vect 
	\longrightarrow 
	\mathcal{D}_{(\mathcal{C}^\times_G)^G} 
\end{equation}
and show that its associated orbifold datum~$\mathbb{A}_{\widehat R}$ coincides with the one in \Cref{prop:Deequivariantisation as orbifold datum1}. 
For this we use the functor $I\colon\Rep(G)\longrightarrow (\mathcal{C}^\times_G)^G$ again and compose it with a functor $\widehat{G}\text{-}\Vect\longrightarrow \Rep(G)$. 
Delooping their composition twice then gives the 1-form symmetry~$\widehat R$.
Note that its codomain is well-defined because $(\mathcal{C}^\times_G)^G$ is a ribbon category. 

Given a finite group $G$, we have a braided monoidal functor 
\begin{align}
	\Phi \colon \widehat{G}\text{-}\Vect&\longrightarrow \Rep(G) 
	\nonumber
	\\
	\C_\phi & \longmapsto (\C,\phi)
\end{align}
where $\C_\phi$ is the 1-dimensional vector space concentrated in degree $\phi\in\widehat{G}$. 
The monoidal structure of~$\Phi$ is given by
\begin{align}
	\Phi^2_{\phi,\psi}\colon (\C,\phi)\otimes_\C (\C,\psi) \stackrel{\cong}{\longrightarrow} (\C,\phi\psi) 
	\, , \quad 
	\Phi^0:=\id\colon (\C,1) \stackrel{\cong}{\longrightarrow} \Phi(\C_1)= (\C,1) 
\end{align}
where the former is given by the canonical identification $\C\otimes_\C\C\cong\C$. 
It is straightforward to check that these data satisfy the axioms of a braided monoidal functor as almost all of the involved morphisms are identities. 
From this we define the 3-functor~$\widehat{R}$ as the double-delooping of the composition of~$I$ and~$\Phi$, i.e.\ such that the following diagram commutes: 
\begin{align}
	\label{eq:dual 1-form symmetry}
	\begin{tikzcd}[column sep= 40,row sep= 20, ampersand replacement=\&]
		\mathrm{B}^2 \widehat{G}\text{-}\Vect 
		\arrow[dr, swap, "{\mathrm{B}^2\Phi}"]
		\arrow[rrr, "{\widehat{R}}"]
		\& \& \& \mathcal{D}_{(\mathcal{C}^\times_G)^G} 
		\\
		\&
		\mathrm{B}^2\Rep(G)
		\arrow[r, "{\mathrm{B}^2I}"{below}]
		\& \mathrm{B}^2(\mathcal{C}^\times_G)^G 
		\arrow[ur, hookrightarrow] 
		\&
	\end{tikzcd}
\end{align}

\begin{remark}
	For a finite commutative group~$G$, the functor~$\Phi$ is a braided equivalence $\widehat{G}\text{-}\Vect\cong\Rep(G)$. 
	On the right-hand side, irreducible representations of finite commutative groups are $1$-dimensional, and they are indexed by characters, whereas on the left-hand side the simple objects are $1$-dimensional vector spaces associated to elements of $\widehat{G}$, which are exactly the characters. 
	The regular representation $\C(G)$ is the sum over all simple objects in $\Rep(G)$, and so is $\C[\widehat{G}]$ in $\widehat{G}\text{-}\Vect$.
\end{remark}

The braided functor $I\circ\Phi \colon \widehat{G}\text{-}\Vect \longrightarrow (\mathcal{C}^\times_G)^G$ gives rise to a commutative $\Delta$-separable symmetric Frobenius algebra~$A_{I\circ \Phi}$ in $(\mathcal{C}^\times_G)^G$ via \Cref{prop:orbifold datum from Frobenius algebra}, from which in turn we obtain an orbifold datum $\mathbb{A}_{A_{I\circ \Phi}}$ in $\mathcal D_{(\mathcal{C}^\times_G)^G}$ via \Cref{prop:FrobeniusAlgebra from 1form symmetry}. 
Thus we have: 

\begin{lemma}[label=lem:ARhat]
	The 1-form symmetry $\widehat R = \mathrm{B}^2(I\circ\Phi)$ gives rise to the orbifold datum
	\begin{equation}
		\mathbb{A}_{\widehat{R}} := \mathbb{A}_{A_{I\circ \Phi}} 
	\end{equation}
	in $\mathcal{D}_{(\mathcal{C}^\times_G)^G}$. 
\end{lemma}

The above orbifold datum~$\mathbb{A}_{\widehat{R}}$ associated to the $1$-form symmetry~$\widehat{R}$ in~\eqref{eq:dual 1-form symmetry} is the same as the orbifold datum~$\mathbb{A}_B$ that was constructed in \Cref{prop:Deequivariantisation as orbifold datum1}, which is part of the main result of this section:

\begin{theorem}[label=thm:inversion of 0form by 1form]
	Let~$G$ be a finite commutative group and let~$\widehat{R}$ be the 1-form symmetry~\eqref{eq:dual 1-form symmetry} for a $G$-crossed braided fusion category~$\mathcal C_G^\times$. 
	Orbifolding its equivariantisation with~$\mathbb{A}_{\widehat{R}}$ extracts the neutral component~$\mathcal C_e$: 
	\begin{equation}
		\left((\mathcal{C}^\times_G)^G\right)_{\mathbb{A}_{\widehat{R}}}
		\;\cong\; 
		\mathcal{C}_e \, .
	\end{equation}
\end{theorem}

\begin{proof}
	Recall that our construction of~$\mathbb{A}_{\widehat{R}}$ for a $1$-form symmetry consists of two steps (\Cref{prop:Orbifold datum from 1form symmetry as composite}): First we construct a commutative separable symmetric Frobenius algebra according to \Cref{prop:FrobeniusAlgebra from 1form symmetry} and then turn it into an orbifold datum using \Cref{prop:orbifold datum from Frobenius algebra}. 
	In order to apply \Cref{prop:orbifold data for compositions of functors}, we observe first that the intermediate Frobenius algebra~$A_{\Phi}$ is precisely $\C[\widehat{G}]$. 
	In fact, each subspace $\Phi(\phi)$ is identified with the ${\C}$-span of~$\phi$, and the algebraic structures follow according to \Cref{prop:FrobeniusAlgebra from 1form symmetry}.
	By \Cref{lem:equivalence of group algebras}, we identify it with $\C(G)$ before applying \Cref{prop:orbifold datum from Frobenius algebra}. 
	Lastly, we apply \Cref{prop:orbifold data for compositions of functors} to transport our orbifold datum into $(\mathcal{C}^\times_G)^G$, obtaining 
	\begin{equation}
		\mathbb{A}_{\widehat{R}} 
		= \big( I(\C(G)),\, I(\C(G)), \,\Delta\circ \mu, \, \id_{B}, \, 1\big)
		= \mathbb{A}_B
	\end{equation}
	Its orbifold category is equivalent to~$\mathcal{C}_e$ by~\eqref{eq:orbifold category of a condensable algebra} and~\eqref{eq:deequivariantisation as local modules}.
\end{proof}

Using the equivalence $F\colon \mathcal{C}_{\mathbb{A}_R} \cong (\mathcal{C}^\times_G)^G$ of ribbon categories of \Cref{thm:eq is orbifolding} below, we can formulate the above result more conceptually. 
Indeed, since~$F$ is pivotal (cf.\ \Cref{lem:Phi is well-defined}), we can compose~$\widehat{R}$ from \eqref{eq:dual 1-form symmetry} with the double delooping $\textrm{B}^2F$ to obtain an orbifold datum in $\mathcal{C}_{\mathbb{A}_R}$. 
By \Cref{prop:orbifold data for compositions of functors}, the resulting orbifold datum~$\mathbb{A}_{F\circ \widehat{R}}$ coincides with the transport of $\mathbb{A}_{\widehat{R}}$ along the equivalence~$F$. 
Due to \cite[Prop.\,5.12]{Mulevicius2022a} and the fact that~$F$ is an equivalence, the resulting orbifold categories are equivalent as well, and we arrive at the fact that one can undo the orbifolding of the $0$-form symmetry by orbifolding the emerging $1$-form symmetry: 

\begin{corollary}
	With the above notation, and recalling $\mathcal C := \mathcal C_e$, we have 
	\begin{equation}
		\big(\mathcal{C}_{\mathbb{A}_R}\big)_{\mathbb{A}_{ F\circ \widehat{R}}} 
		\;\cong\; 
		\mathcal{C} \, .
	\end{equation}
\end{corollary}
	
\subsection{The TQFT Perspective}
\label{subsec:TQFTPerspective}

Here we discuss our algebraic results in terms of defect TQFTs. 

\medskip 

By an $n$-dimensional defect TQFT we mean a symmetric monoidal functor 
\begin{equation}
	\mathcal Z \colon \textrm{Bord}^{\textrm{def}}_{n,n-1}(\mathds D) \longrightarrow \Vect 
\end{equation}
where the domain is a category of labelled stratified bordisms, and~$\mathds D$ consists in particular of sets~$D_j$ to label $j$-dimensional defects for $j\in\{0,\dots,n\}$, see \cite{Carqueville2017} for details. 
If $n=2$ then one extracts a pivotal 2-category~$\mathcal D_{\mathcal Z}$ from such a defect TQFT~$\mathcal Z$, as explained in \cite{Davydov2011}. 
Objects are ``bulk theories'' labelled by~$D_2$, 1-morphisms are line defects labelled by (lists of composable) elements of~$D_1$, and 2-morphism spaces are given by state spaces computed with~$\mathcal Z$. 

It is expected that from an $n$-dimensional defect TQFT~$\mathcal Z$ one obtains an $n$-category~$\mathcal D_{\mathcal Z}$ with coherent adjoints for any value of~$n$, with $k$-morphisms given by $(n-k)$-dimensional defects, i.e.\ (composites of) elements of~$D_{n-k}$ for $k\leqslant n-1$, while $n$-morphism spaces come from state spaces associated to defect $(n-1)$-spheres computed with~$\mathcal Z$. 
For $n=3$ an explicit construction of~$\mathcal D_{\mathcal Z}$ was given in \cite{Carqueville2016} whose result is a Gray category with duals. 
In particular, for the Reshetikhin--Turaev defect TQFT~$\mathcal Z^{\mathcal C}$ \cite{kapustin2010surfaceoperators3dtopological, Carqueville2020} based on a modular fusion category~$\mathcal C$ there is a 3-category~$\mathcal D_{\mathcal Z^{\mathcal C}}$ as explained in \cite{Fuchs_2013, Koppen2021, Carqueville2023}. 
The 3-category~$\mathcal D_{\mathcal C}$ of Example~\ref{exa:3dCategoriesWithAdjoints} is the full sub-3-category of~$\mathcal D_{\mathcal Z^{\mathcal C}}$ that features only one single type of bulk theory, labelled directly by the modular fusion category~$\mathcal C$.\footnote{The other objects in~$\mathcal D_{\mathcal Z^{\mathcal C}}$ can be thought of as those other Reshetikhin--Turaev bulk TQFTs which admit topological interfaces with the Reshetikhin--Turaev model obtained from~$\mathcal C$. 
	More precisely, we have $\mathcal D_{\mathcal Z^{\mathcal C}} \cong (\mathcal D_{\mathcal C})_{\textrm{orb}}$, whose objects can be understood as modular fusion categories made of local modules over commutative $\Delta$-separable symmetric Frobenius algebras in~$\mathcal C$.}

For any defect TQFT~$\mathcal Z$ there is a notion of orbifold datum~$\mathbb{A}$, which can be thought of as the algebraic description of a gaugable (generalised, not necessarily invertible) symmetry. 
It consists of $(n-j)$-morphisms~$\mathbb{A}_j$ in~$\mathcal D_{\mathcal Z}$ that can label defects of codimension~$j$ which are dual to an $(n-j)$-simplex in a triangulation, as well as two $n$-morphisms $\mathbb{A}_0^+, \mathbb{A}_0^-$ (related to the two orientations an $n$-simplex can have). 
The result of gauging the symmetry~$\mathbb{A}$ is the orbifold TQFT~$\mathcal Z_{\mathbb{A}}$ as introduced in \cite{Carqueville2017}. 
Succinctly, $\mathcal Z_{\mathbb{A}}$ is the left Kan extension of the functor which labels stratified bordisms with the orbifold datum~$\mathbb{A}$ and then evaluates those with~$\mathcal Z$ along the functor which forgets triangulations: 
\begin{equation}
	\begin{tikzcd}[column sep= 40,row sep= 40, ampersand replacement=\&]
		\textrm{Bord}^\Delta_{n,n-1}
			\arrow[d]
			\arrow[r, "\mathcal P_{\mathbb{A}}"]
		\&
		\textrm{Bord}^{\textrm{def}}_{n,n-1}(\mathds D)
			\arrow[r, "\mathcal Z"]
		\&
		\Vect \, . 
		\\
		\textrm{Bord}^{\textrm{or}}_{n,n-1}
			\arrow[rru, swap, "\mathcal Z_{\mathbb{A}}"]
		\& \& 
	\end{tikzcd}
\end{equation}
Here $\textrm{Bord}^\Delta_{n,n-1}$ is the category of bordisms that are equipped with oriented triangulations, $\mathcal P_{\mathbb{A}}$ converts triangulations into Poincar\'e dual stratifications that are labelled by~$\mathbb{A}$, and the unnamed vertical arrow simply forgets triangulations. 
Thus~$\mathcal Z_{\mathbb{A}}$ satisfies a universal property. 
In practice, orbifold TQFTs are computed by splitting idempotents, see \cite{Carqueville2017} and the survey \cite{carqueville2023orbifoldstopologicalquantumfield}. 

The orbifold datum~$\mathbb{A}$ in the 3-category~$\mathcal D_{\mathcal C}$ associated to the modular fusion category~$\mathcal C$ (recall Example~\ref{exa:3dOrbifoldData}) is thus a (generalised) symmetry of the Reshetikhin--Turaev theory~$\mathcal Z^{\mathcal C}$. 
Gauging it produces the orbifold TQFT $(\mathcal Z^{\mathcal C})_{\mathbb{A}}$. 
One may wonder whether this is again a TQFT of Reshetikhin--Turaev-type. 
This was conjectured to be the case in \cite{Mulevicius2022}, where the modular fusion category~$\mathcal C_{\mathbb{A}}$ of Remark~\ref{rem:MFCCA} was also constructed. 
As shown in \cite{CMRSS2021, Carqueville2021}, the gauged theory is indeed the Reshetikhin--Turaev TQFT constructed from~$\mathcal C_{\mathbb A}$, i.e.
\begin{equation}
	\label{eq:RTclosed}
	(\mathcal Z^{\mathcal C})_{\mathbb{A}} 
		\;\cong\; 
	\mathcal Z^{\mathcal C_{\mathbb A}} \, . 
\end{equation}

\medskip 

The results of the present paper have a straightforward TQFT interpretation. 
Given a 2-group~$\mathcal G$, a 3-functor $R\colon \textrm{B}\underline{\mathcal G} \longrightarrow \mathcal D_{\mathcal C}$ as in Definition~\ref{def:3d 2-group symmetry} describes a symmetry of the Reshetikhin--Turaev TQFT~$\mathcal Z^{\mathcal C}$ (viewed as the single object of $\mathcal D_{\mathcal C}$). 
Since~$R$ is a 3-functor, the constituents of~$\mathcal G$ act via surface and line defects of~$\mathcal Z^{\mathcal C}$, which in turn satisfy fusion rules up to higher morphisms. 
To this we naturally associate the associative algebra~$\mathbb A_R$ in~$\mathcal D_{\mathcal C}$ of Lemma~\ref{lem:0form orb dat from Gcbrc}. 
The $\mathcal G$-symmetry is \textsl{gaugable} iff~$\mathbb A_R$ is an orbifold datum. 

In the special case that~$R$ is a 0-form symmetry obtained from a $G$-crossed braided fusion category~$\mathcal C^\times_G$ (and assuming its neutral component $\mathcal{C}_e$ is modular), we obtain the orbifold datum~$\mathbb A_R$ in~$\mathcal D_{\mathcal C_e}$ as described in Section~\ref{sec:equivariantisation is orbifolding}. 
As we prove in Theorem~\ref{thm:eq is orbifolding} below, gauging the Reshetikhin--Turaev TQFT~$\mathcal Z^{\mathcal C_e}$ by this symmetry results in the Reshetikhin--Turaev TQFT associated to the equivariantisation of~$\mathcal C^\times_G$: 
\begin{equation}
	\label{eq:TQFTIso}
	\mathcal Z^{(\mathcal C^\times_G)^G} 
		\;\cong\; 
		\mathcal Z^{(\mathcal C_e)_{\mathbb A_{\scaleto{R}{3pt}}}} \, . 
\end{equation}

After gauging, there is a Pontryagin dual 1-form symmetry $\widehat R\colon \textrm{B}^2(\widehat G\text{-}\Vect) \longrightarrow \mathcal D_{(\mathcal C^\times_G)^G}$ which gives rise to the orbifold datum~$\mathbb A_{\widehat R}$ of Section~\ref{sec:Deeq is 1-form orbifolding}. 
According to Theorem~\ref{thm:inversion of 0form by 1form}, gauging this symmetry undoes the equivariantisation if~$G$ is commutative. If~$G$ is not commutative, then the $\Rep(G)$-symmetry is non-invertible and therefore not a 1-form symmetry, but we can still construct an orbifold datum which undoes the equivariantisation, see \eqref{eq:non-commutative deequivariantisation via orbifold}.
Hence together with~\eqref{eq:RTclosed} and~\eqref{eq:TQFTIso} we proved that 
\begin{equation}
	\big( (\mathcal Z^{\mathcal C_e})_{\mathbb A_R} \big)_{\mathbb A_{\widehat R}}
		\;\cong\; 
	\mathcal Z^{\mathcal C_e} \, . 
\end{equation}
Gauging dual 0- and 1-form symmetries of Reshetikhin--Turaev TQFTs coming from $G$-crossed braided fusion categories is mutually inverse.

\section{Equivalence of Equivariantisation and Orbifolding}
\label{sec:equivalence of equivariantisation and orbifolding}

In this section, we establish the equivalence between equivariantisation and orbifolding of the $0$-form symmetry~$R$ from a $G$-crossed braided fusion category~$\Gcbc$, where~$G$ is an arbitrary finite group. 
In \Cref{subsec:ingredients} we briefly collect some tools for dealing with balanced morphisms that are used throughout the rest of the paper. 
In \Cref{sec:equivalence of Eq and Orb} we provide all the data necessary to build an adjoint equivalence between $(\Gcbc)^G$ and $\CAR$ as braided pivotal fusion categories. 
In \Cref{sec:proof of Equivalence} we prove that this indeed is an equivalence. 

\subsection{Balanced Morphisms in the Orbifold Category}
\label{subsec:ingredients}

\Cref{lem:balancedmorphism composition} and the subsequent observation (balanced morphisms absorb the idempotent of relative tensor products) allow us to ``cut and glue'' strings in our diagrams. 

\begin{lemma}[label=lem:balancedmorphism composition]
	Let~$A$ be a Frobenius algebra in a monoidal category, let~$X_2$ and~$Y_1$ be left $A$-modules, $X_1$ a right $A$-module, and let $Y_2, Z\in\mathcal C$. 
	Let further $f\colon X_2\longrightarrow Y_1\otimes Y_2$ be a module morphism and $f^\prime\colon X_1\otimes Y_1\longrightarrow Z$ a balanced morphism. 
	Then 
	\begin{equation}
		(f^\prime\otimes \id_{Y_1})\circ (\id_{X_1}\otimes f)=
 
		\quad \text{is balanced on $X_1\otimes X_2$.}
	\end{equation}
\end{lemma}

\begin{proof}
	This follows immediately from the fact that~$f^\prime$ is balanced, and that~$f$ is a module morphism: 
	\begin{equation}			
		\label{eq:balanced morphisms from composition}
		%
 
		\, . 
	\end{equation}
\end{proof}

Recall that a balanced morphism $f^\prime\colon X\otimes Y\longrightarrow Z$ ``absorbs'' the idempotent \eqref{eq:relative tensor idempotent} whose splitting gives the relative tensor product: 
\begin{equation}
	\label{eq:idempotent absorption by balanced morphisms}
	%
	\, ,
\end{equation}
and the analogously cobalanced morphisms.

Recall the orbifold datum~$\mathbb{A}_R$ from \Cref{lem:0form orb dat from Gcbrc} and the definition of~$\mathcal{C}_{\mathbb{A}_R}$ from \Cref{rem:MFCCA} as a rigid braided monoidal category. 
In the notation of \Cref{sec:equivariantisation is orbifolding}, the idempotent for $m_g\otimes_{A_g}m_g^*$ is given by 
\begin{equation}
	e_{m_g,A_g,m_g^*} := 
	%
 
	\, . 
\end{equation}
Hence for a balanced morphism $f^\prime\colon m_g\otimes m_g^*\longrightarrow X$, we can ``cut and glue'':
\begin{equation}
	%
 
	\, .
\end{equation}

\subsection{Ingredients of the Equivalence}
\label{sec:equivalence of Eq and Orb}

In this section, we spell out the data for the adjoint equivalence between equivariantisation and orbifolding of the $0$-form $G$-symmetry for $G$-crossed braided fusion categories. 

\medskip

Let~$\mathcal{C}^\times_G$ be a $G$-crossed braided fusion category. 
Let $\pr_1\colon \mathcal{C}^\times_G\longrightarrow \mathcal{C}_e$ be the functor which projects onto the part graded by the unit $e\in G$ which we denote as 
\begin{equation}
	\mathcal{C}\equiv\mathcal{C}_e	
\end{equation}
for simplicity. 
Recall the orbifold datum~$\mathbb{A}_R$ from \Cref{lem:0form orb dat from Gcbrc}, featuring in particular the underlying algebra $A=\bigoplus_{g\in G} A_g$ and the isomorphism~$\psi$ of~\eqref{eq:psiInAR}. 
For an object $(X, u)\equiv (\bigoplus_{g\in G}X_g,\{u_h\}_{h\in G})$ in the equivariantisation $(\mathcal{C}^\times_G)^G$, we set
\begin{align}
	\label{eq:FonX}
	 F(X)&:=\pr_1 \Big(\bigoplus_{g,h\in G} m_g^*\otimes X\otimes m_h\Big)
	 	=
	 	\bigoplus_{g,h\in G}m_{gh}^*\otimes X_g\otimes m_h \, , 
\end{align}
which carries the structure of an $(A,A)$-bimodule via \eqref{eq:mmodulestructures}. We equip it with the additional $\tau$-structure (writing $\bar{h}$ instead of $h^{-1}$)
\begin{align}
	\tau_1^{ F(X)}&:=\bigoplus_{g,h,k\in G}
\, .
	\label{eq:tau structure on Phi}
\end{align}
The monoidal unit in $\mathcal{C}_{\mathbb{A}_R}$ is $\mathds{1}_{\mathcal{C}_{\mathbb{A}_{\scaleto{R}{3pt}}}}\equiv (A, \tau_i^{\mathds{1}})$ with the $\tau$-structure given by 
\begin{align}
	\tau_1^{\mathds{1}}&=\bigoplus_{g,h\in G}
	%
\, ,
	\label{eq:tau structure of monoidal unit in CAR}
\end{align}
cf.\ \cite[Eq.\,(3.5)]{Mulevicius2022}. 
For $(X,u)=(\mathds{1},\{\rho^0_g\})$, the monoidal unit in~$(\mathcal{C}^\times_G)^G$, we have 
\begin{equation}
	F(\mathds{1},\{\rho^0_g\}) 
		:= 
		\Big(A,\tau_i^{ F(\mathds{1})},\overline{\tau}_i^{ F(\mathds{1})}\Big)
\end{equation}
whose $\tau$-structure is given by that of $\mathds{1}_{\mathcal{C}_{\mathbb{A}_{\scaleto{R}{3pt}}}}$ \eqref{eq:tau structure of monoidal unit in CAR} without the $\psi$-insertions.
We can then define morphisms
\begin{align}
	 F^2_{X,X^\prime}&:=\bigoplus_{g,h,k\in G}
	\begin{tikzpicture}[very thick,color=red!50!black, baseline=1cm]
		\draw[line width=0pt] 
		(2.5,0) node[below] (D) {{\scriptsize$hk\vphantom{ghkX_gX^\prime_h}$}}
		(4,0) node[below] (s) {{\scriptsize$hk\vphantom{ghkX_gX^\prime_h}$}}; 
		\draw[redirectedred] (s) .. controls +(0,1) and +(0,1) .. (D)node[pos=0.3, line width=-4pt](u){};
		\draw[line width=0pt] 
		(5,0) node[below] (re) {{\scriptsize$k\vphantom{ghkX_gX^\prime_h}$}}
		(1.5,0) node[below] (li) {{\scriptsize$ghk\vphantom{ghkX_gX^\prime_h}$}}; 
		\draw[line width=0pt, color=blue!50!black]
		(3,2) node[above] (oX) {{\scriptsize$X_g\vphantom{ghkX_gX^\prime_h}$}}
		(3.5,2) node[above] (oY) {{\scriptsize$X^\prime_h\vphantom{ghkX_gX^\prime_h}$}}
		(2,0) node[below] (uX) {{\scriptsize$X_g\vphantom{ghkX_gX^\prime_h}$}}
		(4.5,0) node[below] (uY) {{\scriptsize$X^\prime_h\vphantom{ghkX_gX^\prime_h}$}}; 
		\draw[line width=0pt] 
		(4,2) node[above] (ore) {{\scriptsize$k\vphantom{ghkX_gX^\prime_h}$}}
		(2.5,2) node[above] (oli) {{\scriptsize$ghk\vphantom{ghkX_gX^\prime_h}$}}; 
		\draw[directedred] (oli) .. controls +(0,-1.5) and +(0,1) .. (li);
		\draw[redirectedred] (ore) .. controls +(0,-1.5) and +(0,1) .. (re);
		\draw[color=blue!50!black] (uX) .. controls +(0,1.25) and +(0,-1.75) .. (oX)
		(uY) .. controls +(0,1.25) and +(0,-1.75) .. (oY);
		\fill (u) circle (2.3pt)
		($(u)+(-.05,-.3)$) node {{\scriptsize$\psi$}};
	\end{tikzpicture}\, , & F^0&:=
	\begin{tikzpicture}[very thick,color=red!50!black, baseline=1cm]
		\draw[line width=0pt] 
		(3.4,0) node[line width=0pt] (ho) {{\scriptsize$g\vphantom{ghk}$}}
		(3,0) node[line width=0pt] (go) {{\scriptsize$g\vphantom{ghk}$}}
		(3.4,2) node[line width=0pt] (hu) {{\scriptsize$g\vphantom{ghk}$}}
		(3,2) node[line width=0pt] (gu) {{\scriptsize$g\vphantom{ghk}$}};
		\draw[redirectedred] (hu) -- (ho);
		\draw[redirectedred] (go) -- (gu);
		\fill (3,1.3) circle (2.3pt)
		(2.6,1.3) node {{\scriptsize$\psi^{-1}$}};
	\end{tikzpicture}
	\, .
	\label{eq:monoidal structure of Phi}
\end{align}

\begin{lemma}[label=lem:Phi is well-defined]
	The assignment
	\begin{align}
		F\colon (\mathcal{C}^\times_G)^G&\longrightarrow \mathcal{C}_{\mathbb{A}_R}
		\nonumber
		\\
		(X,u)&\longmapsto \Big( F(X),\tau_1^{ F(X)},\tau_2^{ F(X)},\overline{\tau}_1^{ F(X)},\overline{\tau}_2^{ F(X)}\Big)
		\nonumber
		\\
		\big(f\colon (X,u)\longrightarrow (X^\prime, u^\prime)\big) &\longmapsto \pr_1\Big(\bigoplus_{g,h\in G} \id_{m_g^*}\otimes f\otimes \id_{m_h}\Big)
	\end{align}
	together with the monoidal structure $(F^2,F^0)$ of \eqref{eq:monoidal structure of Phi} is a ribbon functor.
\end{lemma}

\begin{proof}
	We first show that $ F$ is well-defined. 
	The axioms which involve $\tau_1^{ F(X)}$ and $\overline{\tau}_1^{ F(X)}$ are satisfied as argued at \cite[Eq.\,(5.9)]{Mulevicius2022a}. 
	There, it is shown that for a left $A$-module~$L$, the object $T\otimes_1 L$ can be equipped with morphisms $\tau_1$ and $\overline{\tau}_1$ which satisfy the required axioms for an object in $\CAR$, and their construction coincides with ours in the case of $L=\bigoplus_{g\in G}m_g^*\otimes X_g$ up to canonical isomorphism given by the braiding. 
	The axioms for the other $\tau$-maps that~$F$ assigns to $(X,u) \in (\mathcal{C}^\times_G)^G$ follow analogously: 
	the conditions \eqref{eq:T2}, \eqref{eq:T4}, \eqref{eq:T5}, \eqref{eq:T6} and \eqref{eq:T7} in \Cref{fig:tau-structure axioms} can be shown by simply inserting our definitions \eqref{eq:string alphaG} and \eqref{eq:tau structure on Phi} and cancelling loops with $\psi^2$-insertions, while for \eqref{eq:T3}, we additionally employ the defining property of equivariantisations, \eqref{eq:defining property of equivariant structure}. 
	
	On morphisms, compatibility with $\tau_1$-maps follows from naturality of the braiding, whereas compatibility with $\tau_2$-maps it follows from the defining property of morphisms in the equivariantisation \eqref{eq:defining property of morphisms in equivariantisation}, i.e.\ they are compatible with the equivariant structures. 
	Lastly, $ F$ is clearly compatible with composition and maps identity morphisms to identity morphisms. 
	Hence~$F$ is indeed a functor. 
	
	Next we discuss the monoidal structure. We identify $ F(\mathds{1})$ with~$A$, i.e.\ blue lines in previous diagrams will be omitted when they would be labelled by~$\mathds{1}$.
	The inverses of the monoidal structure morphisms are given by 
	\begin{align}
		( F^2_{X,X^\prime})^{-1} :=\bigoplus_{g,h,k\in G}
\, .
		\label{eq:inverse monoidal structure of Phi}
	\end{align}
	Note that $ F^0\colon (A,\tau^\mathds{1}_i)\longrightarrow (A,\tau^{ F(\mathds{1})}_i)$ is a morphisms in $\CAR$ since it satisfies \eqref{eq:morphisms in CAA}:
	\begin{equation}
		\tau^{ F(\mathds{1})}_i \circ ( F^0\otimes_A \id_T) = (\id_T \otimes_i  F^0)\circ\tau^{\mathds{1}}_i\, .
	\end{equation} 
	On the right-hand side of this equation, $\psi$ in $\tau^{\mathds{1}}_i$ and $\psi^{-1}$ from $ F^0$ cancel, while on the left-hand side a $\psi^{-1}$ is introduced by $ F^0$, matching the remaining $\psi^{-1}$-insertion in~$\tau^{\mathds{1}}_i$. 
	For $F^2$ we need to show that it is compatible with the  $\tau$-structures $\tau^{ F(X), F(X^\prime)}_i$ from \eqref{eq:CAA monoidal product tau} and $\tau^{ F(X\otimes X^\prime)}_i$, i.e.\ 
	\begin{equation}
		\tau^{ F(X\otimes X^\prime)}_i \circ ( F^2_{X,X^\prime}\otimes_A \id_T) = (\id_T \otimes_i  F^2_{X,X^\prime})\circ\tau^{ F(X), F(X^\prime)}_i \, . 
	\end{equation}
	This follows by spelling out the definitions and cancelling loops with $\psi^2$-insertions, similar to the calculation for the braiding which we show in \eqref{eq:braiding in CAR}.
	The unitality axioms for the monoidal structure follow by observing that the $\psi$ in~$ F^2$ cancels the~$\psi^{-1}$ in~$ F^0$ and the result is the unitor in $\CAR$, which is given by \eqref{eq:mmodulestructures} up to additional monoidal products with identity morphisms. 
	The associativity axiom is immediate when drawing the diagrams. 
	Note also that~$ F^2$ is natural, since $ F(f)\otimes_A F(f^\prime)= F(f\otimes f^\prime)$ under the identification $m_{hk}\otimes_{A_{hk}}m_{hk}^*=\mathds{1}$.
	
	Lastly, we show that~$ F$ is braided and pivotal. 
	For the braiding, we need to show 
	\begin{equation}
		c^{\CAR}_{ F(X), F(X^\prime)}=( F^2_{X^\prime,X})^{-1}\circ  F(c^{G}_{X,X^\prime})\circ  F^2_{X,X^\prime}
	\end{equation}
	for which we recall the braiding~$c^G$ of \Cref{lem:eq of Gcbc is braided}.
	Inserting our definitions into $c^{\CAR}_{ F(X), F(X^\prime)}$ from \eqref{eq:CAA braiding} and moving all $T$-crossings to the left side of the loop, this reads
	\begin{equation}
		\label{eq:braiding in CAR}
		\bigoplus_{g,h,k,j\in G}\frac{1}{|G|}
\, ,
	\end{equation}
	where $x:=\bar{k}\bar{h}ghkj$. 
	The actions of $\rho^{-1}$ appear due to $G$-crossed naturality of the braiding when removing the loops labelled $x$ and $hk$.
	Note that by \eqref{eq:defining property of equivariant structure}, we have $u_k\circ\rho_k(u_h)=u_{hk}\circ(\rho^2_{k,h})_X$ which we can rearrange to $u_h=\rho_k^{-1}(u_k^{-1}\circ u_{hk})$, where we ignored~$\rho^2$, since it is suppressed in diagrams. 
	Then, by $G$-crossed naturality of the braiding, we can identify $ F(c^{G}_{X,X^\prime})$ defined in \eqref{eq:braiding on equivariantisation} in the centre of the diagram on the right-hand side. 
	The complete diagram thus depicts $( F^2_{X^\prime,X})^{-1}\circ  F(c^{G}_{X,X^\prime})\circ  F^2_{X,X^\prime}$ where the sum over~$j$ cancels the factor $|G|^{-1}$. 
	Hence $( F, F^2, F^0)$ is braided.
	
	To show pivotality, note first that duals in the equivariantisation are the duals in the underlying category with equivariant structure $u_g^{X^*}=(u_g^{-1})^*$. 
	In $\CAR$, the underlying object is again the dual in $\mathcal{C}$, and recall that the $\tau$-structure on duals is given by \eqref{eq:dual tau}. 
	We see immediately that $ F(X)^*= F(X^*)$ on the level of objects and similarly for morphisms. 
	Inserting the dual equivariant structures into \eqref{eq:tau structure on Phi} recovers \eqref{eq:dual tau} exactly. 
	The (left) evaluation morphism on $ F(X^*)$ is given by the composition $( F^0)^{-1}\circ F(\ev_X)\circ  F^2_{X^*,X}$ which reduces to the evaluation in $\CAR$, with the $\psi$-insertions coming from the monoidal data. 
	This holds similarly for the coevaluation and right adjunction data. 
	Therefore, since~$ F$ preserves left and right duals, it is pivotal.
\end{proof}

We proceed with a description of the building blocks for the inverse of~$F$, which we denote~$E$ as it reconstructs the equivariantisation. 
Given objects $(Y,\tau)$ and $(Y^\prime,\tau^\prime)$ in~$\mathcal{C}_{\mathbb{A}_R}$ and a bimodule morphism $f\colon Y\longrightarrow Y^\prime$, we set
\begin{align}
	Y_{g,h}&:=m_g\otimes_{A_g}Y\otimes_{A_h} m_h^*\quad \in\;\mathcal{C}^\times_G \,,
	\label{eq:Ygh}
	\\
	f_{g,h}&:=m_g\otimes_{A_g}f\otimes_{A_h} m_h^*\colon Y_{g,h}\longrightarrow Y_{g,h}^\prime\,.
	\label{eq:relative tensor products of morphisms in CAR}
\end{align}
Next, we define 
\begin{equation}
	\gamma_{(Y,\tau)}^{g,h,k}\colon Y_{gk,hk}\longrightarrow Y_{g,h}
\end{equation}
to be given by the balanced morphism
\begin{equation}
	\label{eq:component def of gamma}
	\gamma_{(Y,\tau)}^{g,h,k}:=	
	\begin{tikzpicture}[very thick,scale=0.6,color=red!50!black, baseline]
		\coordinate (Yt) at (1.05, 3);
		\coordinate (ht) at (1.7, 3);
		\coordinate (gt) at (0.4, 3);
		\coordinate (Y) at (-1.05, -3);
		\coordinate (gk) at (-1.7, -3);
		\coordinate (hk) at (-0.4, -3);
		\draw[draw=blue!50!black,very thick] (Y) -- (-1.05, -0.5)
		(Yt) -- (1.05, 0.5);
		\draw[directedred] (gk) .. controls +(0,1) and +(0,-2) .. (-2,0.5);
		\draw[directedred] (2,0.5) -- (2,-0.5);
		\draw[redirectedred] (0.4,2.4) .. controls +(0,-1) and +(0,1) .. (-0.75,0.5);
		\draw[directedred] (ht) .. controls +(0,-1.5) and +(0,1.5) .. (3,-0.5);
		\draw[directedred] (0.3,-0.5) .. controls +(0,-1.5) and +(0,1.5) .. (hk);
		\draw[draw=white,line width=4pt] (2,0.5) .. controls +(0,1) and +(0,1) .. (-0.3,0.5);
		\draw
		(-2,0.5) .. controls +(0,0.65) and +(0,0.65) .. (-1.2,0.5)
		(2,-0.5) .. controls +(0,-0.65) and +(0,-0.65) .. (1.2,-0.5)
		(2,0.5) .. controls +(0,1) and +(0,1) .. (-0.3,0.5)
		(0.4,2.4) -- (0.4,3);
		\draw[directedred] (3,-0.5) .. controls +(0,-1.3) and +(0,-1.3) .. (0.75,-0.5);
		\draw[draw=black, line width=0.5pt, rounded corners=2pt](-1.5,-.5)rectangle(1.5,.5);
		\draw (1.8,0.45) node {{\scriptsize$k$}};
		\draw[color=black] (0,0) node (tau) {{\scriptsize$\tau_1$}};
		\fill[color=blue!50!black]
		(Y) node[below] {{\scriptsize$Y\vphantom{gYhk}$}}
		(Yt) node[above] {{\scriptsize$Y\vphantom{gYhk}$}};
		\fill[color=red!50!black]
		(gk) node[below] {{\scriptsize$gk\vphantom{gYhk}$}}
		(hk) node[below] {{\scriptsize$hk\vphantom{gYhk}$}}
		(gt) node[above] {{\scriptsize$g\vphantom{gYhk}$}}
		(ht) node[above] {{\scriptsize$h\vphantom{gYhk}$}};
		\foreach \x/\y in {2/-0.5,3/-0.5,0.4/2.4}  
		{\fill(\x,\y)  circle (2.9pt);
			\draw ($(\x,\y)+(0.38,0)$) node {{\tiny$\psi^2$}};}
	\end{tikzpicture} 
	\, ,
\end{equation}
where by abuse of notation here and in the following, the coupon labelled~$\tau_1$ denotes the restriction $\tau_1\vert_{Y\otimes_{A_hk}T_{h,k}}\colon Y\otimes_{A_hk}T_{h,k} \longrightarrow T_{g,k}\otimes_{A_{gk}}Y$.

\begin{lemma}[label=lem:gamma is natural]
	The morphisms \eqref{eq:component def of gamma} assemble into a natural isomorphism
	\begin{equation}
		\gamma^{g,h,k}\colon m_{gk}\otimes_{A_{gk}}(-)\otimes_{A_{hk}}m_{hk}^*\Longrightarrow m_g\otimes_{A_g}(-)\otimes_{A_h} m_h^*
	\end{equation}
	for each triple of group elements $g,h,k\in G$, where we understand the functors calculating relative tensor products as functors $\mathcal{C}_{\mathbb{A}_R}\longrightarrow \mathcal{C}^\times_G$. 
	Moreover, the family of these natural isomorphisms is compatible with the group multiplication in the sense that
	\begin{equation}
		\label{eq:compatibility of gamma with group multiplication}
		\gamma^{g,h,k_2k_1}=\gamma^{g,h,k_2}\circ\gamma^{gk_2,hk_2,k_1}\,.
	\end{equation}
	\begin{proof}
		We claim that the inverse is given by the natural transformation with components
		\begin{equation}
			\label{eq:inverse of gamma}
			(\gamma^{g,h,k})^{-1}_{(Y,\tau)}:=(\gamma_{(Y,\tau)}^{g,h,k})^{-1}=	
 \, .
		\end{equation}
		Recall that the full morphisms are given by the composition of these diagrams with retracts and sections from \eqref{eq:RetractionSection} which are balanced and cobalanced, respectively. 
		Thus, omitting the section and the coevaluations from $\gamma^{g,h,k}$ leaves us with
		\begin{equation}
			%
 \, .
		\end{equation}
		Note that this is balanced on $(m_h\otimes m_k)\otimes (m_k^*\otimes m_h^*)$ with respect to $A_h\otimes A_k$ by \Cref{lem:balancedmorphism composition}.
		To prove that \eqref{eq:inverse of gamma} is the inverse, we consider the composition $\gamma_{(Y,\tau)}^{g,h,k}\circ (\gamma_{(Y,\tau)}^{g,h,k})^{-1}$ and use \eqref{eq:idempotent absorption by balanced morphisms} twice: 
		\begin{equation}
			%
 \, .
		\end{equation}
		In the last step, the~$\psi^2$ on~$m_g$ was cancelled by the~$\psi^{-2}$ coming from the defining condition \eqref{eq:T5}. 
		The remaining contributions from loops are cancelled by their $\psi^2$-insertions. Thus, we obtain $\gamma_{(Y,\tau)}^{g,h,k}\circ (\gamma_{(Y,\tau)}^{g,h,k})^{-1}=\id_{Y_{g,h}}$, and a similar calculation shows that the composition in the opposite order also reduces to the identity (using \eqref{eq:T4}). 
		The transformation and its inverse are natural, since morphisms in $\mathcal{C}_{\mathbb{A}_R}$ commute with~$\tau_1$ and~$\overline{\tau}_1$.
		
		To prove \eqref{eq:compatibility of gamma with group multiplication}, we use \eqref{eq:T1} along similar lines.
		As before, we are left with loops carrying $\psi^2$-insertions which cancel once again. 
	\end{proof}
\end{lemma}

Using the maps $\gamma^{g,h,k}_Y$, we define an equivariant structure $u^{(Y,\tau)}$ on
\begin{equation}
	\label{eq:YGe}
	Y_{G,e}:=\bigoplus_{g\in G}Y_{g,e} \, , 
\end{equation}
whose components are given by
\begin{align}
	u^{(Y,\tau)}_h&:=\bigoplus_{g\in G}
 \, ,
	\label{eq:equivariant structure on preimage}
\end{align}
where we write $\bar{h}$ for $h^{-1}$, as before.

\begin{lemma}[label=lem:equivariant structure on Y is well-defined]
	The pair $(Y_{G,e}, u^{(Y,\tau)})$ is an equivariant object, i.e.\ an object of $(\Gcbc)^G$.
\end{lemma}

\begin{proof}
	The fact that the right-hand side of \eqref{eq:equivariant structure on preimage} is the inverse of the left-hand side follows from the fact that the coupons labelled by $\gamma$ and its inverse cancel each other, the remaining parts are inverses of each other by the same argument as in the proof of \Cref{lem:gamma is natural}. 
	
	It remains to check the defining property~\eqref{eq:defining property of equivariant structure} of an equivariant structure, i.e.\ 
	\begin{equation}
		u^{(Y,\tau)}_g\circ\rho_g(u^{(Y,\tau)}_h)=u^{(Y,\tau)}_{hg}\circ(\rho^2_{g,h})_X\,. 
	\end{equation}
	Note that if we insert the string diagram for $\gamma^{\bar{h}gh,e,\bar{h}}_Y$, then the composition $u^{(Y,\tau)}_g\circ\rho_g(u^{(Y,\tau)}_h)$ contains~$\tau_1$ and~$\tau_2$ twice each, in alternating order. 
	We apply the identity \eqref{eq:T2} to swap the order of the two morphisms in the middle, then we use \eqref{eq:T3} on the two~$\tau_2$ at the bottom and \eqref{eq:T1} on the two~$\tau_1$ at the top, leaving us with the desired result containing one copy of each. 
	In doing so, loops with $\psi^2$-insertions have to be added and removed to create the necessary copies of~$\alpha$.
\end{proof}

Using the notation of~\eqref{eq:Ygh}, \eqref{eq:relative tensor products of morphisms in CAR} and~\eqref{eq:YGe}, we have: 

\begin{lemma}
	The following assignment is a functor: 
	\begin{align}
		 E\colon  \mathcal{C}_{\mathbb{A}_R}&\longrightarrow (\mathcal{C}^\times_G)^G
		\nonumber\\
		(Y,\tau_i,\overline{\tau}_i)&\longmapsto \big(Y_{G,e},\, u^{(Y,\tau)}\big)
		\nonumber\\
		\big(f\colon Y\longrightarrow Y^\prime\big) &\longmapsto f_{G,e} :=\bigoplus_{g\in G} f_{g,e} \, . 
	\end{align} 
\end{lemma}

\begin{proof}
	\Cref{lem:equivariant structure on Y is well-defined} states that~$E$ is well-defined on objects. 
	For a morphism $f\colon Y\longrightarrow Y^\prime$ in $\mathcal{C}_{\mathbb{A}_R}$, we have $u^{(Y^\prime,\tau^\prime)}_g\circ \rho_g(f_{G,e})=f_{G,e}\circ u^{(Y,\tau)}_g$. 
	The $G$-crossed naturality of the braiding removes the action~$\rho_g$, and since~$f$ is a morphism in $\mathcal{C}_{\mathbb{A}_R}$, it commutes with~$\tau_2$. 
	Therefore, $f_{G,e}$ is a well-defined morphism in $(\Gcbc)^G$. 
	Clearly, $E$ preserves composition and identity morphisms.
\end{proof}

Lastly, we describe the evaluation and coevaluation morphisms for the adjunction between~$E$ and~$F$.  
Identifying $m_g\otimes_{A_g}m_g^*=\mathds{1}$ as before, we have $ E( F(X,u))=(X,u^{ F(X,u)})$ and similarly on morphisms, $ E( F(f))=f$. 
Hence we define a natural isomorphism $\eta\colon \id_{(\Gcbc)^G}\Longrightarrow  E\circ F$ whose components are morphisms in $(\Gcbc)^G$ whose underlying morphisms in $\Gcbc$ are all identity morphisms. 
To show that these are morphisms $(X,u)\longrightarrow (X,u^{ F(X,u)})$ in $(\Gcbc)^G$, one checks that $u^{ F(X,u)}_g=u_g$ for all $g\in G$.

To show that~$\eta$ is a coevaluation witnessing the adjunction $F \dashv E$, we provide a natural isomorphism $\epsilon\colon F\circ E\Longrightarrow \id_{\mathcal{C}_{\mathbb{A}_{\scaleto{R}{3pt}}}}$ which is the associated evaluation. 
Its components are given by
\begin{equation}
	\label{eq:counit in adjunction}
	\epsilon_{(Y,\tau)}:=\bigoplus_{g,h\in G}
	\begin{tikzpicture}[very thick,scale=0.6,color=blue!50!black, baseline=.2cm]
		\coordinate (Yt) at (.5, 3.5);
		\coordinate (X1) at (-1, -2.5);
		\coordinate (X2) at (0.5, -2.5);
		\coordinate (X3) at (2, -2.5);
		\coordinate (s1) at (.5, 3);
		\coordinate (s2) at (.5, 2.6);
		\draw
		(X2) -- (0.5, -1.5)
		(Yt) -- (0.5, -0.5);
		\draw[directedred] (-.5,.5)--(-.5,1.5);
		\draw[redirectedred] (1.5,.5)--(1.5,1.5);
		\draw[directedred] (-1,1.5)--(X1);
		\draw[redirectedred] (2,1.5)--(X3);
		\draw[color=green!50!black] 
		(-.75,1.5) .. controls +(0,.5) and +(-.2,-.2) .. (s2)
		(1.75,1.5) .. controls +(0,.5) and +(.3,-.3) .. (s1)
		(-.75,2.4) node {{\scriptsize$A_{gh}$}}
		(1.75,2.4) node {{\scriptsize$A_h$}};
		\draw[color=red!50!black, dotted] 
		(-1,1.5) -- (-.5,1.5)
		(2,1.5) -- (1.5,1.5);
		\draw[draw=black, line width=0.5pt, rounded corners=2pt](1.5,-.5)rectangle(-0.5,-1.5);
		\draw[color=black] (1.8,.5) -- (-.8,.5);
		\draw[color=black] (0.6,-1) node (Y1) {{\scriptsize$\gamma_{(Y,\tau)}^{gh,h,\bar{h}}$}};
		\fill[color=blue!50!black]
		(X2) node[below] {{\scriptsize$Y_{g,e}\vphantom{gh}$}}
		(Yt) node[above] {{\scriptsize$Y\vphantom{Y_2Z}$}}
		(1.14,0) node (Y1) {{\scriptsize$Y_{gh,h}$}}
		(s1) circle (3.9pt)
		(s2) circle (3.9pt);
		\draw[color=red!50!black]
		(X1) node[below] {{\scriptsize$gh\vphantom{Y_{g,e}}$}}
		(X3) node[below] {{\scriptsize$h\vphantom{Y_{g,e}}$}}
		(0.06,1) node {{\scriptsize$gh$}}
		(1.1,1) node {{\scriptsize$h$}};
	\end{tikzpicture} \, .
\end{equation}%
Note that the upper third of this diagram involving the bimodule structure of~$Y$ is the retraction of $A_{gh}\otimes Y\otimes A_h$ onto $Y\cong A_{gh}\otimes_{A_{gh}}Y\otimes_{A_h} A_h$. 
The fact that the components of~$\epsilon$ are well-defined morphisms in $\CAR$ (i.e.\ they are compatible with the $\tau$-structures) follows by the $\tau$-axioms \eqref{eq:T1} and \eqref{eq:T2} for~$\tau_1$ and~$\tau_2$, respectively. 
That~$\epsilon$ is natural follows from naturality of $\gamma^{g,h,k}$ for all $g,h,k\in G$ and the fact that morphisms in $\CAR$ are bimodule morphisms. 
\begin{proposition}
	The data $(F,E,\eta,\epsilon)$ form an adjunction
	\begin{equation}
		\begin{tikzcd}[column sep=3em, >=stealth]
			(\mathcal{C}^\times_G)^G   \arrow[rr, pos=0.5, out=20, in=160, "F"] 
			& \perp & 
			\mathcal{C}_{\mathbb{A}_R} \, .   \arrow[ll, out=-160, in=-20, "E"] 
		\end{tikzcd}
	\end{equation}
\end{proposition}

\begin{proof}
	We have to check that~$\eta$ and~$\epsilon$ satisfy zigzag identities. 
	This is straightforward because~$\eta$ only consists of identity morphisms. 
	By \Cref{lem:identification of relative tensor products of morphisms} below, $\gamma_{ F(X,u)}^{g,h,k}$ and therefore $\epsilon_{ F(X,u)}$ also reduce to identities, and we obtain $\epsilon F\circ  F\eta=\id_ F$. 
	Similarly, by \eqref{eq:compatibility of gamma with group multiplication} we have $\gamma^{g,e,h}=\gamma^{g,e,e}\circ\gamma^{g,e,h}$, and hence $\gamma_{(Y,\tau)}^{g,e,e}=\id_{Y_{g,e}}$. 
	It follows that $ E(\epsilon_{(Y,\tau)})=(\epsilon_{(Y,\tau)})_{G,e}=\bigoplus_{g\in G}\gamma_{(Y,\tau)}^{g,e,e}=\id_{Y_{G,e}}$, and we obtain the other zigzag identity $E \epsilon\circ\eta E=\id_ E$.
\end{proof}

It is now straightforward to arrive at the main result of \Cref{sec:equivalence of equivariantisation and orbifolding}:

\begin{theorem}[label=thm:eq is orbifolding]
	Let~$\mathcal{C}^\times_G$ be a $G$-crossed braided fusion category. 
	There is an equivalence of ribbon categories between the equivariantisation $(\mathcal{C}^\times_G)^G$ and the fusion category~$\mathcal{C}_{\mathbb{A}_R}$ associated to the orbifold datum~$\mathbb{A}_R$ (of \Cref{lem:0form orb dat from Gcbrc}) given by the functor
	\begin{equation}
		 F\colon \left(\mathcal{C}^\times_G\right)^G\stackrel{\cong}{\longrightarrow}\mathcal{C}_{\mathbb{A}_R}
	\end{equation}
	of \Cref{lem:Phi is well-defined}. 
\end{theorem}

To prove \Cref{thm:eq is orbifolding}, it remains to be shown that $\eta$ and $\epsilon$ are indeed isomorphisms and that everything is well-defined. 
Invertibility of $\eta$ is obvious, invertibility of $\epsilon$ follows from invertibility of~$\gamma$, and the inverse is given by the mirror image of \eqref{eq:counit in adjunction} along the horizontal line, replacing~$\gamma$ with its inverse and the module structure with the induced comodule structure.
We shall not add to our previous comments or expand on the proof here, as it involves the same techniques that have been used throughout this section. 
Instead, we prove directly that~$F$ is an equivalence in \Cref{sec:proof of Equivalence}.

\begin{remark}
	Given an adjoint equivalence where one functor is equipped with a (braided) strong monoidal structure, there is a canonical (braided) strong monoidal structure on the adjoint functor \cite[Rem.\,2.4.10]{EGNO2015}, such that the natural transformations involved in the adjunction are monoidal. 
	In our case, this equips~$ E$ with a braided monoidal structure, so we omit its explicit description here. 
	In short, we have
	\begin{align}
		 E^2_{Y,Y^\prime}&:= E(\epsilon_Y\otimes_{A_G}\epsilon_{Y^\prime})\circ  E( F^2_{ E(Y), E(Y^\prime)})^{-1}\circ \eta_{ E(Y)\otimes E(Y^\prime)}\,,
		 \\
		 E^0&:= E( F^0)^{-1}\circ \eta_{\mathds{1}}\,.
	\end{align}
\end{remark}

\subsection{Proof of Equivalence}
\label{sec:proof of Equivalence}

In this section we prove that the functor $F\colon (\mathcal{C}^\times_G)^G \to \mathcal{C}_{\mathbb{A}_R}$ of \Cref{lem:Phi is well-defined} is an equivalence. 
In \Cref{subsec:fullyfaithful} we argue that~$F$ is fully faithful, and in \Cref{subsec:surjective} we use a result of \cite{Mulevicius2022} to show that it is essentially surjective.

\subsubsection[$ F$ is Fully Faithful]{$\boldsymbol{ F}$ is Fully Faithful}
\label{subsec:fullyfaithful}

We first show that~$F$ is full. 
Then we use pivotality of~$F$ together with a result of \cite{Mulevicius2022a} to find that~$F$ is also faithful. 

\begin{lemma}[label=lem:identification of relative tensor products of morphisms]
	For all $g,h,k \in G$ and all $(X,u) \in (\Gcbc)^G$, we have
	\begin{equation}
		\gamma_{ F(X,u)}^{g,h,k}=\id_{X_{gh^{-1}}}
	\end{equation}
	under the identification $ F(X)_{g,h}=X_{gh^{-1}}= F(X)_{gk,hk}$. 
\end{lemma}

\begin{proof}
	Inserting the relevant ingredients \eqref{eq:FonX}--\eqref{eq:tau structure of monoidal unit in CAR} of $F(X,u)$ into the expression~\eqref{eq:component def of gamma} with $Y=F(X,u)$, one obtains an explicit string diagrammatic presentation for $\gamma_{ F(X,u)}^{g,h,k}$. 
	We claim that it is equal to $\id_{X_{gh^{-1}}}$. 
	This can be seen by explicitly including section and retraction for the relative tensor products in domain and codomain, which for $m_g\otimes_{A_g}m_g^*$ are evaluation and coevaluation, where the latter comes with an insertion of~$\psi^2$. 
	We are then left with a number of loops which cancel their $\psi^2$-insertions, leaving only the identity on~$X_{gh^{-1}}$.
\end{proof}

As a consequence of \Cref{lem:identification of relative tensor products of morphisms}, for any morphism $f\colon  F(X,u)\longrightarrow F(X^\prime,u^\prime)$ in $(\Gcbc)^G$, naturality of $\gamma^{g,e,h}$, i.e.\ $\gamma_{ F(X,u)}^{g,e,h}\circ f_{gh,h}=f_{g,e}\circ \gamma_{ F(X^\prime,u^\prime)}^{g,e,h}$, implies
\begin{equation}
	\label{eq:identification of relative tensor products of morphisms}
	f_{gh,h}=f_{g,e}
\end{equation}
for the components of~$f$ introduced in~\eqref{eq:relative tensor products of morphisms in CAR}. 
We claim that $f_{G,e} = \bigoplus_{g\in G}f_{g,e}$ is a preimage of~$f$ under~$F$. 
The is straightforward to check:
\begin{align*}
	 F(f_{G,e})&\stackrel{\hphantom{\eqref{eq:identification of relative tensor products of morphisms}}}{=}\bigoplus_{g,h\in G} m^*_{gh}\otimes f_{g,e}\otimes m_h\\
	&\stackrel{\eqref{eq:identification of relative tensor products of morphisms}}{=} \bigoplus_{g,h\in G} m^*_{gh}\otimes f_{gh,h}\otimes m_h\\
	&\stackrel{\eqref{eq:relative tensor products of morphisms in CAR}}{=}\bigoplus_{g,h\in G} m^*_{gh}\otimes (m_{gh}\otimes_{A_{gh}} f\otimes_{A_h}m_h^*)\otimes m_h\\
	&\stackrel{\hphantom{\eqref{eq:identification of relative tensor products of morphisms}}}{=}\bigoplus_{g,h\in G} A_{gh}\otimes_{A_{gh}} f\otimes_{A_h}A_h\\
	&\stackrel{\hphantom{\eqref{eq:identification of relative tensor products of morphisms}}}{=} f \, .
\end{align*}
Thus~$F$ is full. 
Since it is also pivotal by \Cref{lem:Phi is well-defined}, it follows from \cite[Prop.\,2.3]{Mulevicius2022a} that~$F$ is fully faithful.

\subsubsection[$F$ is Essentially Surjective]{$\boldsymbol{F}$ is Essentially Surjective}
\label{subsec:surjective}

We prove that~$F$ is essentially surjective by showing that the induced map between isomorphism classes of simple objects is surjective. 
We do so by showing that the images of simple objects saturate the global dimension, i.e.\ adding the squares of their dimensions gives the global dimension of $\CAR$.

\medskip

Recall from \cite[Thm.\,3.17]{Mulevicius2022} the formula for the global dimension of the category associated to an orbifold datum:
\begin{equation}
	\dim(\CAR)=\frac{\dim(\mathcal{C})}{\phi^8(\tr_{\mathcal{C}}(\omega_A^2))^2}\,,
\end{equation}
where we continue to denote $\mathcal{C}\equiv\mathcal{C}_e$. 
To calculate the denominator, we note that
\begin{equation}
	\tr_{\mathcal{C}}(\omega_A^2)=\bigoplus_{g\in G}\,
	\begin{tikzpicture}[very thick,scale=0.6,color=red!50!black, baseline]
		\draw[directedred] 
		(2,0) .. controls +(0,-1) and +(0,-1) .. (1,0) node[pos=0.8](psi){};
		\draw[directedred] 
		(1,0) .. controls +(0,1) and +(0,1) .. (2,0) node[pos=1](g){};
		\draw[directedred] 
		(0,0) .. controls +(0,-2) and +(0,-2) .. (3,0) node[pos=0.2](psi2){};
		\draw[redirectedred] 
		(0,0) .. controls +(0,2) and +(0,2) .. (3,0)node[pos=1](g2){};
		\fill(psi)  circle (2.9pt)
		(psi2) circle (2.9pt);
		\draw ($(psi)+(0.4,0.15)$) node {{\tiny$\psi^2$}}
		($(psi2)+(0.4,0.15)$) node {{\tiny$\psi^2$}}
		($(g)+(0.3,0)$) node{{\scriptsize$g$}}
		($(g2)+(0.3,0)$) node{{\scriptsize$g$}};
	\end{tikzpicture}	\!
	=\bigoplus_{g\in G}\id_\mathds{1}= |G|
\end{equation}
because of~\eqref{eq:PsiMaps} and~\eqref{eq:psiInAR}, and we have $\phi^{-8}=|G|^4$ by~\eqref{eq:phiInAR}. 
Thus 
\begin{equation}
	\label{eq:dimCAR}
	\dim(\CAR)=|G|^2\dim(\mathcal{C}) \, . 
\end{equation}

For $(\Gcbc)^G$, recall from \cite[Eq.\,(6.7)]{Mulevicius2022a} that for a commutative $\Delta$-separable symmetric haploid Frobenius algebra~$B$ we have
\begin{equation}
	\dim \Big(B\text{-}\Mod^\textrm{loc}\big((\Gcbc)^G\big)\Big)= \frac{\dim\big((\Gcbc)^G\big)}{\big(\!\dim_{(\Gcbc)^G}(B)\big)^2}\,.
\end{equation}
On the other hand, by \cite[Prop.\,2.23]{BB2017}, pivotal functors preserve dimensions of objects. Since the functor~$I\colon\Rep(G)\longrightarrow (\Gcbc)^G$ (cf.~\eqref{def:inclusion of RepG}, with~$\mathcal C$ replaced by $\Gcbc$) is pivotal, the dimension of $B=I(\C(G))$ is
\begin{equation}
	\label{eq:dimBIsG}
	\dim_{(\Gcbc)^G}(B)=\dim_{\Rep(G)}(\C(G))=|G|\,.
\end{equation}
As a result, we find that the global dimensions of $\CAR$ and $(\Gcbc)^G$ are the same, 
\begin{align}
	\dim(\CAR)&\stackrel{\eqref{eq:dimCAR}}{=}|G|^2\dim(\mathcal{C}) 
	\nonumber 
	\\
	&\stackrel{\hphantom{(}\ref{thm:inversion of 0form by 1form}\hphantom{)}}{=}|G|^2\dim\big(\big((\Gcbc)^G\big)_{\mathbb{A}_B}\big) 
	\nonumber 
	\\
	&\stackrel{\eqref{eq:orbifold category of a condensable algebra}}{=}|G|^2\frac{\dim\big((\Gcbc)^G\big)}{\big(\!\dim_{(\Gcbc)^G}(B)\big)^2}
	\nonumber 
	\\
	&\stackrel{\eqref{eq:dimBIsG}}{=}\dim\big((\Gcbc)^G\big)\,.
\end{align}

By \cite[Prop.\,6.3.1]{EGNO2015}, a fully faithful functor between fusion categories sends simple objects to simple ones. 
Combining this with the fact that~$F$ preserves dimensions of objects (because~$F$ is pivotal), denoting the set of chosen representatives of isomorphism classes of simple objects in a category~$\mathcal{B}$ by $\mathcal{I}_\mathcal{B}$, and abbreviating $F(X) \equiv F(X,u)$, we have
\begin{equation}
\label{eq:dimension of Geq}
	\sum_{X\in\mathcal{I}_{(\mathcal{C}^{\scaleto{\times}{3pt}}_{\scaleto{G}{3pt}})^{\scaleto{G}{3pt}}}}\dim_{\mathcal{C}_{\mathbb{A}_{\scaleto{R}{3pt}}}}\big( F(X)\big)^2=\sum_{X\in\mathcal{I}_{(\mathcal{C}^{\scaleto{\times}{3pt}}_{\scaleto{G}{3pt}})^{\scaleto{G}{3pt}}}}\dim_{(\Gcbc)^G}(X)^2=\dim \big((\Gcbc)^G\big)=\dim(\CAR)\,.
\end{equation}
Note that since~$ F$ is fully faithful, non-isomorphic simple objects have non-isomorphic images. Therefore we have 
\begin{equation}
\label{eq:dimension of CAR dissected}
	\dim(\CAR)=\sum_{X\in\mathcal{I}_{(\mathcal{C}^{\scaleto{\times}{3pt}}_{\scaleto{G}{3pt}})^{\scaleto{G}{3pt}}}}\dim_{\mathcal{C}_{\mathbb{A}_{\scaleto{R}{3pt}}}}\big( F(X)\big)^2 +\sum_{Y\in \mathcal{I}_{\mathcal{C}_{\mathbb{A}_{\scaleto{R}{3pt}}}}\setminus F(\mathcal{I}_{(\mathcal{C}^{\scaleto{\times}{3pt}}_{\scaleto{G}{3pt}})^{\scaleto{G}{3pt}}})}\dim_{\mathcal{C}_{\mathbb{A}_{\scaleto{R}{3pt}}}}(Y)^2 \,.
\end{equation}
By \cite[Thm.\,7.21.12]{EGNO2015}, dimensions of simple objects in fusion categories have strictly positive squares. Comparing \eqref{eq:dimension of Geq} and \eqref{eq:dimension of CAR dissected},  the second sum in \eqref{eq:dimension of CAR dissected} must therefore be zero. 
This implies that there is no isomorphism class of simple objects in $\CAR$ which does not have a representative of the form $ F(X)$. 
It follows that the map between sets of isomorphism classes of simple objects induced by~$ F$ is surjective, and hence~$ F$ is essentially surjective.

\medskip

In summary, the results of \Cref{subsec:fullyfaithful,subsec:surjective} show that~$F$ is an equivalence, thus completing the proof of \Cref{thm:eq is orbifolding}.

\nocite{*}
\printbibliography
\end{document}